\documentclass[11pt,a4paper]{article}
\usepackage[margin=1in]{geometry}
\usepackage[T1]{fontenc}
\usepackage{lmodern}
\usepackage{amsmath,amssymb,amsthm,mathtools}
\usepackage{microtype}
\usepackage{enumitem}
\usepackage{graphicx}
\usepackage[hidelinks]{hyperref}
\setlist[enumerate]{itemsep=2pt,topsep=4pt}
\setlist[itemize]{itemsep=2pt,topsep=4pt}
\newtheorem{theorem}{Theorem}[section]
\newtheorem{lemma}[theorem]{Lemma}
\newtheorem{proposition}[theorem]{Proposition}
\newtheorem{corollary}[theorem]{Corollary}
\newtheorem{claim}[theorem]{Claim}
\theoremstyle{definition}

\newtheorem{remark}[theorem]{Remark}

\title{Exact local spectral thresholds for loose and tight Hamilton cycles in $3$-graphs}
\author{Pei Liu\thanks{Department of Mathematics, Sungkyunkwan University,
Suwon, 16419, Republic of Korea. E-mail: \href{mailto:liupei2023@g.skku.edu}{\texttt{liupei2023@g.skku.edu}}. Corresponding author. }}
\date{\today}
\begin{document}
\maketitle
\begin{abstract}
For an $n$-vertex $3$-uniform hypergraph $H$, let $\sigma(H)$ denote
the minimum adjacency spectral radius of its vertex links.
We determine the exact thresholds for $\sigma(H)$ that force loose
and tight Hamilton cycles for all sufficiently large admissible orders,
and give constructions showing that each bound is best possible.
The loose threshold is $(1+\sqrt{13})n/8+O(1)$, with two forms
according to $n$ modulo four.
The tight threshold is $2n/3+O(1)$, with four forms determined by
$n$ modulo three and, when $3\mid n$, modulo nine.
The resulting spectral conditions and the corresponding minimum
vertex degree conditions do not imply one another.
For loose cycles, we combine fractional packing stability and absorption
with exact completion near a space barrier.
For tight cycles, spectral stability and rooted sampling reduce the
problem to three extremal configurations.
We resolve these configurations by controlling the spectral loss caused
by missing triples and constructing short paths that extend to Hamilton
cycles, including in the presence of exceptional vertices.
\end{abstract}
\section{Introduction}
\label{sec:introduction}
A $3$-uniform hypergraph, or $3$-graph, is a pair $H=(V(H),E(H))$ with $E(H)\subseteq\binom{V(H)}3$. For $v\in V(H)$, its link $L_H(v)$ is the graph on $V(H)\setminus\{v\}$ whose edges are the pairs $xy$ with $vxy\in E(H)$. We write $\rho(G)$ for the adjacency spectral radius of a graph $G$ and put $\sigma(H)=\min_{v\in V(H)}\rho(L_H(v))$. A loose cycle has distinct vertices $a_1,b_1,\ldots,a_t,b_t$ and edges $a_ib_ia_{i+1}$, where $t\ge3$ and indices are cyclic. A tight path is an ordering of distinct vertices in which every three consecutive vertices form an edge; a tight cycle has a cyclic ordering $v_1,\ldots,v_t$ of distinct vertices, where $t\ge4$ and every $v_iv_{i+1}v_{i+2}$ is an edge, with indices taken modulo $t$. A cycle is Hamilton if it covers all vertices; a loose Hamilton cycle therefore requires even order.

Write $d_H(v)=|E(L_H(v))|$, $\delta_1(H)=\min_v d_H(v)$
and $\Delta_1(H)=\max_v d_H(v)$.
The minimum codegree $\delta_2(H)$ is the minimum number
of edges containing a pair of vertices.
Minimum degree conditions for Hamilton cycles are central
in extremal hypergraph theory; see Zhao~\cite{ZhaoSurvey}
for a survey.
For $3$-graphs of sufficiently large order, R\"odl,
Ruci\'nski and Szemer\'edi~\cite{RRSExact} determined the
exact minimum codegree threshold $\lfloor n/2\rfloor$
for tight Hamilton cycles.
For loose Hamilton cycles, K\"uhn and
Osthus~\cite{KuhnOsthusLoose} proved the asymptotic codegree
threshold $n/4$, and Czygrinow and Molla~\cite{CzygrinowMolla}
established the exact threshold $\lceil n/4\rceil$
for sufficiently large even $n$.
Further exact codegree thresholds for Hamilton $\ell$-cycles
in uniform hypergraphs were determined
in~\cite{GanHanXu,HanZhaoCodegree,LuYuanHamilton}.

For minimum vertex degree, Bu\ss, H\`an and
Schacht~\cite{BussHanSchacht} proved that the asymptotic
threshold for loose Hamilton cycles in $3$-graphs is
$\frac7{16}\binom n2$.
Han and Zhao~\cite{HanZhao} subsequently determined the
exact threshold for all sufficiently large even orders.
For tight Hamilton cycles, Reiher, R\"odl, Ruci\'nski,
Schacht and Szemer\'edi~\cite{ReiherEtAl} proved the
asymptotically sharp threshold $\frac59\binom n2$.
Space barriers and constructions excluding one edge type
across a partition play an essential role in these problems;
see also Han and Zhao~\cite{HanZhaoForbidding}.
For extensions to higher uniformities and more general Hamilton
cycles, see~\cite{LangSanhuezaDegree,LangSanhuezaLoose,
LangSchachtVolec,PolcynEtAl}.

Other spectral Hamiltonicity results concern a different
spectrum or impose additional distribution assumptions.
Brooks, Linz and Luo~\cite{BrooksLinzLuo} use the
adjacency-tensor spectral radius to force Hamiltonian
Berge cycles.
Lenz, Mubayi and Mycroft~\cite{LenzMubayiMycroft} proved
loose Hamiltonicity in dense quasirandom hypergraphs
with a minimum vertex degree condition.
H\`an, Han and Morris~\cite{HanHanMorris} extended this
line of work to sparse pseudorandom hypergraphs, with
consequences under spectral-gap assumptions. For tight cycles, Ara\'ujo, Piga and Schacht~\cite{AraujoPigaSchacht}
obtained a sharp density condition for uniformly dense $3$-graphs
with minimum vertex degree of order $n^2$. 
In our results, the only density assumption is a lower
bound on the largest adjacency eigenvalue of every
vertex link.

The degree of a vertex counts the edges in its link,
whereas the spectral radius also depends on their arrangement.
This leads to the question of whether a lower bound on
$\sigma(H)$ can force a Hamilton cycle.
For perfect matchings, Lin, Lu, Yuan and
Zhao~\cite{LinLuYuanZhao} proved that, for fixed $\gamma>0$
and sufficiently large $n$ divisible by three,
$\sigma(H)>(2/3+\gamma)n$ suffices, and that the leading
constant is best possible. 
Lu and Yuan~\cite{LuYuanPartite} proved an asymptotically sharp local
spectral condition for perfect matchings in balanced $3$-partite
$3$-graphs. 
Liu and O~\cite{LiuO} determined the exact matching
threshold $2n/3-2$ for sufficiently large orders divisible
by three and the corresponding thresholds in balanced
$3$-partite $3$-graphs.
Here we determine the exact local spectral thresholds for loose
and tight Hamilton cycles for all sufficiently large admissible
orders. The tight threshold has leading term $2n/3$, as in the
matching problem, but its exact value also depends on $n$ modulo
nine. This dependence comes from a triangle decomposition in
the internal link of an extremal construction; see
Remark~\ref{rem:mod-nine}. Moreover, the spectral hypotheses and
the corresponding asymptotic minimum vertex degree hypotheses
do not imply one another, even with fixed positive margins;
see Remark~\ref{rem:incomparable}. Thus the degree theorems do
not yield the spectral results by an implication between
their hypotheses.

Our fractional-cover normalization and stability arguments for
split graphs and cliques adapt the matching methods of
Liu and O~\cite[Theorem~3.4 and Lemmas~3.5--3.6]{LiuO}.
For loose cycles, we use oriented-edge packings with load vector
$(2,1,1)$ to reflect the overlap of consecutive edges.
For tight cycles, we use the Hamilton framework of Lang and
Sanhueza-Matamala~\cite{LangSanhuezaDegree} in the form supplied
by their bandwidth theorem~\cite{LangSanhueza}. This requires
tight connectivity, aperiodicity and consistency under vertex
deletion. The exact argument then has to construct paths that
correct the part sizes or join distinct tight components,
while incorporating exceptional vertices. These steps account
for the different completion arguments near the three barriers.

Write $I_s$ for an independent graph of order $s$, $P_3$ for a path on three vertices, and $\vee$ for graph join. For even $n$, set $a_n=\lceil n/4\rceil-1$ and
\begin{equation}\label{eq:loose-threshold}
\tau_n=\begin{cases}
\rho\bigl(K_{a_n}\vee(K_2\cup I_{n-a_n-3})\bigr),&n\equiv0\pmod4,\\
\rho\bigl(K_{a_n}\vee I_{n-a_n-1}\bigr),&n\equiv2\pmod4.
\end{cases}
\end{equation}
\begin{theorem}\label{thm:loose}
For every sufficiently large even integer $n$, an $n$-vertex $3$-graph $H$ with $\sigma(H)>\tau_n$ contains a loose Hamilton cycle. The bound is best possible for each such $n$.
\end{theorem}

For the tight threshold, put $W_m=\rho(K_{m-1}\vee(P_3\cup I_{2m-3}))$, $U_m=\rho(K_1\vee(K_m\cup K_{2m-2}))$, and $V_m=\rho(K_m\vee I_{2m})=(m-1+\sqrt{9m^2-2m+1})/2$. Define
\begin{equation}\label{eq:tight-threshold}
T_n=\begin{cases}
W_m,&n=3m,\quad 3\mid m,\\
U_m,&n=3m,\quad 3\nmid m,\\
V_m,&n=3m+1,\\
2m,&n=3m+2.
\end{cases}
\end{equation}
\begin{theorem}\label{thm:tight}
For every sufficiently large integer $n$, an $n$-vertex $3$-graph $H$ with $\sigma(H)>T_n$ contains a tight Hamilton cycle. The bound is best possible for each such $n$.
\end{theorem}
Thus $\tau_n$ and $T_n$ are the maximum values of $\sigma(H)$ among $3$-graphs without the corresponding Hamilton cycle. They satisfy $\tau_n=cn+O(1)$ and $T_n=2n/3+O(1)$, where $c=(1+\sqrt{13})/8$. The hypergraphs attaining these values need not be unique;
see Remark~\ref{rem:nonunique}.

\begin{proposition}\label{prop:sharp}
For every sufficiently large even $n$, there is an $n$-vertex $3$-graph without a loose Hamilton cycle and with minimum link radius $\tau_n$. For every sufficiently large $n$, there is an $n$-vertex $3$-graph without a tight Hamilton cycle and with minimum link radius $T_n$.
\end{proposition}
\begin{proof}
For a set $A$ of size $a$, the space-barrier construction
$H_A$ consists of all triples meeting $A$. Its $A$-links are complete and its other links are $K_a\vee I_{n-a-1}$, so
\begin{equation}\label{eq:space-radius}
\sigma(H_A)=\frac{a-1+\sqrt{(a-1)^2+4a(n-a-1)}}2.
\end{equation}
A loose Hamilton cycle has $n/2$ edges and uses each vertex at most twice, requiring $n/2\le2a$. A tight Hamilton cycle has $n$ edges and uses each vertex three times, requiring $n\le3a$. Taking $a=a_n$ when $n\equiv2\pmod4$ gives the second loose branch; taking $a=m$ when $n=3m+1$ gives the third tight branch.

When $4\mid n$, take $|A|=a_n$, fix $p,q\notin A$, and add to $H_A$ all triples containing $p,q$ and lying outside $A$. A loose Hamilton cycle would need at least $n/2-2a_n=2$ edges outside $A$, but two such edges share two vertices. At every vertex outside $A\cup\{p,q\}$ the link is $K_{a_n}\vee(K_2\cup I_{n-a_n-3})$, and the other links contain this graph. Hence $\sigma(H)=\tau_n$. In both loose constructions, the minimum vertex degree is $\binom{n-1}{2}-\binom{\lfloor3n/4\rfloor}{2}+\mathbf1_{\{4\mid n\}}$, one below the forcing threshold in~\cite{HanZhao}. Thus they also give the sharp degree examples in~\cite[Proposition 1.3]{HanZhao}.

For $n=3m$ with $3\mid m$, take $|A|=m-1$ and $B=V\setminus A$. Fix $z\in B$, partition $B\setminus\{z\}$ into triples, and add to $H_A$ each edge $zxy$ for which $x,y$ lie in the same triple. A tight Hamilton cycle would need at least $3m-3(m-1)=3$ edges inside $B$. All contain $z$, so all three cycle edges through $z$ would be internal. Their five consecutive vertices give a four-vertex path in the internal link of $z$, which is a disjoint union of triangles. This is impossible. For $b\in B\setminus\{z\}$ the link is $K_{m-1}\vee(P_3\cup I_{2m-3})$; the link at $z$ contains this graph, and the $A$-links are complete. Thus $\sigma(H)=W_m$.

For $n=3m$, partition the vertices into $\{h\},A,B$ of sizes $1,m,2m-1$. Include all triples through $h$ and precisely the types $AAA,AAB,BBB$ on $A\cup B$. Deleting $h$ from a tight Hamilton cycle would leave a spanning tight path. Since any two consecutive edges of a tight path share exactly two vertices, its edges cannot move between the $BBB$ component and the component of types $AAA,AAB$. Thus they must all contain at least two $A$ vertices. Consecutive $B$ vertices in its ordering would then be at distance at least three, allowing at most $\lceil\frac{3m-1}{3}\rceil=m$ of them, a contradiction. The $B$-links are $K_1\vee(K_m\cup K_{2m-2})$, the $A$-links are $K_m\vee I_{2m-1}$ of radius $2m-1$, and the $h$-link is complete. Lemma~\ref{old:2.1} gives $U_m<2m-1$ and $W_m>U_m$. This construction gives the second tight branch.

Finally, for $n=3m+2$, take $|A|=2m+1$, $|B|=m+1$, and all triples with at least two $A$ vertices, together with all $BBB$ triples. The latter again form a separate tight component. A tight Hamilton cycle would have at least two $A$ vertices in every edge, requiring $3|A|\ge2n$, which fails. A $B$-link is $K_{2m+1}\cup K_m$ of radius $2m$, whereas an $A$-link is $K_{2m}\vee I_{m+1}$ of larger radius. This proves the last branch.
\end{proof}

\begin{remark}\label{rem:nonunique}
The extremal hypergraphs need not be unique. For $n=3m+2$,
take disjoint sets $A,B$ with $|A|=2m+1$ and $|B|=m+1$.
Include every triple with at least two vertices in $A$, and
the edges of an arbitrary $3$-graph $F$ on $B$.
For $b\in B$, let $L_F(b)$ be the graph on
$B\setminus\{b\}$ whose edges are the pairs $xy$ with
$bxy\in E(F)$. Then $L_H(b)=K_{2m+1}\cup L_F(b)$,
where the union is disjoint. Since $L_F(b)$ has $m$
vertices, its spectral radius is at most $m-1$, and hence
$\rho(L_H(b))=2m$. For $a\in A$, the link is
$K_{2m}\vee I_{m+1}$ and has larger radius. A tight Hamilton cycle cannot move between the
$BBB$ edges and the other edge types. It would therefore
use at least two vertices of $A$ in every edge, contradicting
$3|A|<2n$. Thus every choice of $F$ gives
$\sigma(H)=T_n$, and choices with different numbers of edges
give nonisomorphic extremizers.

Nonuniqueness also occurs in both loose constructions of
Proposition~\ref{prop:sharp}. Delete one triple inside $A$.
The minimum links outside $A$ are unchanged, while every
affected link is a complete graph on $n-1$ vertices with one
edge deleted. Its radius is at least
$n-2-2/(n-1)>\tau_n$ for sufficiently large $n$, by the
Rayleigh quotient of the constant unit vector. Deletion
cannot create a Hamilton cycle, so the resulting hypergraph
is another extremizer. We do not classify all hypergraphs
attaining either threshold.
\end{remark}

\begin{remark}\label{rem:mod-nine}
The dependence on $n$ modulo nine reflects the internal link
in the construction attaining $W_m$ in Proposition~\ref{prop:sharp}.
Write $F=L_{H[B]}(z)$. Minimum degree at least two in $F$
ensures that the internal link of each vertex in
$B\setminus\{z\}$ contains a copy of $P_3$, whereas the absence
of a four-vertex path in $F$ obstructs a tight Hamilton cycle.
A graph with these two properties is necessarily a disjoint
union of triangles. Since $F$ has $2m$ vertices, these
requirements force $3\mid m$, or equivalently $9\mid n$.
\end{remark}

\begin{remark}\label{rem:incomparable}
The local spectral condition and the corresponding minimum vertex degree condition do not imply each other, even with a fixed asymptotic margin. Write $(c_0,d_0)=((1+\sqrt{13})/8,7/16)$ for loose cycles and $(c_0,d_0)=(2/3,5/9)$ for tight cycles. In both cases $d_0<c_0<\sqrt{d_0}$. Here we compare $\sigma(H)\ge(c_0+\gamma)n$ with $\delta_1(H)\ge(d_0+\gamma)\binom n2$, where $\gamma>0$ is fixed and sufficiently small. Choose $c_0<q<\sqrt{d_0}$, take a complete $3$-graph on $n-1$ vertices, and add a vertex whose link is $K_{\lfloor qn\rfloor}\cup I_{n-1-\lfloor qn\rfloor}$. Then $\sigma(H)=\lfloor qn\rfloor-1$, while $\delta_1(H)=\binom{\lfloor qn\rfloor}{2}$, so the spectral condition holds and the degree condition fails. Conversely, choose $d_0<p<c_0$, take $V(H)=\mathbb Z_n$, and fix $S\subseteq\mathbb Z_n$ of size $\lfloor pn\rfloor$. Include each triple of distinct vertices whose sum belongs to $S$. Every pair has between $|S|-2$ and $|S|$ completions. The row-sum bounds on each link give $\sigma(H)=pn+O(1)$ and $\delta_1(H)=p\binom{n-1}{2}+O(n)$, so these hypergraphs satisfy the degree condition and fail the spectral condition. Choose $\gamma$ smaller than $q-c_0$, $d_0-q^2$, $p-d_0$ and $c_0-p$, and then take $n$ sufficiently large. Restrict to even $n$ for loose cycles.
\end{remark}

For loose cycles, fractional packing stability and absorption
handle the nonextremal case, followed by exact completion near
the space barrier. For tight cycles, spectral stability and
rooted sampling leave three extremal configurations. Their
completion arguments share a final assignment step, proved
once in Lemma~\ref{lem:common-completion}. The preceding path
constructions depend on the barrier: they correct a part-size
imbalance or provide crossings between tight components.
Since competing link radii can differ by order $1/n$, these
constructions require bounds on the spectral loss caused by
missing triples that retain this precision.

Section~\ref{sec:tools} collects spectral and rooted-sampling estimates,
and Section~\ref{sec:loose} proves Theorem~\ref{thm:loose}.
Section~\ref{sec:stability} reduces the tight-cycle problem to three
extremal configurations. Section~\ref{sec:space-twohit} resolves the
space and two-hit barriers, while Section~\ref{sec:smallpart} treats
the remaining missing-layer barrier and completes the proof of
Theorem~\ref{thm:tight}. All graphs and hypergraphs are finite and simple.
The notation $x\ll y$ means that $x$ is chosen sufficiently small
in terms of $y$.

\section{Spectral estimates and rooted sampling}\label{sec:tools}
Matrix norms are operator norms unless the Frobenius norm $\|\cdot\|_F$ is specified. Put $S_m=\rho(K_{m-1}\vee(K_2\cup I_{2m-2}))$ and $R_m=\rho(K_{m-1}\vee(2K_2\cup I_{2m-4}))$.

We will use the bound $\rho(G)(\rho(G)+1)\le2e(G)$. To see this, choose a vertex $v$ attaining the maximum coordinate of a nonnegative Perron vector. The eigenvalue equations give $\rho(G)\le d(v)$ and $\rho(G)^2\le\sum_{u\in N(v)}d(u)\le2e(G)-d(v)$. We will also use
$
\|A(G)-A(G')\|
\le \|A(G)-A(G')\|_F
=\sqrt{2|E(G)\triangle E(G')|}.
$

\begin{lemma}\label{old:2.1}
For every integer $m\ge3$,
\[
2m-2<S_m<U_m<R_m<W_m,\qquad U_m<2m-1.
\]
Moreover, $U_m>2m-2+1/(3m)$ and $2m-2/3<V_m<2m$.
\end{lemma}
\begin{proof}
Put $k=m-1$. Equitable partitions give quotient matrices with characteristic polynomials
\begin{align*}
s(t)&=t^3-kt^2-(2k^2+k+1)t+2k^2,\\
q(t)&=t^3-(3k-1)t^2+(2k^2-4k-1)t+4k^2+k-1
\end{align*}
for $S_m,U_m$, respectively, while $r(t)=s(t)-2k$ is the corresponding polynomial for $R_m$. The translates
\[
s(2k+x)=x^3+5kx^2+(6k^2-k-1)x-2k,\quad
q(2k+x)=x^3+(3k+1)x^2+(2k^2-1)x-k-1
\]
are strictly increasing for $x\ge0$ and negative at $x=0$, so $s,q$ have unique roots above $2k$, namely $S_m,U_m$. For $k\ge2$, substitution gives $s(2k)<0<s(2k+2/(5k))$ and $q(2k+2/(5k))<0$, proving $2k<S_m<U_m$. Also
\[
q\left(2k+\frac1{k+1}\right)=\frac{k(k^3-2k-2)}{(k+1)^3}>0,\qquad
q\left(2k+\frac2{3k}\right)=\frac k3-1+\frac2{3k}+\frac4{9k^2}+\frac8{27k^3}>0.
\]
The second expression is $4/27$ at $k=2$ and is positive for $k\ge3$. Thus $U_m<2k+1/(k+1)<2k+1$. Writing $U_m=2k+x$, we also have $0<x<2/(3k)$, and $q(U_m)=0$ gives
\[
r(U_m)=1-3k+(4k^2-k)x+(2k-1)x^2
< -\frac{(k-2)(3k^2+3k-2)}{9k^2}\le0.
\]
Hence $U_m<R_m$. Replacing disjoint edges $ab,cd$ in the graph for $R_m$ by $ab,bc$ preserves the Rayleigh quotient of its positive Perron vector, since the four coordinates agree, but that vector ceases to be an eigenvector. Therefore $W_m>R_m$. Finally, $q(2m-2+1/(3m))<0$, and the quadratic for $V_m$ has values $-2/9$ and $2m$ at $2m-2/3$ and $2m$, respectively.
\end{proof}

\begin{lemma}\label{old:2.2}
Let $A,E$ be real symmetric matrices. Suppose $A$ has largest eigenvalue $\lambda_1$, unit eigenvector $u$, and second eigenvalue $\lambda_2$, with $g=\lambda_1-\lambda_2>2\|E\|$. Then
\[
\lambda_1+u^TEu\le\lambda_{\max}(A+E)
\le\lambda_1+u^TEu+\frac{\|Eu\|^2}{g-2\|E\|}.
\]
\end{lemma}
\begin{proof}
The lower bound follows by evaluating the Rayleigh quotient at $u$. For the upper bound, decompose the space as $\mathbb Ru\oplus u^\perp$ and write $A+E=\bigl(\begin{smallmatrix}\lambda_1+\epsilon&b^T\\b&K\end{smallmatrix}\bigr)$, where $\epsilon=u^TEu$, $\|b\|\le\|Eu\|$, and $\lambda_{\max}(K)\le\lambda_2+\|E\|$. The largest eigenvalue $\mu\ge\lambda_1-\|E\|$ of this matrix exceeds $\lambda_{\max}(K)$. We may therefore eliminate the perpendicular coordinates of its eigenvector to obtain $\mu=\lambda_1+\epsilon+b^T(\mu I-K)^{-1}b\le\lambda_1+\epsilon+\|b\|^2/(\mu-\lambda_{\max}(K))$. Since the denominator is at least $g-2\|E\|$, the claimed bound follows.
\end{proof}

\begin{lemma}\label{old:2.3}
Let $0<a<1$, $0<\eta<a$ and $m\ge4/\eta$. If a graph $J$ on $N\ge m$ vertices satisfies $\rho(J)\ge aN$, then a uniform $m$-subset $U$ satisfies
\[
\Pr\bigl(\rho(J[U])<(a-\eta)m\bigr)\le2\exp(-a^4\eta^2m/128).
\]
\end{lemma}
\begin{proof}
Take a nonnegative unit Perron vector $x$ and put $y_i=\sqrt N x_i$. The eigenvalue equation gives $\rho(J)x_i\le\sum_jx_j\le\sqrt N$, so $0\le y_i\le a^{-1}$ and $\sum_i y_i^2=N$. For $A=A(J)$, put $D=\sum_{i\in U}y_i^2$ and $W=\sum_{i,j\in U}A_{ij}y_iy_j$. Then $\mathbb ED=\frac{m}{N}\sum_i y_i^{2}=m$ and $\mathbb EW=m(m-1)\rho(J)/(N-1)\ge am(m-1)$. Replacing one sampled vertex changes $D,W$ by at most $a^{-2},2ma^{-2}$, respectively. The same bounds hold for the exposure-martingale differences: given two possible next choices, couple the remaining samples by exchanging those two vertices.

Azuma's inequality~\cite{Azuma} gives $\Pr(D>(1+\eta/2)m)\le e^{-a^4\eta^2m/8}$ and $\Pr(W<(a-\eta/2)m^2)\le e^{-a^4\eta^2m/128}$. For the latter, the distance below the expectation is at least $\eta m^2/4$, since $m\ge4/\eta$. Outside these events, $W,D>0$, and the Rayleigh quotient of $y|_U$ gives $\rho(J[U])\ge W/D\ge(a-\eta/2)m/(1+\eta/2)\ge(a-\eta)m$.
\end{proof}

Rooted sampling means choosing a uniform vertex set conditioned on containing a prescribed set of vertices.

\begin{lemma}[Rooted Sampling]\label{old:2.4}
Let $0<b<b+\gamma<1$. For all sufficiently large $s$ and then all sufficiently large $n$, if $\sigma(H)\ge(b+\gamma)n$, every six-set $Q\subseteq V(H)$ is contained in at least $(1-s^{-2})\binom{n-6}{s-6}$ sets $S$ of size $s$ with $\sigma(H[S])\ge(b+\gamma/2)s$. The failure proportion can also be taken to be $1/(2s^2)$.
\end{lemma}
\begin{proof}
Fix a six-set $Q$, and put $a=b+\gamma$, $\eta=\gamma/4$ and
$\kappa=a^4\eta^2/128$. Choose $s$ sufficiently large that
$s-7\ge4/\eta$, $(a-\eta)(s-7)\ge(b+\gamma/2)s$, and
$4se^{-\kappa(s-7)}<1/(2s^2)$. Let $n\ge12s+1$, and choose $S$
uniformly among the $s$-sets containing $Q$.

Conditional on $v\in S$, the set $U=S\setminus(Q\cup\{v\})$ has
size $m=s-6$ if $v\in Q$, and $m=s-7$ otherwise. Its distribution
is that of a uniform $m$-subset of $V(H)\setminus\{v\}$ conditioned
to avoid $Q\setminus\{v\}$. Let $\mathcal A$ denote this avoidance event. By the union bound,
$\Pr(\mathcal A)\ge1-6m/(n-1)\ge1/2$. Thus, for any event
$\mathcal E$ in the unconditioned sampling space,
\[
\Pr(\mathcal E\mid\mathcal A)
=\frac{\Pr(\mathcal E\cap\mathcal A)}{\Pr(\mathcal A)}
\le\frac{\Pr(\mathcal E)}{\Pr(\mathcal A)}
\le2\Pr(\mathcal E).
\]
Since $\rho(L_H(v))\ge an\ge a(n-1)$, Lemma~\ref{old:2.3} bounds
the probability of $\rho(L_H(v)[U])<(a-\eta)m$ by
$2e^{-\kappa m}$ before conditioning, and hence by
$4e^{-\kappa m}$ afterwards.

Since $L_H(v)[U]\subseteq L_{H[S]}(v)$, monotonicity of the
spectral radius and $(a-\eta)m\ge(b+\gamma/2)s$ show that,
conditional on $v\in S$, the probability of
$\rho(L_{H[S]}(v))<(b+\gamma/2)s$ is at most
$4e^{-\kappa(s-7)}$. Let $X$ count the vertices of $S$ violating
this bound. By linearity of expectation,
\[
 \mathbb{E}X
 \le 4e^{-\kappa(s-7)}
       \sum_{v\in V(H)}\Pr(v\in S)
 =4se^{-\kappa(s-7)}
 <\frac{1}{2s^2}.
\]
Here the equality uses $|S|=s$. Markov's inequality gives
$\Pr(X\ge1)<1/(2s^2)$. Consequently, at least a proportion
$1-1/(2s^2)$ of the $s$-sets containing $Q$ satisfy
$\sigma(H[S])\ge(b+\gamma/2)s$, proving both assertions.
\end{proof}
\section{Loose Hamilton cycles}\label{sec:loose}
A loose path has distinct vertices $v_1,\ldots,v_{2r+1}$ and edges $v_{2i-1}v_{2i}v_{2i+1}$, $1\le i\le r$; its ends are $v_1,v_{2r+1}$. Throughout this section, $c=(1+\sqrt{13})/8$.

\subsection{Fractional stability}
An oriented edge has one distinguished vertex and load vector $(2,1,1)$. A fractional packing assigns nonnegative weights to oriented edges so that the total load at each vertex is at most one. It is perfect if every load is one.

\begin{lemma}\label{old:3.1}
There is an absolute constant $C$ such that, for sufficiently small $\delta>0$, if $\sigma(H)\ge(c-\delta)n$ and $H$ has no perfect fractional oriented-edge packing, then $H$ has a vertex cover of size at most $(1/4+C\delta)n$.
\end{lemma}
\begin{proof}
Each oriented edge has total load four, so the maximum total
packing weight is at most $n/4$, with equality precisely for a
perfect packing. The dual program minimizes $\sum_v t_v$ over
nonnegative vertex weights satisfying $2t_u+t_v+t_z\ge1$ for
every oriented edge $uvz$, where $u$ is the distinguished vertex.
Since no perfect packing exists, strong duality gives such
weights with $\sum_v t_v<n/4$. Put $\alpha=\min_v t_v<1/4$ and
$\widetilde z_v=(t_v-\alpha)/(1-4\alpha)$. These numbers have
minimum zero, total less than $n/4$, and satisfy
$2\widetilde z_u+\widetilde z_v+\widetilde z_z\ge1$ for every
oriented edge. Truncating each coordinate at one preserves
these inequalities. Write $z$ for the resulting vector and
choose $u$ with $z_u=0$.

Let $G=L_H(u)$ with $u$ adjoined as an isolated vertex. Taking $u$ as the distinguished vertex in each edge containing it shows that $z$ is a fractional vertex cover of $G$. For $t\in(0,1/2)$, round a coordinate to $0,1/2,1$ according as it is below $t$, between $t$ and $1-t$, or above $1-t$. This preserves every graph-cover inequality, and for uniform $t$ each rounded coordinate has expectation $z_v$. Some rounding therefore has total less than $n/4$. Let $A,Y,Z$ be its weight-one, weight-half and weight-zero classes, of sizes $a,b,n-a-b$. Then
\[
 G\subseteq K_a\vee(K_b\cup I_{n-a-b}),\qquad a+b/2<n/4.
\]
Put $x=a/n$ and $y=b/n$, so $2x+y<1/2$. Adding loops to the
two clique classes gives the symmetric normalized quotient
\[
 Q(x,y)=
 \begin{pmatrix}
 x&\sqrt{xy}&\sqrt{x(1-x-y)}\\
 \sqrt{xy}&y&0\\
 \sqrt{x(1-x-y)}&0&0
 \end{pmatrix},
\]
whose largest eigenvalue $s$ satisfies $s \ge\rho(G)/n$.
For $x,y>0$, this eigenvalue $s>y$ satisfies $f(s,y)=0$, where
$f(s,y)=s-x-xy/(s-y)-x(1-x-y)/s$.
Since $\partial_s f>0$ and
$\partial_y f=x/s-xs/(s-y)^2\le0$, we have $ds/dy=-\partial_y f/\partial_s f \ge 0$, which shows that $s$ is nondecreasing in $y$. The boundary cases
follow by continuity. We may therefore increase $y$ to
$1/2-2x$; write $r_x\ge \rho(G)/n\ge c-\delta$ for the resulting largest
eigenvalue. Its characteristic polynomial is
\[
 F_x(s)=s^3-(1/2-x)s^2-x(1/2+x)s
       +x(1/2-2x)(1/2+x).
\]
Deleting the first row and column of $Q(x,1/2-2x)$ leaves
$\operatorname{diag}(1/2-2x,0)$, so interlacing bounds the
second eigenvalue by $1/2-2x<c$. Moreover,
$F_x(c)=(1/4-x)(2x^2+(c+1)x+c/2-3/16)$, where the second
factor is at least $(\sqrt{13}-2)/16>1/10$.
Thus $F_x(c)\ge0$, and, since the other two eigenvalues are
less than $c$, we obtain $r_x\le c$.
Finally, $F_x(r_x)=0$ and $|F'_x(s)|<2$ for $0\le s\le c$.
The mean value theorem therefore gives
$(1/4-x)/10\le F_x(c)\le2(c-r_x)\le2\delta$.
Hence $a\ge(1/4-20\delta)n\ge n/5$ for sufficiently small
$\delta$, and $b<n/2-2a\le40\delta n$.

Delete the edges of $G[Y]$ to obtain a graph $K\subseteq T=K_a\vee I_{n-a}$. The deleted matrix has norm at most $b$, so $\rho(K)\ge(c-41\delta)n$. Writing $\lambda=\rho(T)\le cn$ and $D=\lambda-\rho(K)$ gives $D\le41\delta n$. Let $p,q$ be nonnegative unit Perron vectors of $T,K$, respectively. All other eigenvalues of $T$ are nonpositive, whence $1-(p^\top q)^2\le D/\lambda$ and
\[
 \|p-q\|^2\le 2D/\lambda\le5D/n,
\]
since $\lambda\ge\sqrt{a(n-a)}\ge2n/5$. Every coordinate of $p$ is at least $1/(5\sqrt n)$: the clique coordinate is at least $1/\sqrt n$, and the independent coordinate is $a/\lambda$ times it. Consequently, only $O(D)$ coordinates satisfy $q_v<p_v/2$. Every missing edge avoiding those vertices contributes at least $1/(50n)$ to $q^\top(A(T)-A(K))q\le D$. Including the edges meeting those vertices, the number $M$ of missing edges is $O(Dn)=O(\delta n^2)$.

Return to the unrounded vector $z$, put $B=V(H)\setminus A$, and write $Z_A=\sum_{v\in A}z_v$, $Z_B=\sum_{v\in B}z_v$ and $S=Z_A+Z_B<n/4$. Summing the cover inequalities over present $A$--$B$ edges yields $(n-a)Z_A+aZ_B\ge a(n-a)-M$. Hence
\[
 (n-2a)Z_B\le(n-a)(S-a)+M=O(\delta n^2),
\]
using $a\ge n/4-20\delta n$. Since $n-2a\ge n/2$, we have $Z_B=O(\delta n)$. The set $A'=A\cup\{v\notin A:z_v\ge1/4\}$ has size at most $a+4Z_B\le(1/4+C\delta)n$. Any edge outside $A'$ would have all three coordinates below $1/4$, contradicting its oriented-edge inequality. Thus $A'$ is a vertex cover.
\end{proof}

\subsection{Absorption and the nonextremal case}
\begin{theorem}\label{old:3.2}
For every $\eta>0$ there are $\varepsilon>0$ and $n_0$ such that the following holds. If $n\ge n_0$ is even, $\sigma(H)\ge(c-\varepsilon)n$, and every $\lfloor3n/4\rfloor$-set $B\subseteq V(H)$ satisfies $e(H[B])>\eta n^3$, then $H$ contains a loose Hamilton cycle.
\end{theorem}

For distinct vertices $x,y$, put $N_H(xy)=\{z:xyz\in E(H)\}$ and $\operatorname{codeg}(x,y)=|N_H(xy)|$. Set $N_H(xx)=\varnothing$ and $\operatorname{codeg}(x,x)=0$. For $U\subseteq V(H)$, write $d_U(xy)=|N_H(xy)\cap U|$, with the analogous notation for other ambient hypergraphs.

\begin{figure}[htbp]
  \centering
  \includegraphics[width=0.85\linewidth]{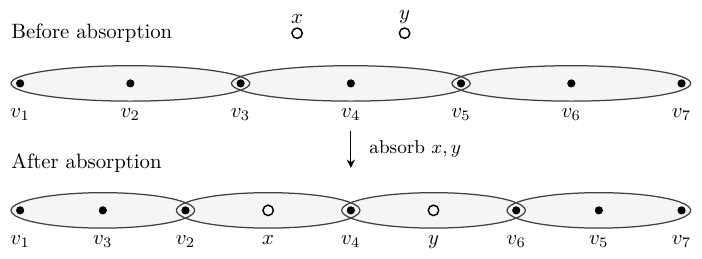}
  \caption{A seven-vertex absorber for $x,y$, before (top)
  and after (bottom) absorption. Each oval represents a
  hyperedge; both paths have ends $v_1,v_7$.}
  \label{fig:absorber}
\end{figure}

\begin{lemma}\label{old:3.3}
For every $\theta\in(4/7,1)$ there are $\kappa>0$ and $n_0$ such that, if $\sigma(H)\ge\theta n$ and $n\ge n_0$, each pair of distinct vertices $x,y$ has at least $\kappa n^7$ ordered absorbing seven-tuples. Such a tuple consists of distinct vertices $v_1,\ldots,v_7$ outside $\{x,y\}$ with edges $v_1v_2v_3$, $v_3v_4v_5$, $v_5v_6v_7$, $v_2xv_4$ and $v_4yv_6$.
\end{lemma}
\begin{proof}
Choose $\eta>0$ so small that, for $h=1-\theta+\eta$, we have $h<\theta$ and $f(h)<\theta$, where $f(h)=(h+\sqrt{h^2+(1-h)^2})/2$. This is possible because $f(1-\theta)<\theta$ is equivalent to $\theta>4/7$. Choose $\tau>0$ with $2\sqrt\tau<\eta$, $\theta-2\sqrt\tau>\max\{f(h),h\}$ and $\sqrt{2\tau}<\eta/2$.

For each $x$, let $S_x=\{v:\operatorname{codeg}(x,v)\ge\tau n\}$. Deleting the other vertices from its link removes at most $\tau n^2$ edges. The Frobenius bound gives $|S_x|\ge \rho(L_{H}(x)[S_x]) \ge (\theta-\sqrt{2\tau})n$. For distinct vertices $x,y$, let $D(x,y)$ be the set of all
$v\in S_x\cap S_y$ such that
$\max\{\operatorname{codeg}(x,v),\operatorname{codeg}(y,v)\}\ge hn$. We claim that $|D(x,y)|\ge\tau n$.

Suppose otherwise and, exchanging the vertices if necessary, assume $|S_x|\le|S_y|$. Put $U=S_x\setminus D(x,y)$, $V=S_y\setminus D(x,y)$, $a=|U|/n$ and $s=|U\setminus V|/n$. Passing from $L_H(x)$ to $J=L_H(x)[U]$ removes at most $2\tau n^2$ edges, so $\rho(J)\ge(\theta-2\sqrt\tau)n> f(h)n$. Every vertex in $U\cap V$ has degree at most $hn$ in $J$. Also $a>h$, $s\le1-a$ and $s\le1-\theta+2\sqrt\tau<h$.

If $s=0$, then $U \subseteq V$, so $\rho(J)\le hn$, a contradiction. Otherwise the matrix $\bigl(\begin{smallmatrix}s&a-s\\s&h-s\end{smallmatrix}\bigr)$ is a quotient majorant: its Perron vector has first coordinate at least its second because $a>h$, so the weighted row bound also holds at each low-degree vertex. Collatz--Wielandt therefore gives
\[
 \frac{\rho(J)}n\le\frac{h+\sqrt{h^2+4s(a-h)}}2\le f(h),
\]
where $s(a-h)\le s(1-s-h)\le(1-h)^2/4$. This contradicts the lower bound and proves the claim.

At least $\tau n/2$ choices of $v_4\in D(x,y)$ have the same high-codegree orientation. Exchange $x,y$ if necessary so that $\operatorname{codeg}(y,v_4)\ge hn$ and $\operatorname{codeg}(x,v_4)\ge\tau n$. Choose $v_2\in N_H(xv_4)$, then $v_3\in D(v_2,v_4)$, followed by $v_1\in N_H(v_2v_3)$ and $v_5\in N_H(v_3v_4)$. There are at least $(\eta-\sqrt{2\tau})n\ge\eta n/2$ choices of $v_6\in N_H(yv_4)\cap S_{v_5}$, and at least $\tau n$ choices of $v_7\in N_H(v_5v_6)$. Excluding $x,y$ and all earlier choices loses only a bounded number at each step. The product of the numbers of choices is at least $\kappa n^7$ for large $n$. If $x,y$ were exchanged, reverse the tuple at the end.

The original order is $(v_1,\ldots,v_7)$. The absorbing order $(v_1,v_3,v_2,x,v_4,y,v_6,v_5,v_7)$ uses, in sequence, the first, fourth, fifth and third asserted edges and preserves both ends. The two paths are illustrated in Figure~\ref{fig:absorber}.
\end{proof}

The same support argument gives a connecting count. For every $\theta\in(1/2,1)$ there are $\kappa_0>0$ and $n_0$ such that, whenever $n\ge n_0$ and an $n$-vertex $3$-graph $H$ satisfies $\sigma(H)\ge\theta n$, each pair of distinct vertices $a,b$ has at least $\kappa_0n^3$ ordered triples $(x,y,z)$ of distinct vertices outside $\{a,b\}$ with $axy,byz\in E(H)$. Indeed, choose $\tau$ with $2\theta-1-2\sqrt{2\tau}>0$ and put $S_u=\{v:\operatorname{codeg}(u,v)\ge\tau n\}$. The support bound gives at least $(2\theta-1-2\sqrt{2\tau})n$ choices of $y\in S_a\cap S_b$, followed by at least $\tau n$ choices each of $x\in N_H(ay)$ and $z\in N_H(by)$. Excluding coincidences removes only $O(n^2)$ triples, giving the assertion for large $n$.

\begin{lemma}\label{old:3.4}
For every $\theta>4/7$ and $\xi>0$ there are $\beta>0$ and $n_0$ such that, if $n\ge n_0$ and $\sigma(H)\ge\theta n$, then $H$ contains a loose path $P$ on at most $\xi n$ vertices with the following property: for every even set $U\subseteq V(H)\setminus V(P)$ of size at most $\beta n$, there is a loose path on $V(P)\cup U$ with the same ends as $P$.
\end{lemma}
\begin{proof}
Let $\kappa$ and $\kappa_0$ be given by
Lemma~\ref{old:3.3} and the connecting count, respectively.
Choose
$0<\zeta\le\min\{\kappa/392,\kappa_0/120,\xi/20\}$,
set $\beta=\kappa\zeta/4$, and take $n$ sufficiently large.

Independently select each ordered seven-tuple forming the three
original path edges with probability $p=\zeta n^{-6}$.
Let $X$ count all selected tuples, and let $X_{xy}$ count those
absorbing a fixed pair of distinct vertices $x,y$.
There are at most $n^7$ eligible tuples in total and at least
$\kappa n^7$ absorbing each pair, so
$\mathbb E X\le\zeta n$ and
$\mathbb E X_{xy}\ge\kappa\zeta n$.
Chernoff's bound gives
$\Pr(X>2\zeta n)\le e^{-\zeta n/3}$ and
$\Pr(X_{xy}<\kappa\zeta n/2)\le e^{-\kappa\zeta n/8}$.
Since there are fewer than $n^2$ vertex pairs, a union bound
shows that, with probability $1-o(1)$, at most $2\zeta n$
tuples are selected and every pair has at least
$\kappa\zeta n/2$ selected absorbers.

We also need to control intersections among the selected tuples.
Let $Z$ count ordered pairs of distinct selected tuples that
share a vertex. Choosing the two positions of a common vertex
and then the vertices gives at most $49n^{13}$ possible pairs.
Each pair is selected with probability $p^2$, so
$\mathbb E Z\le49n^{13}p^2=49\zeta^2n$.
By Markov's inequality, $Z\le98\zeta^2n$ with probability
at least $1/2$. Thus this bound and all the preceding bounds
hold simultaneously with probability at least $1/2-o(1)>0$.
Fix such a selection.

Successively delete one tuple from each remaining intersecting
pair. At most $98\zeta^2n$ tuples are deleted, leaving a
vertex-disjoint family $\mathcal F$ of size $m\le2\zeta n$.
By our choice of $\zeta$, every vertex pair still has at least
$
 \kappa\zeta n/2-98\zeta^2n
 \ge \kappa\zeta n/4
$
absorbers in $\mathcal F$.

Join the paths in $\mathcal F$ consecutively using the connecting
count. At each step, avoid all vertices of $\mathcal F$ except
the two prescribed ends, together with all previously used
connecting vertices. There are at most
$7m+3m\le20\zeta n$ forbidden vertices.
Each forbidden vertex rules out at most $3n^2$ ordered
connecting triples, so at most
$60\zeta n^3\le\kappa_0n^3/2$ triples are excluded.
Since the connecting count supplies at least $\kappa_0n^3$
triples for each pair of ends, every required connection
can be made. Each connection adds three vertices, so the
resulting loose path $P$ has
$7m+3(m-1)\le20\zeta n\le\xi n$ vertices.

Now let $U\subseteq V(H)\setminus V(P)$ be even with
$|U|\le\beta n$, and partition $U$ into $q=|U|/2$ pairs.
Each pair has at least $\kappa\zeta n/4$ candidate absorbers
in $\mathcal F$, whereas
$q\le\beta n/2=\kappa\zeta n/8$.
We can therefore assign a distinct absorber to each pair
greedily: at each step, fewer than $q$ candidates have
already been used.
Replace each assigned seven-vertex path by its nine-vertex
absorbing alternative. The assigned paths are vertex-disjoint,
the added vertices lie outside $P$, and each replacement
preserves the ends of its segment. The resulting loose path
has vertex set $V(P)\cup U$ and the same ends as $P$.
\end{proof}

The next lemma provides a small reservoir of vertices for
joining vertex-disjoint loose paths.
\begin{lemma}\label{lem:reservoir}
For every $\theta\in(1/2,1)$ there is $\kappa_0>0$ such that,
for every fixed $r\in(0,1/4)$ and all sufficiently large $n$,
the following holds.
Let $H$ be an $n$-vertex $3$-graph with $\sigma(H)\ge\theta n$,
and let $F\subseteq V(H)$ satisfy $|F|\le\kappa_0n/6$.
Then there is a set $R\subseteq V(H)\setminus F$ with
$|R|\le2rn$ such that every pair of distinct vertices
$a,b\in V(H)$ has at least $\kappa_0r^3n^3/4$ ordered triples
$(x,y,z)$ of distinct vertices in $R\setminus\{a,b\}$ with
$axy,byz\in E(H)$. In particular, the ends $a,b$ may belong to $F$.
\end{lemma}

\begin{proof}
Take $\kappa_0$ from the preceding connecting count.
For any fixed pair of distinct ends $a,b$, excluding connecting
triples with a vertex in $F$ removes at most
$3|F|n^2\le\kappa_0n^3/2$ triples. Thus at least
$\kappa_0n^3/2$ connecting triples have all three internal
vertices outside $F$.

Select each vertex of $V(H)\setminus F$ independently with
probability $r$, and let $R$ be the resulting set.
For fixed ends $a,b$, let $Y$ count the connecting triples
whose internal vertices all lie in $R$.
Each available triple is retained with probability $r^3$,
so $\mathbb EY\ge\kappa_0r^3n^3/2$.
Changing the selection of one vertex affects at most $3n^2$
triples. Put $t=\kappa_0r^3n^3/4$, so that $\mathbb EY\ge2t$.
Revealing the vertex choices one at a time and applying
Azuma's inequality~\cite{Azuma} gives
\[
 \Pr(Y<t)
 \le \exp\left(-\frac{t^2}{2n(3n^2)^2}\right)
 =\exp\left(-\frac{\kappa_0^2r^6n}{288}\right).
\]
Also, Chernoff's bound gives
$\Pr(|R|>2rn)\le e^{-rn/3}$.
A union bound over fewer than $n^2$ pairs of ends shows that
the probability of any failure is at most
$n^2\exp(-\kappa_0^2r^6n/288)+e^{-rn/3}=o(1)$.
Hence there is a choice of $R$ satisfying all the required
bounds simultaneously.

Finally, any fixed number $k$ of pairwise disjoint prescribed
end pairs can be connected greedily through $R$.
At step $i$, avoid all $2k$ prescribed ends and the
$3(i-1)$ internal vertices used earlier.
These exclusions rule out at most
$3(2k+3(i-1))n^2\le15kn^2$ connecting triples.
For sufficiently large $n$, this is smaller than
$\kappa_0r^3n^3/4$, so a connection remains available at
every step. The resulting connecting paths are
pairwise vertex-disjoint.
\end{proof}

For pairwise disjoint vertex sets $A,B,C$, let $d(A,B,C)$
denote the proportion of triples in $A\times B\times C$
that are edges.
We call $(A,B,C)$ $\epsilon$-regular if
$|d(A',B',C')-d(A,B,C)|\le\epsilon$
whenever $A'\subseteq A$, $B'\subseteq B$ and $C'\subseteq C$
satisfy $|A'|\ge\epsilon|A|$, $|B'|\ge\epsilon|B|$ and
$|C'|\ge\epsilon|C|$.

\begin{lemma}\label{old:3.5}
Suppose that $(A,B,C)$ is $\epsilon$-regular with density $d\ge d_0>0$ and class sizes $r,r,2r$. If $\epsilon$ is sufficiently small compared with $d_0$ and $r$ is sufficiently large, a loose path leaves at most $16\epsilon r/d_0+O(1)$ vertices uncovered.
\end{lemma}
\begin{proof}
Put $\delta=4\epsilon/d_0$. Construct a path whose ends and intersection vertices, in their
path order, alternate between $A$ and $B$, with one vertex of
$C$ in each edge.
Write $A',B',C'$ for the sets of unused vertices in the three classes. Require its current end to have at least $(d-2\epsilon)|B'||C'|$ choices into the unused opposite class and unused $C$-class, or the analogous condition when its end lies in $B$. Start with a typical vertex in $A$.

Throughout the construction,
$\bigl||A'|-|B'|\bigr|\le1$ and
$|C'|=|A'|+|B'|+1$.
Thus $|A'|,|B'|\ge\delta r$ also ensures $|C'|\ge2\delta r$. 
While $|A'|,|B'|\ge\delta r$, suppose that the current end
lies in $A$; the other case is symmetric.
Let $B_0\subseteq B'$ consist of vertices whose degree into
$A'\times C'$ is less than $(d-\epsilon)|A'||C'|$.
If $|B_0|\ge\epsilon r$, averaging these degrees gives
$d(A',B_0,C')<d-\epsilon$.
Since $|A'|\ge\epsilon r$ and $|C'|\ge2\epsilon r$,
this contradicts the $\epsilon$-regularity of $(A,B,C)$.
Thus $|B_0|<\epsilon r$. The current end has an edge through a vertex of $B'\setminus B_0$
because $(d-2\epsilon)|B'|>\epsilon r$. Choose such an edge. Deleting its private vertex reduces the new end's degree by at most $|A'|$, and
\[
 (d-\epsilon)|A'||C'|-|A'|\ge(d-2\epsilon)|A'|(|C'|-1)
\]
for large $r$, since $\epsilon|C'|+d-2\epsilon-1\ge0$. The invariant persists. Stop when one of $A',B'$ has fewer than $\delta r$ vertices.
The above identities then give
$|A'|+|B'|+|C'|<4\delta r+3
=16\epsilon r/d_0+3$,
as required.
\end{proof}

\begin{lemma}\label{old:3.6}
For every $\eta>0$ there is $\varepsilon_0>0$ such that,
for every $\alpha>0$, there are $T,n_0$ with the following
property. If an $n$-vertex $3$-graph $H$ satisfies
$n\ge n_0$, $\sigma(H)\ge(c-\varepsilon_0)n$, and every
$\lfloor3n/4\rfloor$-set spans more than $\eta n^3$ edges,
then $H$ contains at most $T$ vertex-disjoint loose paths
leaving at most $\alpha n$ vertices uncovered.
\end{lemma}

\begin{proof}
Choose $\delta_0>0$ with $C\delta_0\ll\eta$, where $C$ is from
Lemma~\ref{old:3.1}, and then $\varepsilon_0\ll\delta_0$.
These choices depend only on $\eta$.
Given $\alpha$, first choose a density cutoff $d>0$ and a lower bound $t_0$
on the number of classes, and then choose a regularity parameter
$\epsilon>0$ such that
$d+1/t_0+\sqrt{\epsilon}\ll\delta_0^4,\eta,\alpha^2$ and
$\sqrt{\epsilon}/d\ll\min\{\alpha,1\}$.

Apply weak regularity~\cite[Theorem~2.4]{HanZhao}, obtaining an
exceptional set $V_0$ of size at most $\epsilon n$ and
$t_0\le t\le T$ equal classes of size $m$.
Let $H_0$ contain only the edges across triples of distinct classes
that are $\epsilon$-regular and have density at least $d$.
Edges meeting $V_0$, using a class twice, or lying in irregular
or low-density triples account for at most $\mu n^3$ deletions,
where $\mu=10(d+1/t_0+\epsilon)$ suffices.
Call a vertex bad if it loses more than $\sqrt\mu n^2$ incident
edges. Counting deleted incidences gives at most $3\sqrt\mu n$
bad vertices. Discard classes consisting entirely of bad vertices;
since $tm\ge(1-\epsilon)n$, there are at most $4\sqrt\mu t$ such
classes. Let $R$ have the retained $t'$ classes as vertices,
with edges corresponding to their regular triples of density
at least $d$.

For each retained class $V_i$, choose a nonbad vertex $v$.
Its link in $H_0$, restricted to the retained classes, differs
from its original link by at most $10\sqrt\mu n^2$ deleted edges,
with omitted vertices viewed as isolated.
The Frobenius bound gives spectral radius at least
$(c-\varepsilon_0-20\mu^{1/4})n$.
Ignoring isolated vertices, this link is a subgraph of the equal
$m$-fold blow-up of $L_R(i)$, whose radius is
$m\rho(L_R(i))$.
Indeed, the blow-up adjacency matrix is
$A(L_R(i))\otimes\mathbf J_m$, with $\mathbf J_m$ the all-ones matrix.
Since $mt'\le n$, this gives
\[
 \sigma(R)\ge(c-\varepsilon_0-20\mu^{1/4})t'
 \ge(c-\delta_0)t'.
\]
Only one nonbad representative per class is needed; no regularity
of individual links is assumed.

Suppose $R$ has no perfect fractional oriented-edge packing.
Lemma~\ref{old:3.1} gives a vertex cover of size at most
$(1/4+C\delta_0)t'$.
The union $B_0$ of the other classes contains no edge of $H_0$,
so it spans at most $\mu n^3$ edges of $H$, while
$|B_0|\ge(3/4-C\delta_0-O(\sqrt\mu)-\epsilon)n$.
Enlarge it to a $\lfloor3n/4\rfloor$-set if necessary, adding at most
$(C\delta_0+O(\sqrt\mu)+\epsilon)n$ vertices, each contributing
at most $n^2/2$ edges. If it is already larger, take a subset.
Our parameter choices give a set spanning fewer than $\eta n^3$
edges, a contradiction.

Thus $R$ has a perfect fractional packing.
Choose one with minimum support. Its load vectors are linearly
independent, since otherwise a load-preserving perturbation would
eliminate a positive weight. Hence its support has size at most $t'$.
Discard weights below $\sqrt\epsilon$, losing total vertex load
at most $4\sqrt\epsilon t'$.
For each remaining oriented edge $e$ of weight $w_e$, put
$r_e=\lfloor w_em\rfloor$ and allocate $2r_e,r_e,r_e$ vertices
in its distinguished and other classes.
The load constraints allow disjoint allocations, and rounding
loses at most $4t'$ vertices.

Every allocated part has size at least $(\sqrt\epsilon/2)m$
for large $n$. Slicing gives regularity parameter
$4\sqrt\epsilon$ and density at least $d/2$:
test subsets of relative size $4\sqrt\epsilon$ have at least
$2\epsilon m$ vertices, and two applications of the original
regularity bound the density difference by $2\epsilon$.
Apply Lemma~\ref{old:3.5}, with the distinguished class as its
class of twice the size.
Each allocation yields one path leaving an
$O(\sqrt\epsilon/d)$ proportion uncovered.
There are at most $t'\le T$ disjoint paths, and the total
number of uncovered vertices is at most
\[
 \bigl(O(\sqrt\mu)+\epsilon+4\sqrt\epsilon
       +O(\sqrt\epsilon/d)\bigr)n+O(T)\le\alpha n.
\]
After regularity fixes $T$, take $n_0$ large enough for these
applications and the bounded rounding errors.
\end{proof}

\begin{proof}[Proof of Theorem~\ref{old:3.2}]
Since $c>4/7$, fix $4/7<\theta<c$ and take $\kappa_0$ from
Lemma~\ref{lem:reservoir}.
Obtain $\varepsilon_0$ from Lemma~\ref{old:3.6} with
nonextremality parameter $\eta/2$, and choose
$\varepsilon\ll\varepsilon_0,c-\theta$ and
$\xi\ll\varepsilon_0^2,\eta,\kappa_0$.
Lemma~\ref{old:3.4} gives an absorbing path $P$ on at most
$\xi n$ vertices, with capacity $\beta n$.
Choose $r\ll\beta,\varepsilon_0^2,\eta$.
Applying Lemma~\ref{lem:reservoir} with $F=V(P)$ gives a
reservoir $R\subseteq V(H)\setminus V(P)$ of size at most $2rn$.

Put $H'=H-(V(P)\cup R)$ and $N=|V(H')|=(1-q)n$, where
$q\le\xi+2r$.
Each remaining link loses at most $qn^2$ edges, so the
Frobenius bound and our parameter choices give
\[
 \sigma(H')\ge(c-\varepsilon-\sqrt{2q})n
 \ge(c-\varepsilon_0)N.
\]
If a $\lfloor3N/4\rfloor$-set in $H'$ spans at most
$(\eta/2)N^3$ edges, enlarging it by at most $qn$ vertices
gives a $\lfloor3n/4\rfloor$-set in $H$ spanning at most
$(\eta/2)N^3+qn^3/2<\eta n^3$ edges, a contradiction.
Thus $H'$ also satisfies nonextremality with parameter $\eta/2$.

Apply Lemma~\ref{old:3.6} with uncovered proportion
$\alpha\ll\beta$ to obtain a bounded family of vertex-disjoint
loose paths in $H'$ leaving at most $\alpha N$ vertices uncovered.
Using disjoint connectors from $R$, join these paths and $P$
cyclically to obtain a loose cycle $C$.
The uncovered set $U$ consists of the vertices missed in $H'$
and the unused vertices of $R$, so
$|U|\le\alpha N+2rn\le\beta n$.
Since both $n$ and $|V(C)|$ are even, $|U|$ is even.
The absorbing property of $P$ therefore gives a path on
$V(P)\cup U$ with the same ends as $P$.
Replacing $P$ by this path in $C$ produces a loose Hamilton cycle.
\end{proof}

\subsection{Exact completion near the space barrier}
For a set $U$, let $d_H(v,U)$ count edges containing $v$ and two vertices of $U\setminus\{v\}$, and let $\bar d_H(v,U)$ count missing pairs in $\binom{U\setminus\{v\}}2$. For disjoint $U,W$, let $\bar d_H(v,UW)$ count pairs $(u,w)\in(U\setminus\{v\})\times(W\setminus\{v\})$ for which $vuw\notin E(H)$. Omit the subscript when the hypergraph is clear. Write $S_a(N)=K_a\vee I_{N-a}$ and
\begin{equation}\label{eq:old:4}
 \lambda_a(N)=\rho(S_a(N))=\frac{a-1+\sqrt{(a-1)^2+4a(N-a)}}2.
\end{equation}
Throughout this subsection, $N=n-1$ and $a_0=\lceil n/4\rceil-1$.

\begin{theorem}\label{old:3.7}
There are $\beta>0$ and $n_0$ such that, for every even $n\ge n_0$, an $n$-vertex $3$-graph $H$ with $\sigma(H)>\tau_n$ contains a loose Hamilton cycle whenever some set $B$ of size $\lfloor3n/4\rfloor$ satisfies $e(H[B])\le\beta n^3$.
\end{theorem}

The hypothesis permits missing edges meeting $V(H)\setminus B$. The following estimate controls them directly from the link spectra.

\begin{lemma}\label{old:3.8}
Let $N$ be sufficiently large, $N/5\le h\le N/3$, and $J=S_h(N)$. Obtain $G$ by deleting $M\subseteq E(J)$ and adding a graph $T$ inside its independent part. If $\rho(G)\ge\rho(J)-2$ and $\rho(G)\ge N/3$, then $|M|\le C(e(T)+N)$ for an absolute constant $C$.
\end{lemma}
\begin{proof}
Put $\lambda_0=\rho(J)$ and take nonnegative unit Perron vectors $x,u$ of $G,J$. The eigenvalue equation and Cauchy--Schwarz give $x_i\le\sqrt N/\rho(G)\le3/\sqrt N$. If $p,q$ are the clique and independent coordinates of $u$, then $1\le p/q=\lambda_0/h\le5$, so every $u_i\ge1/(5\sqrt N)$.

All eigenvalues of $J$ other than $\lambda_0$ are nonpositive: its other quotient eigenvalue is negative, and the remaining eigenvalues are $-1$ and $0$. Also $\lambda_0\ge N/3$, by the constant-vector Rayleigh quotient. Put $\alpha=\langle x,u\rangle$ and $R=2+18e(T)/N$.
Writing $x=\alpha u+y$ with $y\perp u$ gives
$x^\top A(J)x\le\lambda_0\alpha^2$.
Since $A(G)=A(J)-A(M)+A(T)$ and
$x^\top A(T)x\le18e(T)/N$, the hypothesis implies
\[
 \lambda_0(1-\alpha^2)+x^\top A(M)x\le R.
\]
Thus $\|x-u\|^2=2(1-\alpha)\le6R/N$. Only $O(R)$ vertices have $x_i<u_i/2$, and at most $O(RN)$ missing edges meet them. Each other missing edge contributes at least $1/(50N)$ to $x^\top A(M)x\le R$. There are therefore $O(RN)$ missing edges in total, proving the assertion.
\end{proof}

We use the prescribed-end completion lemma of Han and Zhao~\cite[Lemma~3.7]{HanZhao}, with the same partition relation $|Z|=3(|X|-1)$ and prescribed distinct ends in $X$. Their complement count for $x\in X$ is exactly $\bar d_F(x,Z)$. For $z\in Z$, it also counts the $|X|$ impossible pairs containing $z$, so their count is exactly $\bar d_F(z,XZ)+|X|$. If $|Z|\ge1/\varrho$, the two bounds below therefore imply their hypotheses with error parameter $2\varrho$. This proves the stated version for sufficiently small $\varrho$ and sufficiently large order.

\begin{lemma}[Han--Zhao]\label{old:3.9}
For every sufficiently small $\varrho>0$ and all sufficiently large orders, the following holds. Let $F$ be a $3$-graph with partition $V(F)=X\mathbin{\dot\cup}Z$, where $|Z|=3(|X|-1)$. Suppose that $\bar d_F(x,Z)\le\varrho\binom{|Z|}2$ for every $x\in X$ and $\bar d_F(z,XZ)\le\varrho|X||Z|$ for every $z\in Z$. For any distinct $x_0,x_1\in X$, there is a loose Hamilton path in $F$ with ends $x_0,x_1$.
\end{lemma}

\begin{proof}[Proof of Theorem~\ref{old:3.7}]
Fix $C$ from Lemma~\ref{old:3.8} and put $L=2C+10$. Choose $\varrho_*>0$ so that $2\varrho_*$ lies in the range permitted by Han--Zhao's lemma; then Lemma~\ref{old:3.9} applies for $\varrho\le\varrho_*$ and sufficiently large order. Choose $0<\beta\ll\eta\ll\epsilon\ll1/L,\varrho_*$ with $10L\epsilon<\varrho_*$, and then choose $n$ sufficiently large. In particular, require $\beta\ll\eta\epsilon/C$ and $n\gg C/(\eta\epsilon)$.

Write $b=\lfloor3n/4\rfloor$ and $h=n-b=\lceil n/4\rceil$. Among all $b$-sets choose $B$ minimizing $e(H[B])$, and put $A=V(H)\setminus B$. Then $e(H[B])\le\beta n^3$. For large $n$, $\tau_n\ge\lambda_{h-1}(N)>n/2$. Regarding~\eqref{eq:old:4} as a function of real $a$, its derivative is
\[
 \frac12+\frac{4N-2-6a}{4\sqrt{(a-1)^2+4a(N-a)}},
\]
which lies between $1/2$ and $2$ for $N/5\le a\le N/3$. Consequently, whenever $a,a+k\in[N/5,N/3]$ with $k\ge0$,
we have $k/2\le\lambda_{a+k}(N)-\lambda_a(N)\le2k$.
Thus increasing the clique size by one raises the spectral
radius by at most $2$, while increasing it by $k$ raises
it by at least $k/2$.

Define $A'=\{v:d(v,B)\ge(1-\epsilon)\binom b2\}$, $B'=\{v:d(v,B)\le\epsilon\binom b2\}$ and $V_0=V(H)\setminus(A'\cup B')$. The first two sets are disjoint.

\begin{claim}\label{old:3.10}
Each of $|A\setminus A'|$, $|B\setminus B'|$, $|A'\setminus A|$ and $|B'\setminus B|$ is at most $\eta n$, and $|V_0|\le2\eta n$. Moreover:
\begin{enumerate}
\item $A\cap B'\ne\varnothing$ implies $B\subseteq B'$, whereas $B\cap A'\ne\varnothing$ implies $A\subseteq A'$;
\item $\bar d(x,(A'\setminus\{x\})(B'\setminus\{x\}))\le L\epsilon n^2$ for every $x\in B'$;
\item $\Delta_{1}(H[B'])\le\epsilon n^2$.
\end{enumerate}
\end{claim}
\begin{proof}
For $v\in B$, compare $L_H(v)$ with the split graph having clique $A$ and independent set $B\setminus\{v\}$. Its radius satisfies $\lambda_h(N)\le\lambda_{h-1}(N)+2<\rho(L_H(v))+2$. The added internal edges number $d(v,B)$, so Lemma~\ref{old:3.8} bounds its missing split edges by $C(d(v,B)+N)$. Each missing $ABB$ triple is counted twice in the sum over $v\in B$. Thus their number $D$ satisfies $D\le C(3e(H[B])+bN)/2\le2C\beta n^3+Cn^2$.

By definition, each vertex of $A\setminus A'$ accounts for more than $\epsilon\binom b2$ such missing triples. Since $\binom b2\ge n^2/5$, this gives $|A\setminus A'|\le D/(\epsilon n^2/5)\le10C\beta n/\epsilon+5C/\epsilon\le\eta n$. Similarly, $\sum_{v\in B}d(v,B)=3e(H[B])$ gives $|B\setminus B'|\le15\beta n/\epsilon\le\eta n$. Disjointness of $A',B'$ yields the other bounds.

For the first exchange assertion, if $a\in A\cap B'$ and $v\in B\setminus B'$, replacing $v$ by $a$ changes $e(H[B])$ by at most $-d(v,B)+d(a,B)<0$, contradicting minimality. Likewise, if $v\in B\cap A'$ and $a\in A\setminus A'$, then $d(a,B)<d(v,B)$ and the same replacement again contradicts minimality.

For $x\in B'$, compare its link with the split graph on $A\setminus\{x\}$ and $B\setminus\{x\}$. Its clique size is $h$ or $h-1$, and its radius is at most $\rho(L_H(x))+2$. It has at most $\epsilon\binom b2$ added internal edges. Lemma~\ref{old:3.8} and $n\ge2/\epsilon$ bound its missing split edges by $C\epsilon n^2$. Passing to $A',B'$ adds at most $2\eta n^2$ possibly missing pairs, proving the second assertion. Finally, $d_{H[B']}(x)\le d(x,B)+|B'\setminus B|n\le(\epsilon/2+\eta)n^2\le\epsilon n^2$.
\end{proof}

\begin{claim}\label{old:3.11}
For distinct $u,v\in B'$ and $W\subseteq V(H)\setminus\{u,v\}$ with $|W|\le n/8$, there are $a\in A'\setminus W$ and distinct $z_1,z_2\in B'\setminus(W\cup\{u,v\})$ such that $uz_1a,az_2v\in E(H)$.
\end{claim}
\begin{proof}
On $A'\setminus W$ and $B'\setminus(W\cup\{u,v\})$, retain the pairs $az$ for which both $uaz$ and $vaz$ are edges. The two parts have sizes at least $h-\eta n-n/8$
and $b-\eta n-n/8-2$, respectively, and hence at least
$n/9$ and $n/2$ for sufficiently small $\eta$ and large $n$. Claim~\ref{old:3.10} bounds missing pairs by $2L\epsilon n^2$, so at least $n^2/36$ pairs remain. Since this exceeds the size of the first part, some $a$ has two distinct neighbours.
\end{proof}

\begin{claim}\label{old:3.12}
Put $q=|A\cap B'|$. If $q>0$, then $H[B']$ has a collection of vertex-disjoint loose paths consisting of one edge if $q=1$ and $n\equiv2\pmod4$, two edges whose intersection has size at most one if $q=1$ and $n\equiv0\pmod4$, and $2q$ vertex-disjoint edges if $q\ge2$.
\end{claim}
\begin{proof}
Claim~\ref{old:3.10} gives $B\subseteq B'$ and $|B'|=b+q$. Put $a'=n-|B'|=h-q$. First suppose $q\ge2$. For $v\in B'$, let $G=L_H(v)$ and $d=d_{H[B']}(v)$. Deleting its $d$ internal edges leaves a subgraph of $S_{a'}(N)$. A nonnegative unit Perron vector $x$ of $G$ satisfies $\max_i x_i^2\le N/\rho(G)^2<4/n$, so $\rho(G)\le\lambda_{a'}(N)+8d/n$. Since $a_0-a'=q-1$, the derivative bound gives $\lambda_{a_0}(N)-\lambda_{a'}(N)\ge(q-1)/2$. Thus $\rho(G)>\tau_n\ge\lambda_{a_0}(N)$ implies $d>(q-1)n/16$. Summing yields $e(H[B'])>(q-1)n^2/64$, because $|B'|\ge3n/4$. After fewer than $2q$ disjoint edges have been chosen, at most $6q\epsilon n^2$ edges meet their vertices, by Claim~\ref{old:3.10}. This is smaller than $(q-1)n^2/64$ when $\epsilon<1/768$, so greedy selection gives $2q$ disjoint edges.

Now let $q=1$, so $a'=a_0$. If $n\equiv2\pmod4$ and $H[B']$ has no edge, each link rooted in $B'$ is a subgraph of $S_{a_0}(N)$, contradicting the spectral hypothesis. Suppose instead that $n\equiv0\pmod4$ and the required pair does not exist. Then any two distinct edges of $F=H[B']$ meet in two vertices. Either all edges of such a $3$-graph contain a common pair or they are supported on at most four vertices. Indeed, label two edges $123,124$. If not all edges contain $12$, a third is $134$ or $234$. An edge using a vertex outside $\{1,2,3,4\}$ would then need its other two vertices in the common intersection of these three edges, which has size one.

In the common-pair case, a vertex $v\in B'$ outside the pair has at most one internal link edge. Its link is a subgraph of $K_{a_0}\vee(K_2\cup I_{N-a_0-2})$. In the four-vertex case, a vertex outside the support has link contained in $S_{a_0}(N)$. Both conclusions contradict $\rho(L_H(v))>\tau_n$.
\end{proof}

We construct a short loose path $P$ covering $V_0$, with both ends in $B'$, such that
\begin{equation}\label{eq:old:5}
 3|A'\setminus V(P)|-|B'\setminus V(P)|\ge1,
 \qquad |V(P)|\le30\eta n.
\end{equation}
Put $t=|V_0|$. If $q>0$, take the correcting collection
from Claim~\ref{old:3.12}, and let $e$ and $g$ be its numbers
of edges and components. If $q=0$, take the empty collection
and put $e=g=0$. These paths lie in $B'$ and use
$2e+g\le6q$ vertices.

For each $x\in V_0$, choose two disjoint pairs in its link,
using four new vertices of $B'$. They form a two-edge loose
path whose common vertex is $x$ and whose ends lie in $B'$.
Indeed, $d(x,B)>\epsilon\binom b2\ge\epsilon n^2/5$.
Excluding $B\setminus B'$ and all previously used vertices
removes at most $(\eta n+6q+4t)n\le15\eta n^2$ pairs.
Since $\eta\ll\epsilon$, more than $\epsilon n^2/6$ pairs
remain. A fixed pair meets at most $2n$ other pairs, so two
remaining pairs are disjoint for large $n$.

Suppose first that $t+g>0$. Join these $t+g$ paths using
Claim~\ref{old:3.11}, avoiding every initial path and previous
connector except the two prescribed ends.
Each of the $t+g-1$ connectors adds one vertex of $A'$
and two of $B'$. The initial paths use at most $5t+6q$
vertices, and $g\le2q$, so throughout the construction at most
$5t+6q+3(t+g-1)\le8t+12q\le28\eta n$ vertices are used.
Thus every forbidden set has size below $n/8$.
The resulting path $P$ covers $V_0$ and has both ends in $B'$.

Let $\ell=3|A'\setminus V(P)|-|B'\setminus V(P)|$.
Relative to $3|A'|-|B'|$, the paths through $V_0$ contribute
$4t$, the correcting collection contributes $2e+g$,
and each connector contributes $-1$.
Since $|B'|\le b+q$ and $|A'|=n-t-|B'|\ge h-q-t$,
\[
\begin{aligned}
 \ell
 &=3|A'|-|B'|+4t+(2e+g)-(t+g-1)\\
 &=3|A'|-|B'|+3t+2e+1
 \ge3h-b-4q+2e+1\ge1.
\end{aligned}
\]
The last inequality holds because $e=2q$ and $3h-b\ge0$,
except when $q=1$ and $n\equiv2\pmod4$, in which case
$e=1$ and $3h-b=2$.

If $q=t=0$, take instead a two-edge path from
Claim~\ref{old:3.11}, with common vertex in $A'$ and the
other four vertices in $B'$.
Then $\ell=3|A'|-|B'|+1\ge3h-b+1\ge1$.
Thus in all cases $P$ satisfies~\eqref{eq:old:5}.

Since $P$ covers $V_0$ and has odd order,
$\ell\equiv n-|V(P)|\equiv1\pmod2$.
Also $|A'|\le h+\eta n$ and $|B'|\ge b-\eta n$, so
$\ell\le3|A'|-|B'|+|V(P)|
\le2+34\eta n\le36\eta n$.
While the residual difference exceeds one, apply
Claim~\ref{old:3.11} to a current end $u\in B'$ and a fresh
auxiliary end in $B'$, avoiding the other vertices of the
current path.
Keep the first edge $uaz$, where $a\in A'$ and $z\in B'$
are new, and append it in the order $u,a,z$.
The new end remains in $B'$, and the residual difference
decreases by two. After $(\ell-1)/2$ extensions, we obtain
a path $P'$ with residual difference one, both ends in $B'$,
and at most $66\eta n$ vertices.

Extend each end once more in the same way, now taking the
new $A'$ vertex as the end. Let the resulting path be $Q$,
with ends $x_0,x_1\in A'$.
At most $66\eta n+4<70\eta n<n/8$ vertices are used throughout,
so all applications of Claim~\ref{old:3.11} are valid.
Put $X=(A'\setminus V(Q))\cup\{x_0,x_1\}$
and $Z=B'\setminus V(Q)$.
The final extensions use two vertices of each class,
but their $A'$ ends are included in $X$.
Thus $|X|=|A'\setminus V(P')|$ and
$|Z|=|B'\setminus V(P')|-2=3|X|-3=3(|X|-1)$;
also $|X|\ge n/5$ and $|Z|\ge2n/3$.

For $x\in X\subseteq A'$, the definition of $A'$ and
$|Z\setminus B|\le\eta n$ give
$\bar d(x,Z)\le\epsilon\binom b2+\eta n^2
\le\epsilon n^2\le5\epsilon\binom{|Z|}2$.
For $z\in Z$, Claim~\ref{old:3.10} gives
$\bar d(z,XZ)\le L\epsilon n^2\le10L\epsilon|X||Z|$.
Lemma~\ref{old:3.9}, applied to $H[X\cup Z]$ with
$\varrho=10L\epsilon$, supplies a loose Hamilton path
with ends $x_0,x_1$.
This path and $Q$ have only their ends in common and together
cover $V(H)$, so they form a loose Hamilton cycle.
\end{proof}

\begin{proof}[Proof of Theorem~\ref{thm:loose}]
Choose $\beta$ from Theorem~\ref{old:3.7} and apply Theorem~\ref{old:3.2} with $\eta=\beta$, obtaining $\varepsilon>0$. Since $\tau_n=cn+O(1)$, the hypothesis $\sigma(H)>\tau_n$ implies $\sigma(H)\ge(c-\varepsilon)n$ for sufficiently large $n$. If some $\lfloor3n/4\rfloor$-set spans at most $\beta n^3$ edges, apply Theorem~\ref{old:3.7}; otherwise apply Theorem~\ref{old:3.2}. Proposition~\ref{prop:sharp} proves sharpness.
\end{proof}

\begin{remark}\label{old:3.13}
Equality need not determine a unique hypergraph. Deleting all triples wholly contained in $A$ from either loose extremal construction leaves every critical $B$-link unchanged. Each $A$-link still contains the clique on $B$, whose radius exceeds $\tau_n$ for large $n$. These subgraphs therefore also have minimum link radius $\tau_n$ and no loose Hamilton cycle.
\end{remark}

\section{Tight Hamilton cycles: stability}\label{sec:stability}

\subsection{Sampling the barrier densities}
For $\beta\in\{1/3,2/3\}$ and $P\subseteq\{0,1,2,3\}$, define
\[
 b_{\beta,P}(H)=\min_{|A|=\lfloor\beta n\rfloor}
 \frac{|\{e\in E(H):|e\cap A|\in P\}|}{n^3},\qquad n=|V(H)|.
\]
Write $b_1=b_{1/3,\{0\}}$, $b_2=b_{2/3,\{1\}}$, and
$b_3=b_{1/3,\{1\}}$. Moving $O(\delta n)$ vertices between the two
parts changes the density by $O(\delta)$, so the three barriers considered
below are described by $\min_i b_i(H)=o(1)$.

\begin{lemma}\label{lem:barrier-sampling}
Fix $\beta\in\{1/3,2/3\}$ and $P\subseteq\{0,1,2,3\}$.
For every $\delta>0$, all sufficiently large $s$ and then
sufficiently large $n$ have the following property.
If $H$ is an $n$-vertex $3$-graph and $U$ is a uniform
$s$-subset of $V(H)$, then
\begin{equation}\label{eq:old:6}
 \Pr\bigl(|b_{\beta,P}(H[U])-b_{\beta,P}(H)|>\delta\bigr)<\delta.
\end{equation}
\end{lemma}

\begin{proof}
We apply the microcanonical sampling theorem of Karpinski and
Mark\'o~\cite[Theorem~5.4]{KarpinskiMarko}, keeping track of the
prescribed part sizes. Use arity three, two states, and mass vector
$\boldsymbol a=(\beta,1-\beta)$. Partition $[0,1]$ into equal cells
$I_v$, indexed by $v\in V(H)$. For $\boldsymbol i\in\{1,2\}^3$, set
$W_{\boldsymbol i}(x_1,x_2,x_3)=-1/6$ when $x_j\in I_{v_j}$,
$\{v_1,v_2,v_3\}\in E(H)$ and
$|\{j:i_j=1\}|\in P$; set it to zero otherwise, in particular on repeated
cell indices. For measurable $f_1,f_2\ge0$ with $f_1+f_2=1$ and
$\int f_i=a_i$, the energy is
\[
 \sum_{\boldsymbol i\in\{1,2\}^3}
 \int_{[0,1]^3}W_{\boldsymbol i}(x_1,x_2,x_3)
 \prod_{j=1}^3 f_{i_j}(x_j)\,dx_1\,dx_2\,dx_3.
\]
Write $E_{\boldsymbol a}(W)$ for its maximum. Averaging each $f_i$ on the cells preserves both the energy and
the masses. Writing $p_v=n\int_{I_v}f_1$, we obtain vertex weights
$p_v\in[0,1]$ with $\sum_v p_v=\beta n$.
Thus $E_{\boldsymbol a}(W)$ is exactly the finite fractional optimum. For an integral colouring with
$|A|=\lfloor\beta n\rfloor$, each counted edge contributes $-1/n^3$.
The integer mass convention in \cite[Definition~5.1]{KarpinskiMarko}
allows $\bigl||A|-\beta n\bigr|\le1$, so one recolouring gives the
required size at a cost of at most $\binom{n-1}{2}/n^3$. Together with
\cite[Lemma~5.3]{KarpinskiMarko}, this gives
$|E_{\boldsymbol a}(W)+b_{\beta,P}(H)|\le C/n$
for an absolute constant $C$.

Let $G_s$ be the weighted sample obtained from $s$ independent uniform
points of $[0,1]$, and let $\widehat E_{\boldsymbol a}(G_s)$ denote its
integer optimum with the same mass convention. The source theorem applies
to these signed kernels, whose maximum absolute entry is at most $1/6$.
Thus, for any $\eta>0$ and sufficiently large $s$,
$\Pr(|E_{\boldsymbol a}(W)-\widehat E_{\boldsymbol a}(G_s)|>\eta)<\eta$.
The sample-size threshold depends only on $\eta$, the arity $3$, the
number of states $2$ and the uniform kernel bound $1/6$, and is independent
of $W$, $H$ and $n$.
Let $\mathcal D$ be the event that the sampled cells are distinct. Then
$\Pr(\mathcal D)\ge1-\binom{s}{2}/n$, and, conditional on $\mathcal D$, their vertex set
$U$ is uniform. The same recolouring estimate gives
$|\widehat E_{\boldsymbol a}(G_s)+b_{\beta,P}(H[U])|\le C/s$ on $\mathcal D$.
Choose $0<\eta<\min\{\delta/4,1/4\}$, then $s$ so that the sampling
bound holds and $C/s\le\delta/4$, and finally $n$ so that
$C/n\le\delta/4$ and $\Pr(\mathcal D)\ge1/2$. The triangle inequality bounds
the probability in~\eqref{eq:old:6} by $\eta/\Pr(\mathcal D)\le2\eta<\delta$.
\end{proof}

\begin{lemma}\label{old:4.1}
For every fixed $\eta>0$, all sufficiently large $s$ and then sufficiently
large $n$ satisfy, uniformly over all $n$-vertex $3$-graphs $H$ and
$i\in\{1,2,3\}$,
\[
 \Pr\bigl(|b_i(H[U])-b_i(H)|>\eta\bigr)
 \le 2\exp(-c\eta^2s),
\]
where $U$ is a uniform $s$-subset and $c>0$ is an absolute constant.
\end{lemma}
\begin{proof}
Replacing one sampled vertex, while giving its replacement the same state,
preserves the part sizes and affects at most $\binom{s-1}{2}$ triples.
Thus the optimum changes by at most $1/(2s)$. The exposure martingale has
the same difference bound: couple completions after two possible next
vertices by exchanging those vertices, so the final sets agree or differ
by one replacement. Applying Lemma~\ref{lem:barrier-sampling} with error parameter
$\eta/4$ and using $0\le b_i\le1/6$, we obtain
$|\mathbb E b_i(H[U])-b_i(H)|
\le\eta/4+\eta/24<\eta/2$. Azuma's inequality proves the
claim, for instance with $c=1/8$ after weakening constants.
\end{proof}

\begin{lemma}\label{old:4.2}
For every $\eta>0$, all sufficiently large $s$ and then sufficiently large
$n$ have the following property. If $b_i(H)\ge\eta$ for $i=1,2,3$, then
every six-set $Q$ is contained in at least
$(1-1/(2s^2))\binom{n-6}{s-6}$ sets $S$ of size $s$ with
$b_i(H[S])\ge\eta/2$ for all three $i$.
\end{lemma}
\begin{proof}
Deleting or adding $r$ vertices changes $b_i$ by $O(r/n)$: restrict or
extend a balanced partition, recolour at most $r$ vertices to restore its
size, and count the affected triples and the change of normalization.
Hence $b_i(H-Q)\ge7\eta/8$ for large $n$. A uniform $(s-6)$-subset $U$
outside $Q$ satisfies $b_i(H[U])\ge3\eta/4$ simultaneously, except with
probability at most $6\exp(-c'\eta^2(s-6))$, by Lemma~\ref{old:4.1}.
Adding $Q$ changes each density by at most $\eta/8$ for large $s$. Increase $s$ until the
exceptional probability is below $1/(2s^2)$. The estimates are uniform in
$Q$.
\end{proof}

\subsection{The structural trichotomy}

Two triples are tightly connected if a sequence of triples joins them with
consecutive triples sharing two vertices. A tight component is a maximal
class of triples under this relation, together with its incident vertices.
A fractional matching assigns nonnegative edge weights whose sum at each
vertex is at most one; it is perfect when equality holds everywhere. A
fractional vertex cover assigns nonnegative vertex weights summing to at
least one on each edge. Linear programming duality equates the maximum
matching weight and minimum cover weight. For a graph $J$, let $\nu^*(J)$
be its fractional matching number. A closed tight walk is a cyclic vertex
sequence, allowing repetitions, in which every three consecutive positions
form an edge. Its order is the number of positions in the cyclic sequence,
counted with multiplicity, or equivalently the number of ordered-pair
transitions around the walk.

The sequential stability statements below, including Lemma~\ref{old:4.11},
have the following uniform meaning. For every output error $\eta>0$ there
are $\varepsilon>0$ and $N_0$ such that, for orders at least $N_0$, input
errors at most $\varepsilon$ give the stated conclusions with output errors
at most $\eta$, after normalizing each error by the corresponding power of
the order. Indeed, failure for a fixed $\eta$ would give counterexamples
of orders $N_j\ge j$ with input errors at most $1/j$, contradicting the
sequential statement. When alternatives occur, the choice may vary from
one member of the sequence to another.

\begin{lemma}\label{old:4.3}
Let $Q_N$ be either $K_{a_N}\vee I_{N-a_N}$ with $a_N/N\to1/3$, or
$K_{b_N}\cup I_{N-b_N}$ with $b_N/N\to2/3$. If
$|E(J_N)\setminus E(Q_N)|=o(N^2)$ and
$\rho(J_N)\ge(2/3-o(1))N$, then
$|E(J_N)\mathbin\triangle E(Q_N)|=o(N^2)$.
\end{lemma}
\begin{proof}
The Frobenius bound shows that removing the exceptional edges changes
the radius by $o(N)$, so we may assume $J_N\subseteq Q_N$. The largest eigenvalue of
$Q_N$ is $(2/3+o(1))N$ and its second largest eigenvalue is at most zero.
For nonnegative unit Perron vectors $x,y$ of $J_N,Q_N$, respectively,
$x^{\mathsf T}A(Q_N)x\ge\rho(J_N)=\rho(Q_N)-o(N)$.
The spectral gap gives $\|x-y\|^2=o(1)$. Since $\|A(J_N)\|\le N$, it
follows that $y^{\mathsf T}(A(Q_N)-A(J_N))y=o(N)$. The explicit Perron
vector of either template has $y_iy_j\ge c/N$ on every template edge,
for an absolute $c>0$ and all sufficiently large $N$. Thus only $o(N^2)$
edges are missing.
\end{proof}

\begin{lemma}\label{old:4.4}
If $\nu^*(J_N)\le(1/3+o(1))N$ and
$\rho(J_N)\ge(2/3-o(1))N$, then $J_N$ differs by $o(N^2)$ edges from
one of the two templates in Lemma~\ref{old:4.3}.
\end{lemma}
\begin{proof}
A graph has a minimum fractional vertex cover that is half-integral.
Indeed, choose an extreme optimum with coordinates in $[0,1]$ and consider
the graph of tight edge constraints on its nonintegral coordinates.
A bipartite component permits a sufficiently small alternating
perturbation preserving feasibility, contrary to extremality; otherwise
an odd cycle forces all coordinates in the component to be $1/2$.
Let $A,B,C$ be the classes of weights $1,1/2,0$, of sizes $a,b,c$.
Then $J_N\subseteq K_a\vee(K_b\cup I_c)$ and
$a+b/2\le N/3+o(N)$.

Pass to a subsequence on which the normalized sizes converge, with
$a/N\to t$. Enlarging the limiting $B$-part if necessary gives a
dominating template with proportions $(t,2/3-2t,1/3+t)$,
$0\le t\le1/3$. Suppose $0<t<1/3$, and put $q=2/3-2t$. The normalized
quotient matrix of the enlarged graph converges to
\[
 M_t=
 \begin{pmatrix}
 t&q&1/3+t\\
 t&q&0\\
 t&0&0
 \end{pmatrix}.
\]
Since the enlarged graph contains $J_N$, the spectral hypothesis
gives $\rho(M_t)\ge2/3$. On the other hand, the Perron equations
show that $\rho(M_t)>q$ is the unique zero on $(q,\infty)$ of
$F_t(z)=z-t-tq/(z-q)-t(1/3+t)/z$.
Since $F_t'(z)>0$ and
$F_t(2/3)=\frac32(1/3-t)(t+2/3)>0$,
we have $\rho(M_t)<2/3$, a contradiction.
Thus $t=0$ or $t=1/3$.

If $t=0$, deleting the edges incident with the $o(N)$ vertices
of $A$ from the enlarged graph leaves a clique of order
$(2/3+o(1))N$ and isolated vertices.
If $t=1/3$, its enlarged $B$-part has $o(N)$ vertices,
and deleting the edges within this part leaves a split graph
with clique order $(1/3+o(1))N$.
In either case, $J_N$ has only $o(N^2)$ edges outside the
corresponding template, so Lemma~\ref{old:4.3} gives the required
$o(N^2)$ edit bound.
Applying this argument to convergent subsequences of the part
ratios proves the assertion for the original sequence.
\end{proof}

\begin{lemma}\label{old:4.5}
Let $J_N,J'_N$ be edge-disjoint graphs on the same $N$ vertices, both with
radius at least $(2/3-o(1))N$. There is a partition $S\cup T$, with
$|S|=(1/3+o(1))N$, such that, in one of the two orders, $J_N,J'_N$
differ by $o(N^2)$ edges from $K_S\vee I_T$ and $I_S\cup K_T$.
\end{lemma}
\begin{proof}
We give a self-contained proof inspired by the deficit-vector argument of
Agarwal et al.~\cite[Lemma~12]{Agarwal}. Put $A=A(J_N)$ and
$B=\mathbf J-A$, where $\mathbf J$ is the all-ones matrix, assigning
all loops to $B$. Write $\lambda=\rho(A)$, $\mu=\rho(B)$, and scale
nonnegative Perron vectors so that $\sum_i x_i=\lambda$ and
$\sum_i y_i=\mu$. Then $0\le x_i,y_i\le1$, and
\[
 \begin{aligned}
 (\lambda+\mu)\sum_i x_iy_i
 &=\lambda x^\top y+\mu x^\top y=(Ax)^\top y+x^\top By
 =x^\top(A+B)y=\lambda\mu.
 \end{aligned}
\]
Define
\[
 D_N:=\sum_i(1-x_i)(1-y_i)
     =N-\lambda-\mu+\frac{\lambda\mu}{\lambda+\mu}.
\]
Since $D_N\ge0$ and $\lambda\mu\le(\lambda+\mu)^2/4$, we have
$\lambda+\mu\le4N/3$. The adjacency matrix of $J'_N$ is entrywise
at most $B$, so both $\lambda$ and $\mu$ are $(2/3+o(1))N$,
and hence $D_N=o(N)$.

Choose $t=t_N\to0$ slowly enough that $D_N=o(t^2N)$.
The sets $X_0=\{i:x_i\ge1-t\}$ and $Y_0=\{i:y_i\ge1-t\}$
cover all but $D_N/t^2=o(N)$ vertices. Their intersection has size $o(N)$:
for $i\in X_0\cap Y_0$,
$(Bx)_i+(Ay)_i=\lambda(1-x_i)+\mu(1-y_i)\le2Nt$,
while each vertex of the intersection contributes at least $1-t$
to the left side. Put $X=X_0$, $Y=Y_0\setminus X_0$, and let $Z$
be the remaining $o(N)$ vertices. Since $(Bx)_i=\lambda(1-x_i)$ and $(Ay)_i=\mu(1-y_i)$,
the definitions of $X$ and $Y$ give
\[
 \sum_{i\in X}(Bx)_i+\sum_{j\in Y}(Ay)_j
 \le t(\lambda|X|+\mu|Y|)=o(N^2).
\]
Each edge counted by $e_B(X)$ contributes at least $2(1-t)$
to the first sum, since both endpoint coordinates of $x$ are
at least $1-t$. The same argument applies to $e_A(Y)$ and
the second sum. Thus $e_B(X)=o(N^2)$ and $e_A(Y)=o(N^2)$.
Moreover, for $i\in X$ and $j\in Y$, the relations
$A_{ij}+B_{ij}=1$ and $0\le x_j,y_i\le1$ imply
$x_jy_i\le A_{ij}y_i+B_{ij}x_j$. Summing over these pairs gives
\[
 \begin{aligned}
 \left(\sum_{i\in X}y_i\right)
 \left(\sum_{j\in Y}x_j\right)
 &\le \sum_{i\in X,\,j\in Y}
       (A_{ij}y_i+B_{ij}x_j)\le \sum_{j\in Y}(Ay)_j+\sum_{i\in X}(Bx)_i
 =o(N^2).
 \end{aligned}
\]
Since $x_i\ge1-t$ on $X$ and $y_i\ge1-t$ on $Y$, we have
$\sum_{i\in X}x_i=|X|+o(N)$ and
$\sum_{j\in Y}y_j=|Y|+o(N)$.
Also, $\sum_{i\in Z}(x_i+y_i)=o(N)$, since $|Z|=o(N)$.
Thus the two factors satisfy
\[
 \sum_{i\in X}y_i+\sum_{j\in Y}x_j
 =\lambda+\mu-|X|-|Y|+o(N)
 =(1/3+o(1))N.
\] Since the two nonnegative factors have product $o(N^2)$
and sum $(1/3+o(1))N$, their minimum is $o(N)$.
Interchanging $A,B$ and the corresponding eigenvalues,
vectors and sets if necessary, we may therefore assume that
$\sum_Xy_i=o(N)$; then $\sum_Yx_i=(1/3+o(1))N$.
Now
$\mu=\sum_Xy_i+\sum_Yy_i+\sum_Zy_i=|Y|+o(N)$,
so $|Y|=(2/3+o(1))N$ and $|X|=(1/3+o(1))N$. Finally,
\[
 \lambda\sum_{j\in Y}x_j
 \le e_A(X,Y)+2e_A(Y)+|Z|N+O(N)
\]
forces $e_A(X,Y)=|X||Y|-o(N^2)$. Thus, in one of the two orders,
$A$ and $B$ are close to the split and clique templates.
Returning to the original notation, $J_N$ is close to one
template and its complement to the other. Since $J'_N$ is a
subgraph of this complement, Lemma~\ref{old:4.3} applies;
loops affect only $O(N)$ entries.
\end{proof}

\begin{lemma}\label{old:4.6}
Let $F_n$ be $3$-graphs with $\sigma(F_n)\ge(2/3-o(1))n$ and no
perfect fractional matching. Either there is a set $A$ of size
$(1/3+o(1))n$ with $e(F_n[V\setminus A])=o(n^3)$, or there is a set
$C$ of size $(2/3+o(1))n$ such that only $o(n^3)$ edges contain fewer
than two vertices of $C$. Moreover, in the first case all but $o(n^3)$
triples of type $ABB$ belong to $F_n$, where $B=V\setminus A$; in
the second case the same holds for type $CCD$, where $D=V\setminus C$.
\end{lemma}
\begin{proof}
Take a nonnegative fractional cover $w$ with sum less than $n/3$, and
choose $v$ minimizing $w$, with weight $\alpha<1/3$. Replace $w_i$ by
$z_i=(w_i-\alpha)/(1-3\alpha)$ and cap values at one. This remains
a fractional cover, has $z_v=0$, and satisfies $0\le z_i\le1$ and
$\sum_i z_i<n/3$. In particular, it is a fractional vertex cover of $L_{F_n}(v)$, so
Lemma~\ref{old:4.4} applies to that link.

If the link is close to the clique on $C$, where
$|C|=(2/3+o(1))n$, then for every fixed $\delta>0$ only $o(n)$
vertices of $C$ have $z_i<1/2-\delta$, since all pairs among them
would be missing clique edges. It follows that $r:=\sum_{i\in C}(1/2-z_i)_+=o(n)$,
where $a_+=\max\{a,0\}$.
The total weight bound gives
\[
 \sum_{i\in C}(z_i-1/2)_+
 +\sum_{i\notin C}z_i
 \le n/3-|C|/2+r=o(n).
\]
Hence $\sum_{i\in C}|z_i-1/2|=o(n)$ and
$\sum_{i\notin C}z_i=o(n)$.
Except at $o(n)$ vertices, the weights are close to $1/2$ on $C$
and zero outside. Every edge with at most one vertex in $C$ must
therefore meet this exceptional set, yielding $o(n^3)$ such edges.

If the link is close to the split graph with core $A$, summing
$z_i+z_j\ge1$ over the almost-complete $A$--$B$ cut gives
$|B|\sum_Az_i+|A|\sum_Bz_i\ge|A||B|-o(n^2)$.
Put $u=\sum_A(1-z_i)$ and $r=\sum_Bz_i$.
The total weight bound gives $r-u<n/3-|A|=o(n)$,
while the preceding inequality gives
$|B|u\le |A|r+o(n^2)$.
Consequently $(|B|-|A|)u\le o(n^2)$.
Since $|B|-|A|=(1/3+o(1))n$, we obtain
$u=o(n)$ and $r=o(n)$. Every edge inside $B$
again meets an $o(n)$-vertex exceptional set.

For the saturation assertions, choose a set $E$ of $o(n)$ vertices
such that outside $E$ the weights differ from their respective
limiting values by at most $1/10$.
In the first case, for $b\in B\setminus E$, every edge of
$L_{F_n}(b)$ outside $K_A\vee I_{B\setminus\{b\}}$ meets $E$:
otherwise the corresponding triple has total weight at most $3/10$.
In the second case, for $d\in D\setminus E$, every link edge
outside the clique on $C$ meets $E$, since otherwise its triple
has total weight at most $4/5$.
Thus each link has only $o(n^2)$ edges outside its template,
uniformly over these roots. The spectral hypothesis and
Lemma~\ref{old:4.3} imply that only $o(n^2)$ template edges
are missing, again uniformly.
Summing over the roots counts each missing $ABB$ triple
avoiding $E$ twice, and each missing $CCD$ triple avoiding $E$
once. The triples meeting $E$ number at most
$|E|n^2=o(n^3)$, proving both assertions.
\end{proof}

\begin{theorem}\label{old:4.7}
Suppose $\sigma(H_n)\ge(2/3-o(1))n$. For each vertex $v$, choose a maximum-spectral component $C_v$
of $L_H(v)$ and consider the triples $\{v,x,y\}$ with
$xy\in E(C_v)$. Either all these triples lie in one spanning tight
component having a perfect fractional matching and a closed tight
walk of order $1$ modulo three, or one of the following holds:
\begin{enumerate}
\item $|A|=(1/3+o(1))n$ and $e(H_n[V\setminus A])=o(n^3)$;
\item $|A|=(2/3+o(1))n$ and only $o(n^3)$ edges contain exactly one
      vertex of $A$;
\item $|A|=(1/3+o(1))n$ and only $o(n^3)$ edges contain exactly one
      vertex of $A$.
\end{enumerate}
The errors are uniform functions of the spectral error and the order.
Equivalently, for every prescribed $\eta>0$ there are $\varepsilon>0$
and $n_0$ such that the finite statement holds under
$\sigma(H)\ge(2/3-\varepsilon)n$, with size errors at most $\eta n$
and violation counts at most $\eta n^3$.
\end{theorem}
\begin{proof}
Write $H=H_n$. Since
$\rho(C_v)=\rho(L_H(v))\ge\sigma(H)$ and
$\rho(C_v)\le |V(C_v)|-1$, each $C_v$ has at least
$(2/3-o(1))n$ vertices. It is unique for large $n$, since
two such components cannot be vertex-disjoint.
The triples $\{v,x,y\}$ with $xy\in E(C_v)$ all belong
to one tight component, denoted by $D_v$: indeed, $C_v$
is connected, and two of its edges sharing a vertex
give triples sharing that vertex and $v$.
For matrix comparisons, pad each $C_v$ with zero rows
and columns indexed by $V(H)\setminus V(C_v)$; the notation
$V(C_v)$ continues to denote the original component.
If $D_v\ne D_w$, then $C_v$ and $C_w$ are edge-disjoint:
a shared edge $xy$ would give triples $\{v,x,y\}$ and
$\{w,x,y\}$ in the same tight component.
Thus representatives from distinct $D_v$ give pairwise
edge-disjoint graphs on the common padded vertex set.
Using $\rho(G)^2\le2e(G)$, if there are $k$ representatives,
then $k\sigma(H)^2\le n(n-1)$, which forces $k\le2$.

Suppose $D_v\ne D_w$. Lemma~\ref{old:4.5} gives a partition
$S,T$ of sizes $(1/3+o(1))n,(2/3+o(1))n$, with $C_v$
close to $K_S\vee I_T$ and $C_w$ close to $I_S\cup K_T$,
after interchanging the labels.
For all but $o(n)$ vertices $t\in T$,
$d_{C_w}(t)\ge(2/3-o(1))n$.
The pairs $wx$, with $x\in N_{C_w}(t)$, form a star in
$L_H(t)$, and the corresponding triples $\{t,w,x\}$
belong to $D_w$.
This star meets $V(C_t)$, since their sizes sum to more
than $n-1$, and hence lies in $C_t$.
These triples therefore also belong to $D_t$, giving
$D_t=D_w$.

For each such $t$, the graph $C_t$ is edge-disjoint from
$C_v$. Every edge of $C_t$ outside the clique on $T$ is
therefore a missing edge of $C_v$ in the split template,
so there are only $o(n^2)$ such edges.
Lemma~\ref{old:4.3} makes $C_t$ close to the clique on
this same $T$, uniformly in $t$.
Put $R_t=T\setminus V(C_t)$. At least
$|R_t|(|T|-1)/2$ clique edges are missing from $C_t$,
so $|R_t|=o(n)$.
Every $S$--$T$ edge in another component of $L_H(t)$
has its $T$-endpoint in $R_t$, since distinct components
are vertex-disjoint. Hence the whole link has only
$o(n^2)$ edges between $S$ and $T$, uniformly over these
roots. Summing over $t\in T$ counts each $STT$ triple
twice; the $o(n)$ exceptional roots contribute only
$o(n^3)$. This proves alternative~(3).

Now suppose all $D_v$ coincide with a tight component $F$.
Since $C_v$ has an edge for every $v$, the component $F$
contains every vertex and is therefore spanning.
Each $L_F(v)$ contains $C_v$, so $\sigma(F)\ge\sigma(H)$.
If $F$ has no perfect fractional matching, apply
Lemma~\ref{old:4.6}.
The pair-shadow of $F$ consists of the pairs contained
in its edges. Every edge of $H\setminus F$ has all its
pairs outside this shadow: otherwise it would share
two vertices with an edge of $F$ and belong to the same
tight component.

In the first case of Lemma~\ref{old:4.6}, we have
$e(F[B])=o(n^3)$, and all but $o(n^3)$ triples of type
$ABB$ belong to $F$.
Each pair inside $B$ absent from the shadow accounts
for $|A|$ missing $ABB$ triples, and each such triple
contains only one pair inside $B$.
Thus only $o(n^2)$ pairs inside $B$ are absent from
the shadow. Every additional edge inside $B$ contains
such a pair, and each pair lies in at most $n$ triples.
Consequently $e(H[B])=o(n^3)$, proving alternative~(1).

In the second case, all but $o(n^3)$ triples of type
$CCD$ belong to $F$, and $F$ has only $o(n^3)$ edges
of type $CDD$.
Each cross pair absent from the shadow accounts for
$|C|-1$ missing $CCD$ triples, and each such triple
is counted at most twice. Hence only $o(n^2)$ cross
pairs between $C,D$ are absent from the shadow.
Every additional $CDD$ edge contains such a pair,
so there are only $o(n^3)$ further edges of this type.
This proves alternative~(2) with $A=C$.

It remains to find the required walk when $F$ has a
perfect fractional matching.
Each $C_v$ is nonbipartite, since a bipartite graph
on at most $n$ vertices has radius at most
$\sqrt{e}\le n/2$.
Direct an arc $v\to u$ whenever $u\in V(C_v)$.
There are at least $(2/3-o(1))n^2>\binom n2$ arcs,
so some pair $u,v$ satisfies both $v\to u$ and $u\to v$.
Choose an odd closed graph walk
$w_0,w_1,\ldots,w_{2r+1}$ in $C_v$, with
$w_0=w_{2r+1}=u$; such a walk is obtained by travelling
to an odd cycle and back. Then
\[
 v,w_0,w_1,v,w_2,w_3,v,\ldots,v,w_{2r},w_{2r+1},v
\]
is a tight walk in $F$: every three consecutive positions
consist of $v$ and the endpoints of an edge of $C_v$.
It runs from $(v,u)$ to $(u,v)$ and has $3r+4$ positions,
hence $3r+2\equiv2\pmod3$ ordered-pair transitions.
The same construction in $C_u$, based at $v$, gives
a return walk in $F$ from $(u,v)$ to $(v,u)$, again
with a number of transitions congruent to two modulo
three. Their concatenation is therefore a closed
tight walk of order one modulo three.
\end{proof}

\subsection{The nonextremal reduction}

We use the tight Hamilton framework of Lang and Sanhueza-Matamala~\cite{LangSanhuezaDegree} in the form supplied by their bandwidth theorem~\cite{LangSanhueza}.

\begin{lemma}\label{old:4.8}
For each sufficiently large fixed $s$ and all sufficiently large $n$,
the following holds. Let $\mathcal P$ be a family of $s$-vertex
$3$-graphs. Suppose each $G\in\mathcal P$ has a spanning tight
component $F(G)$ with a perfect fractional matching and a closed
tight walk of order $1$ modulo $3$. Suppose also that
$F(J-x)\cup F(J-y)$ is tightly connected as a subgraph of $J$ whenever $x\ne y$ and
both deletions of an $(s+1)$-vertex $3$-graph $J$ belong to
$\mathcal P$. If every six-set of an $n$-vertex $3$-graph $H$
extends to at least $(1-s^{-2})\binom{n-6}{s-6}$ sets $S$ with
$H[S]\in\mathcal P$, then $H$ has a tight Hamilton cycle.
\end{lemma}
\begin{proof}
Apply Lang and Sanhueza-Matamala~\cite[Theorem~3.5]{LangSanhueza}
with $k=t=3$. The spanning component, fractional matching, closed
walk and deletion consistency are precisely the four conditions in
their Definition~3.1. Let $R$ be the $s$-graph on $V(H)$ with
$E(R)=\{S\in\binom{V(H)}s:H[S]\in\mathcal P\}$.
For every $Q\in\binom{V(H)}6$, the hypothesis gives
\[
 d_R(Q)=|\{S\in E(R):Q\subseteq S\}|
 \ge(1-s^{-2})\binom{n-6}{s-6}.
\]
Since $K_3(H)=H$ and $2k=6$, this is precisely the
robust satisfaction condition in
\cite[Definitions~3.3--3.4]{LangSanhueza}. The target in Theorem~3.5 may be
any $n$-vertex blow-up of a tight $3$-uniform cycle whose clusters have
size at most $(\log\log n)^c$, for the constant $c>0$ supplied there.
Choose the tight cycle on $n$ vertices with all clusters singletons.
It exists for every $n\ge4$, and its cluster sizes satisfy the bound
for large $n$. Thus the theorem imposes no divisibility restriction on $n$.
\end{proof}

\begin{theorem}\label{old:4.9}
For every $\eta>0$, there are $\varepsilon>0$ and $n_0$ such that
every $n$-vertex $3$-graph $H$, $n\ge n_0$, with
$\sigma(H)\ge(2/3-\varepsilon)n$ and $b_i(H)\ge\eta$ for
$i=1,2,3$ has a tight Hamilton cycle.
\end{theorem}
\begin{proof}
Apply Theorem~\ref{old:4.7} with output error $\delta>0$
sufficiently small in terms of $\eta$.
Adjusting a barrier partition to the prescribed size moves
at most $\delta s+1$ vertices, so its corresponding density
is at most
$\delta+(\delta s+1)\binom{s-1}{2}/s^3<\eta/2$
for sufficiently large $s$.
Thus there are $\varepsilon_0>0$ and $s_0$ such that every
$s$-vertex $G$, $s\ge s_0$, satisfying
$\sigma(G)\ge(2/3-\varepsilon_0)s$ and
$b_i(G)\ge\eta/2$ for $i=1,2,3$ has a spanning tight
component $F(G)$ containing the triples corresponding to
all maximum-spectral link components and having a perfect
fractional matching and a closed tight walk of order
one modulo three.
Decrease $\varepsilon_0$ so that
$2/3-\varepsilon_0>1/2$, and put
$\varepsilon=\varepsilon_0/2$.

Choose a sufficiently large fixed $s\ge s_0$ for the
sampling and framework lemmas, and then take $n$ sufficiently
large. Let $\mathcal P$ be the family of $s$-vertex
$3$-graphs satisfying the preceding spectral and density
conditions.
To verify consistency, let $G=J-x$ and $G'=J-y$ belong
to $\mathcal P$, where $J$ has $s+1$ vertices and $x\ne y$.
Fix $v\notin\{x,y\}$.
The maximum-spectral components of $L_G(v)$ and $L_{G'}(v)$
each have more than $s/2$ vertices in the common
$s$-element set $V(J)\setminus\{v\}$, so they share
a vertex $u$.
Choose an edge incident with $u$ in each component.
The corresponding triples belong to $F(G)$ and $F(G')$,
respectively, and both contain $\{v,u\}$.
Hence $F(G)\cup F(G')$ is tightly connected in $J$.

Fix a six-set $Q$ and choose $S$ uniformly among the
$s$-sets containing $Q$.
Apply Lemma~\ref{old:2.4} with
$b=2/3-\varepsilon_0$ and $\gamma=\varepsilon_0/2$.
Since $b+\gamma=2/3-\varepsilon$ and
$b+\gamma/2=2/3-3\varepsilon_0/4>2/3-\varepsilon_0$,
the probability that
$\sigma(H[S])<(2/3-\varepsilon_0)s$
is below $1/(2s^2)$.
Lemma~\ref{old:4.2} gives the same bound for the probability
that $b_i(H[S])<\eta/2$ for some $i\in\{1,2,3\}$.
A union bound therefore gives
$\Pr(H[S]\notin\mathcal P)<s^{-2}$.
Since $Q$ was arbitrary, every six-set has at least
$(1-s^{-2})\binom{n-6}{s-6}$ extensions with
$H[S]\in\mathcal P$.
Lemma~\ref{old:4.8} now gives a tight Hamilton cycle in $H$.
\end{proof}

\begin{corollary}\label{old:4.10}
If $H_n$ has no tight Hamilton cycle and
$\sigma(H_n)\ge(2/3-o(1))n$, then
$\min\{b_1(H_n),\allowbreak b_2(H_n),\allowbreak b_3(H_n)\}=o(1)$.
\end{corollary}
\begin{proof}
Otherwise a subsequence has all three densities at least a fixed
$\eta>0$, and Theorem~\ref{old:4.9} applies to its sufficiently
large members.
\end{proof}

For a partition $V(H)=A\cup B$, use $AAA,AAB,ABB,BBB$ for its
four triple types. A type is saturated up to $o(n^3)$ when all but
$o(n^3)$ triples of that type belong to $H$.

\begin{lemma}\label{old:4.11}
Let $\sigma(H_n)\ge(2/3-o(1))n$. The following statements hold,
with uniform epsilon--delta versions, for partitions $A\cup B$.
\begin{enumerate}
\item If $|A|=(1/3+o(1))n$ and $e(H[B])=o(n^3)$, then $AAB$ and
      $ABB$ are saturated up to $o(n^3)$.
\item If $|A|=(2/3+o(1))n$ and only $o(n^3)$ edges have type
      $ABB$, then $AAB$ is saturated up to $o(n^3)$.
\item If $|A|=(1/3+o(1))n$ and only $o(n^3)$ edges have type
      $ABB$, then $AAA,AAB,BBB$ are all saturated up to $o(n^3)$.
      Thus $H$ differs by $o(n^3)$ edges from the complete
      hypergraph of those three types.
\end{enumerate}
No saturation of $AAA$ is asserted in the first case, nor of
$AAA$ or $BBB$ in the second.
\end{lemma}
\begin{proof}
Write $H=H_n$. Double counting shows that a type with
$o(n^3)$ edges contributes $o(n^2)$ link edges at all but
$o(n)$ roots, with a uniform error bound.
Deleting $m=o(n^2)$ edges changes a link radius by at most
$\sqrt{2m}=o(n)$, by the Frobenius norm bound.

In the first case, for each typical $b\in B$, the link
$L_H(b)$ has only $o(n^2)$ edges outside
$K_A\vee I_{B\setminus\{b\}}$.
Lemma~\ref{old:4.3} implies that only $o(n^2)$ edges of
this split template are missing.
Summing over $B$ proves saturation of $AAB$ and $ABB$,
counted once and twice, respectively.

In the second case, delete the $o(n^2)$ cross edges from
a typical $B$-root link, leaving graphs $J_A$ on $A$ and
$J_B$ on $B\setminus\{b\}$.
The remaining radius is at least $(2/3-o(1))n$,
whereas $\rho(J_B)\le(1/3+o(1))n$.
Hence $\rho(J_A)\ge(2/3-o(1))n$, and
\[
 \binom{|A|}{2}-e(J_A)
 \le \binom{|A|}{2}-\frac{\rho(J_A)^2}{2}
 =o(n^2).
\]
Summing over $B$ proves saturation of $AAB$, with each
missing triple counted once.

In the third case, a typical $A$-root link has only
$o(n^2)$ edges outside
$K_{A\setminus\{a\}}\vee I_B$.
Lemma~\ref{old:4.3} and summation over $A$ give saturation
of $AAA$ and $AAB$, counted three times and twice,
respectively.
At a typical $B$-root, deleting the rare cross edges
leaves two blocks. The $A$-block has radius at most
$(1/3+o(1))n$, so the $B$-block has radius at least
$(2/3-o(1))n$.
The same clique-edge estimate as above gives only
$o(n^2)$ missing pairs in $B\setminus\{b\}$.
Summing over $B$ proves saturation of $BBB$, with each
missing triple counted three times.

In each case, triples meeting exceptional roots contribute
only $o(n)O(n^2)=o(n^3)$ errors.
In the third case, the three saturated types together
with the assumed $o(n^3)$ edges of type $ABB$ give the
asserted edit bound.
\end{proof}

These are average defect bounds. The exact arguments below also account
for the exceptional vertices allowed by the $o(n^3)$ errors.
\section{Tight Hamilton cycles near space barriers}\label{sec:space-twohit}

We first prove a completion lemma used for all three barriers.
It fills the remaining positions after the required paths and
exceptional vertices have been placed.

\begin{lemma}\label{lem:common-completion}
Let $0<\eta\ll\beta<1$ and let $n$ be sufficiently large.
Let $H$ be an $n$-vertex $3$-graph, and let $U,B\subseteq V(H)$
be disjoint, with $|U|,|B|\ge\beta n$.
Consider a linear or cyclic word of $n$ positions: $|U|$
unfilled positions are reserved for $U$, $|B|$ positions are
reserved for $B$, and the remaining positions are filled
bijectively by $V(H)\setminus(U\cup B)$. Some $B$-positions
may also be filled by distinct prescribed vertices of $B$.
Suppose at most $\eta n$ positions are filled, and that:
\begin{enumerate}
\item Each $U$-position has two $B$-positions on either side.
Every three consecutive positions containing no $U$-position
are already filled and form an edge of $H$.
\item A graph $\Gamma$ on $B$ satisfies
$\Delta(\overline\Gamma)\le\eta n$; each edge of $\Gamma$
has at least $|U|-\eta n$ completions in $U$; and each
$u\in U$ misses at most $\eta n^2$ pairs of $B$ in its link.
\item For each local pattern $b_1,b_2,*,b_3,b_4$ around a
$U$-position, call $b_1b_2,b_2b_3,b_3b_4$ its required pairs.
Every required pair whose positions are both filled
belongs to $\Gamma$.
\end{enumerate}
Then the word can be filled to give a spanning tight path or
cycle, respectively, preserving all prescribed entries.
\end{lemma}

\begin{proof}
Call the initially unfilled $B$-positions movable, and
randomly assign the unused $B$-vertices to them.
Let $Y$ count required pairs outside $\Gamma$,
and let $I_u$ count $U$-positions incompatible with $u$.
The unfilled $B$-positions number at least $(\beta-\eta)n$.
The degree bound gives $\mathbb EY=O_\beta(\eta n)$.
Only $O(\eta n)$ $U$-positions have a prescribed $B$-vertex
among their four neighbours; at every other position each
required pair is uniform among the unused $B$-vertices.
Hence $\mathbb EI_u=O_\beta(\eta n)$ as well.
Put $\gamma=\sqrt\eta$. A transposition changes $I_u$ by
at most eight. The permutation exposure martingale and
Azuma's inequality, followed by a union bound, give
$I_u\le\gamma n$ for all $u\in U$ with probability tending
to one. Markov's inequality gives $Y\le\gamma n$ with
positive probability simultaneously.

The graph of required pairs on the $B$-positions has maximum
degree two. A bad pair has a movable endpoint. Swap the
vertex there with one at another movable position, requiring
all pairs incident with the two positions to belong to
$\Gamma$. Each forbidden neighbour occupies one position
and affects at most two candidate swap positions.
The degree bound therefore excludes at most $8\eta n+O(1)$
choices, including nearby positions, so such a swap exists.
It removes a bad pair and creates none. After at most
$\gamma n$ swaps, all required pairs belong to $\Gamma$,
and each $u$ is incompatible with at most $9\gamma n$ slots.
Each slot accepts at least $|U|-3\eta n$ vertices.
For sufficiently small $\eta$, both degrees in the balanced
compatibility graph exceed $|U|/2$. Hall's theorem fills the
$U$-positions, proving the lemma.
\end{proof}

\subsection{The space barrier}
Put $a^*=\lceil n/3\rceil-1$ and $\ell^*=n-3a^*\in\{1,2,3\}$. In this subsection let $\widehat T_n=T_n$, except that $\widehat T_{3m}=S_m$ when $3\nmid m$. Thus $\widehat T_n\le T_n$ by Lemma~\ref{old:2.1}. We prove the stronger bound $\widehat T_n$ near the space barrier, allowing missing triples and arbitrary exceptional vertices.

A \emph{tight forest} is a vertex-disjoint union of tight paths; its size is its number of edges. The boundary pairs of an ordered tight path are its first and last two vertices. Relative to a graph $\Gamma$, a forest is \emph{good} if its components can be ordered with both boundary pairs in $\Gamma$. An edge is \emph{usable} if at least two of its pairs belong to $\Gamma$; these are exactly the good one-edge paths. A \emph{triangle windmill} has a vertex $v$ and disjoint triples $T_i$ partitioning the remaining vertices, with edges $vxy$ for $xy\in\binom{T_i}{2}$.

\begin{lemma}\label{old:5.1}
Let $F$ be a $3$-graph on an odd number $N\ge400$ of vertices, with minimum degree at least two. If $F$ has no three-edge tight forest as a subhypergraph, then it is a triangle windmill. Conversely, every triangle windmill has minimum degree at least two and no such forest.
\end{lemma}
\begin{proof}
We use two elementary facts. A graph containing no
(not necessarily induced) $P_4$ has only stars, triangles
and isolated vertices as components.
Also, a linear intersecting $3$-graph either has a common
vertex or has at most seven edges: four edges through one
vertex force every edge to contain that vertex; otherwise
all degrees are at most three, and counting through a fixed
edge gives at most $1+3(3-1)=7$ edges.

A maximal matching has at most two edges, and its vertices
cover all edges of $F$. Since $e(F)\ge2N/3$, some vertex $v$
has $d_F(v)\ge N/9>44$.
The graph $G=L_F(v)$ is $P_4$-free, since a graph path
$a,b,c,d$ gives the tight path $a,b,v,c,d$.
If every edge contains $v$, minimum degree forces every
component of $G$ to be a triangle, giving a triangle windmill.

Suppose some edge avoids $v$.
At most three components of $G$ contain at least two edges:
an edge avoiding $v$ misses one of any four such components,
whose two adjacent edges give a disjoint two-edge tight path.
No star component has five leaves. Indeed, if its centre is
$u$, an edge avoiding both $u,v$ would miss two leaves,
which give a two-edge tight path disjoint from that edge.
Thus every edge avoiding $v$ contains $u$.
A leaf $x$ has only the edge $vux$ through $v$, so minimum
degree gives a second edge $uxz$.
Choosing another leaf $y$ outside this edge produces the
three-edge tight path $z,x,u,v,y$.

The components with at least two edges therefore contribute
at most twelve edges to $G$, leaving at least $M\ge45-12=33$
isolated edges $x_i y_i$.
Every one of their $2M$ endpoints belongs to an edge of $F-v$.
Two edges of $F-v$ cannot be disjoint or meet in two vertices:
an isolated pair $x_i y_i$ disjoint from their union would
complete a forbidden forest.
Hence $F-v$ is linear intersecting.
Its edges cover more than $21$ vertices, so they have a
common vertex $u$.
Their pairs outside $u$ are disjoint and cover at least
$2M-1$ of the endpoints, so there are at least $M$ such edges.
Every two-edge tight path in $F$ must contain $u$, since
otherwise its four vertices miss one of these edges,
again producing a forbidden forest.

Put $W=V(F)\setminus\{u,v\}$.
The pairs of edges containing $u$ but not $v$ form a matching
$M_u$ on $W$, by linearity of $F-v$.
The analogous pairs through $v$ form a matching $M_v$,
since two meeting pairs would give a two-edge tight path
avoiding $u$.
Thus, for some $S\subseteq W$,
\[
 E(F)=\{uxy:xy\in M_u\}\cup\{vxy:xy\in M_v\}
      \cup\{uvx:x\in S\}.
\]
Minimum degree implies that every vertex of $W\setminus S$
is covered by both matchings, and every vertex of $S$ by
at least one.
If $x\in S$, $xy$ belongs to a matching, and
$z\in S\setminus\{x,y\}$, then $y,x,v,u,z$, or the order
with $u,v$ exchanged, is a three-edge tight path.
Thus $|S|\le2$; when $|S|=2$, every matching edge meeting
$S$ joins its two vertices.
In each case, at least one of the matchings covers $S$
and hence all of $W$.
This forces $|W|$ to be even, contrary to odd $N$.

Conversely, all edges of a triangle windmill contain its
centre, so a tight forest cannot have two nonempty
vertex-disjoint components.
A three-edge tight path would have the centre as its
middle vertex and yield a $P_4$ in its link, which is
a disjoint union of triangles.
Every noncentral vertex has degree two, proving the converse.
\end{proof}

\begin{lemma}\label{old:5.2}
Let $K\ge1$ and let $F$ be a $3$-graph on an odd number $N\ge1000K^3$ of vertices, with minimum degree at least two. For a graph $\Gamma$ on its vertex set put $D=\overline\Gamma$ and $t_v=d_D(v)$. If $\Delta(D)\le N/(12K)$ and $Nt_v\le Kd_F(v)$ for every $v$, then either $F$ has a good three-edge tight forest, or $D$ is empty and $F$ is a triangle windmill.
\end{lemma}
\begin{proof}
Suppose that $F$ contains no good three-edge tight forest
as a subhypergraph, and write $e=e(F)$.
Since $\sum_vt_v\le3Ke/N$, the number $e_0$ of unusable
edges satisfies
\[
 e_0\le\sum_v\binom{t_v}{2}
 \le\frac{\Delta(D)}2\sum_vt_v\le e/8.
\]
A maximal matching of usable edges has at most two edges,
so its vertices cover all usable edges and number at most six.
For each cover vertex $c$, the graph $L_F(c)\cap\Gamma$
is $P_4$-free, since a $P_4$ would give a good three-edge
tight path.
It therefore has at most $N$ edges.
Counting through the cover according to whether the
opposite pair lies in $\Gamma$ or $D$ gives
$e-e_0\le6N+6e(D)\le6N+9Ke/N$.
Since $N\ge72K$, this yields $e\le8N$.
Thus $Q=\{v:t_v>0\}$ has size at most $24K$, and
$d_F(q)\ge N/K$ for every $q\in Q$.
Put $L=V(F)\setminus Q$.
Every pair incident with $L$ is good, so every edge
meeting $L$ is usable.

Every pair $uv$ has at most four leaves in $L$.
Otherwise, choose a leaf $x$ and a second edge $f\ne uvx$
through $x$, which is usable since $x\in L$.
If $f$ avoids $\{u,v\}$, two leaves outside $f$ give
a disjoint good two-edge path through $u,v$.
Otherwise $f=uxz$ or $vxz$, and another leaf $y$ outside
$f$ gives the good three-edge path $z,x,u,v,y$, or the
order with $u,v$ exchanged.
Both cases contradict the assumption.

Write $s=|Q|$.
For distinct $q,r\in Q$, the pair $qr$ has at most four
leaves in $L$ and at most $s-2$ in $Q$.
Hence at most $(s-1)(s+2)\le600K^2$ edges through $q$
meet another vertex of $Q$.
The graph $J_q=L_F(q)[L]$ therefore satisfies
\[
 e(J_q)\ge N/K-600K^2\ge400K^2>20,
 \qquad \Delta(J_q)\le4.
\]
Consequently it has an edge avoiding any prescribed
five vertices.

Suppose $D$ is nonempty and choose $cd\in E(D)$.
Every usable edge meets $\{c,d\}$: otherwise edges
chosen successively in $J_c$ and $J_d$ would complete
a good three-edge matching.
If $q\in Q\setminus\{c,d\}$, any edge of $J_q$ gives
a usable triple avoiding $\{c,d\}$, a contradiction.
Thus $Q=\{c,d\}$.
Every edge now meets $L$ and is usable, so every edge
of $F$ meets $\{c,d\}$.

Every two-edge tight path has a good ordering:
exchanging its middle vertices gives disjoint
boundary-pair sets, at least one of which avoids
the sole bad pair $cd$.
A three-edge path cannot have boundary pair $cd$,
since its opposite edge would avoid $\{c,d\}$.
Together with usability of every edge, this shows
that every three-edge tight forest is good.
Lemma~\ref{old:5.1} therefore makes $F$ a triangle
windmill, contradicting $d_F(c),d_F(d)\ge N/K>2$.
Hence $D$ is empty, and Lemma~\ref{old:5.1} finishes
the proof.
\end{proof}

\begin{lemma}\label{old:5.3}
Let $V(H)=A\mathbin{\dot\cup}B$, $|A|=a\ge3$, $a\le|B|-1\le3a$, and $F=H[B]$. Put $J=K_a\vee I_{|B|-1}$. If $\rho(L_H(v))\ge\rho(J)$ for $v\in B$, the number $M_v$ of missing edges meeting $A$ in that link is at most $C d_F(v)$, where $C=4096$. More generally, if $h\ge0$ and $\rho(L_H(v))\ge\max\{a/2,\rho(J)-h\}$, then $M_v\le16384(ah+d_F(v))$.
\end{lemma}
\begin{proof}
Let $G=L_H(v)$, $\lambda=\rho(G)$ and $\lambda_0=\rho(J)$.
Write $A(G)=A(J)-A(M)+A(T)$, where $M$ consists of the
missing split edges and $e(T)=d=d_F(v)$.
For nonnegative unit Perron vectors $x,p$ of $G,J$,
the eigenvalue equation gives $\max_i x_i^2\le4a/\lambda^2$.
Also $\lambda_0\ge a$, and the clique-to-independent coordinate ratio of $p$
lies in $[1,3]$, so $p_i^2\ge1/(12a)$.
Put $\eta=\lambda_0-x^\top A(J)x\ge0$.
Since all other eigenvalues of $J$ are nonpositive,
$\|x-p\|^2\le2\eta/a$.
At most $96\eta$ vertices have $x_i<p_i/2$, and at most
$384a\eta$ missing edges meet them.
Each remaining missing edge contributes at least
$1/(24a)$ to $x^\top A(M)x$. Hence
\[
 M_v\le384a\bigl(\eta+x^\top A(M)x\bigr)
 =384a\bigl(\lambda_0-\lambda+x^\top A(T)x\bigr).
\]
If $\lambda\ge\lambda_0$, then $\max_i x_i^2\le4/a$,
so $M_v\le3072d\le Cd$.
More generally, $\lambda\ge\max\{a/2,\lambda_0-h\}$
gives $\max_i x_i^2\le16/a$ and
$M_v\le384ah+12288d\le16384(ah+d)$.
\end{proof}

Define $\Gamma$ on $B$ by $d_A(xy)\ge(1-\tau)a$, and put $D=\overline\Gamma$. Lemma~\ref{old:5.3} gives $d_D(v)\le Cd_F(v)/(\tau a)$. If every $v\in B$ belongs to at most $\epsilon ab$ missing $ABB$ triples and $\tau=\epsilon^{1/3}$, then $\Delta(D)\le\epsilon^{2/3}b$, where $b=|B|$. Consequently the number of unusable internal edges is at most
\begin{equation}\label{eq:old:7}
 \sum_v\binom{d_D(v)}2\le\frac{3C\epsilon b}{2\tau^2a}e(F).
\end{equation}
This argument uses only the link bounds at vertices of $B$. The next lemma provides the internal edges needed to overcome
the space barrier, arranged as a good $\ell^*$-edge tight forest
for the Hamilton cycle construction.

\begin{lemma}\label{old:5.4}
Let $|A|=a^*$, $B=V(H)\setminus A$, and suppose every $v\in B$ belongs to at most $\epsilon|A||B|$ missing $ABB$ triples. For sufficiently small fixed $\epsilon>0$ and sufficiently large $n$, the inequality $\sigma(H)>\widehat T_n$ forces a good $\ell^*$-edge forest in $H[B]$, with $\tau=\epsilon^{1/3}$.
\end{lemma}
\begin{proof}
The canonical split radius is at most $\widehat T_n$, so Lemma~\ref{old:5.3} applies. Choose $\epsilon$ so that the coefficient in \eqref{eq:old:7} is at most $1/4$. If $\ell^*=1$, the strict spectral inequality makes $F=H[B]$ nonempty and \eqref{eq:old:7} supplies a usable edge.

If $\ell^*=2$, put $a=a^*$ and $\lambda_0=\rho(K_a\vee I_{2a+1})$. The equation for its radius gives $h=2a-\lambda_0>1/4$, since $h(3a+1-h)=a$. The Perron-coordinate estimate gives $\rho(L_H(v))\le\lambda_0+8d_F(v)/a$, so $d_F(v)>a/32$ and $e(F)>ab/96\ge b^2/288$. If there is no good two-edge forest, the usable family $F^+$ is intersecting and has maximum pair degree at most $\max\{1,\Delta(D)\}$. Indeed, consider usable edges through $xy$. If there are at least two leaves, $xy$ must be good; the absence of a good two-edge path then forces every leaf to have the same singleton good neighbour in $\{x,y\}$, bounding their number by $\Delta(D)$. A fixed usable edge covers $F^+$, giving $e(F^+)\le3(b-1)\max\{1,\Delta(D)\}/2$. By \eqref{eq:old:7}, $e(F)\le2b+2\epsilon^{2/3}b^2$, a contradiction.

Finally, let $\ell^*=3$, so $n=3m$, $a=m-1$ and $b=2m+1$. Since $\widehat T_n\ge S_m$, a vertex of $F$ of degree at most one would have link radius at most $S_m$, which is impossible. For $K=13000/\tau$, the split-defect bound gives $bd_D(v)\le Kd_F(v)$. If $\epsilon$ is small and then $m$ large, $\Delta(D)\le b/(12K)$ and $b\ge1000K^3$. In the absence of a good forest, Lemma~\ref{old:5.2} makes $F$ a triangle windmill. Thus $3\mid m$, and a noncentral vertex has link contained in $K_{m-1}\vee(P_3\cup I_{2m-3})$, giving $\sigma(H)\le W_m=\widehat T_n$, a contradiction.
\end{proof}

\begin{lemma}\label{old:5.5}
There is $\mu_0>0$ such that the following holds for sufficiently large $n$. Suppose $|A|=a=a^*-d\ge n/4$, $d\ge1$, $\Delta_1(H[B])\le\mu n^2$, $0<\mu<\mu_0$, and $\sigma(H)>\widehat T_n$. Put $\epsilon=8C\mu$ and $\tau=\epsilon^{1/3}$. Then $H[B]$ has a usable matching of $n-3a$ edges, and every root in $B$ misses at most $C\mu n^2$ triples of type $ABB$.
\end{lemma}
\begin{proof}
Write $J_j=K_j\vee I_{n-1-j}$ and $F=H[B]$.
By~\eqref{eq:old:4}, the derivative of $\rho(J_j)$
with respect to $j$ is at least $1/2$ for
$n/4\le j\le a^*$.
Since $\widehat T_n\ge\rho(J_{a^*})$, we have
$\widehat T_n-\rho(J_a)\ge d/2$.
Lemma~\ref{old:5.3} gives
$M_v\le Cd_F(v)\le C\mu n^2$.
The same Perron-coordinate estimate gives
$\rho(L_H(v))\le\rho(J_a)+8d_F(v)/a$, whence
$d_F(v)>ad/16$ and
$e(F)>a|B|d/48\ge dn^2/384$.

We have $d_D(v)\le Cd_F(v)/(\tau a)$ and
$\Delta(D)\le4C\mu n/\tau$.
Thus the number of unusable edges is at most
$\sum_v\binom{d_D(v)}2\le24C^2\mu e(F)/\tau^2
\le e(F)/2$, for sufficiently small $\mu_0$.
The vertices of a maximal usable matching cover all
usable edges, so this matching has at least
$e(F)/(6\Delta_1(F))>d/(2304\mu)\ge6d$ edges.
Since $n-3a=3d+\ell^*\le6d$, it contains a matching
of the required size.
\end{proof}

\begin{lemma}\label{old:5.6}
For sufficiently small fixed $\alpha>0$ and sufficiently large $n$, suppose $V(H)=A_0\mathbin{\dot\cup}B_0\mathbin{\dot\cup}Z$, where $|A_0|=a^*-t$ and $0\le t\le q:=|Z|\le\alpha n$. Assume each vertex in $A_0\cup B_0$ belongs to at most $\alpha n^2$ missing triples of type $A_0B_0B_0$, each vertex in $B_0$ has degree at most $\alpha n^2$ in $H[B_0]$, and at most $\alpha n^3$ triples of type $A_0A_0B_0$ are missing. If $\sigma(H)>\widehat T_n$, then $H$ has a tight Hamilton cycle. No assumption is made on the links of $Z$.
\end{lemma}
\begin{proof}
We construct the cycle in the cyclic form
\[
 a_1,G_1,a_2,G_2,\ldots,a_a,G_a,
\]
where $a$ is the number of $A$-positions, the $a_i$ denote these
positions, and each $G_i$ is a nonempty sequence of $B$-vertices,
called a gap. The sets occupying the two types of positions will
be specified in each case. Ordinary gaps have two vertices;
longer gaps and singleton gaps adjust the vertex counts.

Define $\Gamma_0$ on $B_0$ by $d_{A_0}(xy)\ge(1-\alpha^{1/3})|A_0|$. Then $\Delta(\overline\Gamma_0)\le4\alpha^{2/3}n$ and $e(\overline\Gamma_0)\le2\alpha^{2/3}n^2$. Call $z\in Z$ rich if $L_H(z)[B_0]\cap\Gamma_0$ contains at least $30q+100$ vertex-disjoint $P_4$'s. Let $R$ be the rich set and $r=|R|$. For a nonrich $z$, a maximal packing of these paths covers at most $120q+396$ vertices; outside them the graph is $P_4$-free and has at most $n$ edges. With $B=B_0\cup Z$, this yields $d_{H[B]}(z)\le2\alpha^{2/3}n^2+(121q+397)n\le3\alpha^{2/3}n^2$.

Put $\mu=10\alpha^{2/3}$ and $\epsilon=8C\mu$, choosing $\alpha$ small enough for the preceding lemmas. Since $n/4\le|A_0|\le a^*$, Lemma~\ref{old:5.3} applies with clique $A_0$. Each nonrich exceptional vertex misses at most $C\mu n^2$ cross triples at that root. The original roots in $B_0$ have at most $\alpha n^2+qn$ missing cross triples after $Z$ is added, and each $u\in A_0$ misses at most $2\alpha n^2$ pairs outside $A_0\cup R$.

Suppose first $r>t$. Use $A_0$ as the set of $A$-positions in a cyclic word and put $s=3(r-t)-\ell^*\ge0$. For each rich vertex choose a disjoint $P_4$, ordered $p_1p_2p_3p_4$, in its rich link, and prescribe $p_1,p_2,z,p_3,p_4$ as a five-vertex $B$-gap. Its three internal triples are edges and its boundary pairs belong to $\Gamma_0$; the greedy choices avoid fewer than $4q$ previously used vertices of $B_0$.

Also choose $s$ disjoint paths $x_1,x_2,u,y,u',x_3,x_4$, with $u,u'\in A_0$, all other vertices in $B_0$, and boundary pairs in $\Gamma_0$, avoiding earlier choices. These exist greedily: in a uniformly chosen ordered tuple outside the $O(q)$ forbidden vertices, the middle triple $uyu'$ is missing with probability $O(\alpha)$ by the average $A_0A_0B_0$ bound; every other triple is missing with probability $O(\alpha)$ by the pointwise bounds; and each boundary pair is bad with probability $O(\alpha^{2/3})$. The sum of these failure probabilities is less than one. Here $y$ forms a singleton gap, with its two adjacent $A$-positions fixed. Use the $r$ long gaps, the $s$ singleton gaps, and ordinary two-vertex gaps elsewhere. The count is $5r+s+2(|A_0|-r-s)=n-|A_0|$. Keep each singleton together with its two prescribed neighbouring ordinary gaps, and separate the resulting blocks and the rich gaps by movable ordinary gaps. All triples through fixed $A$-positions have been checked. On $B^*=B\setminus R$ define $\Gamma$ by $d_{A_0}(xy)\ge(1-\epsilon^{1/3})|A_0|$. Then $\Gamma_0\subseteq\Gamma[B_0]$ and $\Delta(\overline\Gamma)=O(\epsilon^{2/3}n)$.

Now suppose $r\le t$. Put $d=t-r$, $A'=A_0\cup R$ and $B'=B_0\cup(Z\setminus R)$. The internal 3-graph $F=H[B']$ has maximum degree at most $\mu n^2$: an original root gains at most $qn$ edges, and a nonrich root satisfies the preceding bound. For $d=0$, Lemma~\ref{old:5.3} gives at most $C\mu n^2$ missing cross triples at every root in $B'$, so Lemma~\ref{old:5.4} supplies a good $\ell^*$-edge forest. For $d>0$, Lemma~\ref{old:5.5} supplies a good matching of $3d+\ell^*$ edges. In both cases use $\Gamma$ defined by $d_{A'}(xy)\ge(1-\epsilon^{1/3})|A'|$. The forest has $\ell=n-3|A'|$ edges and at most $9q+9$ vertices; its components are the long gaps.

Choose the forest first. For each $z\in R$, select a $P_4$ from its rich packing avoiding the forest and earlier choices; fewer than $9q+9+4q<30q+100$ vertices are forbidden. Prescribe $p_1,p_2,z,p_3,p_4$, now with $z$ in a fixed $A'$-position and ordinary gaps on both sides. A pair of $\Gamma_0$ misses at most $\alpha^{1/3}|A_0|+r<\epsilon^{1/3}|A'|$ vertices of $A'$, so its boundary pairs are good. If the forest has $k$ components, its gaps use $\ell+2k$ vertices, leaving exactly $2(|A'|-k)$ for ordinary gaps. Separate all prescribed blocks by movable ordinary gaps. There are $O(q+1)$ fixed positions; all unassigned $A'$-vertices lie in $A_0$, and $\Delta(\overline\Gamma)=O(\epsilon^{2/3}n)$. This includes $q=0$.

Apply Lemma~\ref{lem:common-completion} with $U$ the unassigned
$A_0$-vertices and with $B$-set $B\setminus R$ in the first
case, or $B'$ in the second. Both sets have size at least
$n/5$. Only $O(q+1)$ positions are prescribed.
Every free $A$-position has gaps of length at least two on
both sides, and all triples not containing a free position
have already been checked. Required pairs with two fixed
endpoints are the prescribed good boundary pairs.
The chosen graph $\Gamma$ has complementary maximum degree
$O(\epsilon^{2/3}n)$, and each of its edges misses at most
$\epsilon^{1/3}n$ vertices of $U$.
Each $u\in U$ misses $O((\epsilon+\alpha)n^2)$ relevant pairs.
Thus the lemma applies with $\beta=1/5$ and defect parameter
$(\epsilon+\alpha)^{1/4}$, completing the cycle.
\end{proof}

\begin{theorem}\label{old:5.7}
There are $\delta>0$ and $n_0$ such that, for $n\ge n_0$, any $3$-graph $H$ admitting a partition $V(H)=\bar A\mathbin{\dot\cup}\bar B$ with $|\bar A|=\lfloor n/3\rfloor$ and $e(H[\bar B])\le\delta n^3$ has a tight Hamilton cycle whenever $\sigma(H)>\widehat T_n$. In particular, $\sigma(H)>T_n$ is sufficient.
\end{theorem}
\begin{proof}
Fix $\alpha$ for Lemma~\ref{old:5.6}, then
$\zeta\ll\alpha^2$. Since $\widehat T_n=(2/3+o(1))n$,
Lemma~\ref{old:4.11} permits $\delta$ small enough and $n_0$
large enough that $H[\bar B]$ has at most $\zeta n^3$ edges
and each of the $\bar A\bar A\bar B$ and
$\bar A\bar B\bar B$ layers misses at most $\zeta n^3$ triples.

Mark roots with more than $\alpha n^2$ missing
$\bar A\bar B\bar B$ triples, and roots in $\bar B$ with
internal degree greater than $\alpha n^2$.
Each missing $\bar A\bar B\bar B$ triple is counted at three
roots, so the first rule marks at most $3\zeta n/\alpha$
vertices. Similarly,
$\sum_{v\in\bar B}d_{H[\bar B]}(v)=3e(H[\bar B])
\le3\zeta n^3$, so the second rule marks at most
$3\zeta n/\alpha$ vertices. Thus there are at most
$6\zeta n/\alpha$ marked vertices.

Put them in $Z$; if necessary move one additional vertex
of $\bar A$ to $Z$ so that $A_0=\bar A\setminus Z$ has size
at most $a^*$. With $B_0=\bar B\setminus Z$ and
$t=a^*-|A_0|$, we have $0\le t\le|Z|\le\alpha n$ for large
$n$. The marking rule and the average missing-layer bound
verify the hypotheses of Lemma~\ref{old:5.6}.
\end{proof}

\begin{corollary}\label{old:5.8}
Under the space hypothesis of Theorem~\ref{old:5.7}, a non-Hamiltonian $3$-graph on $3m$ vertices satisfies $\sigma(H)\le W_m$ if $3\mid m$, and $\sigma(H)\le S_m$ otherwise. Both bounds are attained.
\end{corollary}
\begin{proof}
Only attainment remains. Take $|A|=m-1$, $|B|=2m+1$, and include every triple meeting $A$. A Hamilton cycle would require at least three internal edges forming a tight forest: if its $B$-runs have lengths $g_i$, their total number of internal edges is $\sum_i\max\{g_i-2,0\}\ge|B|-2|A|=3$. When $3\mid m$, put a triangle windmill on $B$; it has no such forest and its noncentral links have radius $W_m$, as in Proposition~\ref{prop:sharp}. Otherwise put exactly the triples $pqx$ inside $B$, for a fixed pair $pq$. They have no three-edge tight forest. Every ordinary $B$-link is $K_{m-1}\vee(K_2\cup I_{2m-2})$, and all other links have at least this radius, so $\sigma(H)=S_m$. Both constructions have only $O(m)$ internal edges outside an appropriate $m$-set.
\end{proof}

\subsection{The two-hit barrier}
For a partition $V(H)=C\mathbin{\dot\cup}D$ with
$0<|C|<2n/3$, the two-hit barrier consists of all triples
of types $CCC,CCD,DDD$, with the $CDD$ layer absent.
Since a $DDD$ edge shares at most one vertex with an edge
of either other type, every edge of a tight Hamilton cycle
would meet $C$ at least twice, forcing $3|C|\ge2n$,
a contradiction. In this subsection we study partitions
with $|C|$ close to $2n/3$. Let $M$ count missing $CCD$
triples and let $F$ be the family of present $CDD$ triples.
We call $CCD$ triples normal and $CDD$ triples rare,
and regard the $C$ vertex of a rare triple as its colour.
No completeness assumption is made on $H[C]$ or $H[D]$. Lemma~\ref{old:4.11}(2) describes the corresponding
asymptotic regime, with $M+|F|=o(n^3)$.

\begin{theorem}\label{old:6.1}
There is $\delta>0$ such that, for all sufficiently large $n=3m+r$, $r\in\{0,1,2\}$, the following holds. Suppose $H$ has a partition $C\mathbin{\dot\cup}D$ with $\bigl||C|-2n/3\bigr|\le\delta n$ and $M+|F|\le\delta n^3$. If $\sigma(H)>U_m$ for $r=0$, $\sigma(H)>V_m$ for $r=1$, or $\sigma(H)>2m$ for $r=2$, then $H$ has a tight Hamilton cycle.
\end{theorem}

We use the balance $b(P)=2|V(P)\cap D|-|V(P)\cap C|$.
We seek a short path with suitable $CC$ end pairs and balance
$2|D|-|C|-2$, which can be completed using normal triples
under the hypotheses of Lemma~\ref{old:6.3}.
The short-path constructions below use rare edges to raise
the balance and $CCC$ triples to lower it.
We first choose a partition in which every vertex has a dense
link in its assigned role. The full template distance counts
missing $CCC,CCD,DDD$ triples and present $CDD$ triples.
We find a local minimum of this distance while controlling
the layer error $M+|F|$.

\begin{lemma}\label{old:6.2}
Let $0<\eta\ll\epsilon\ll\alpha\ll1$. A partition $(C_0,D_0)$ with $\bigl||C_0|-2n/3\bigr|\le\eta n$ and total normal and rare error at most $\eta n^3$ can be changed on at most $3\eta n/\epsilon$ vertices to a local minimum $(C,D)$ of the full template distance. Its size error is at most $(\eta+3\eta/\epsilon)n$, and its layer error is at most $(\eta+3\eta/(2\epsilon))n^3$. For any local minimum, writing $c=|C|$, $d=|D|$, one has
\[
 e(L_H(x)[C\setminus\{x\},D])\ge\tfrac12\left((c-1)d-\binom d2\right)\quad(x\in C),
\]
and $e(L_H(x)[C])\ge e(L_H(x))-c(d-1)/2-3\binom{d-1}{2}/2$ for $x\in D$. If $|c-2n/3|\le\eta n$ and $\sigma(H)\ge(2/3-\eta)n$, both relevant link graphs have at least $n^2/50$ edges for sufficiently small $\eta$ and large $n$.
\end{lemma}
\begin{proof}
For a vertex $x$, let $e_{CC}$, $e_{CD}$ and $e_{DD}$
denote the numbers of edges of $L_H(x)$ with both endpoints
in $C$, one in each part, and both in $D$, respectively.
For $x\in C$, moving $x$ to $D$ changes the distance by
$2e_{CD}-(c-1)d+\binom d2-2e_{DD}$.
Nonnegativity gives the first bound.
For $x\in D$, the reverse move gives
$e_{CD}-e_{DD}\le[c(d-1)-\binom{d-1}{2}]/2$.
Subtracting these link parts and using
$e_{DD}\le\binom{d-1}{2}$ gives the second bound.
Together with $2e(G)\ge\rho(G)(\rho(G)+1)$, these yield
limiting lower bounds $n^2/12$ and $n^2/36$,
respectively, and hence the stated $n^2/50$ bound.

To construct the partition, put $c_0=|C_0|$ and $d_0=|D_0|$.
Mark every vertex incident with at least $\epsilon n^2$
original layer errors. There are at most
$q=3\eta n/\epsilon$ marked vertices, since each error is
counted at three vertices.
The initial move costs at unmarked vertices are at least
$(c_0-1)d_0+\binom{d_0}{2}-2\epsilon n^2$ in $C_0$ and
$c_0(d_0-1)-\binom{d_0-1}{2}-2\epsilon n^2$ in $D_0$.
Both exceed $n^2/7$ for sufficiently small parameters
and large $n$.

Repeatedly move a vertex whenever this decreases the distance.
Inductively, only marked vertices move. Indeed, as long as
this holds, the current partition differs from the original
one on at most $q$ vertices. Changing one other vertex's
class affects at most $n-2$ triples in an unmarked vertex's
move cost, each by at most two. Its current move cost
therefore differs from its initial value by at most
$2qn<n^2/7$, and remains positive.
The distance is a nonnegative integer and strictly decreases,
so the process ends at a local minimum.

At most $q$ vertices have different classes in the final
and original partitions. Thus the size error increases
by at most $q$. Only triples meeting these vertices can
change their layer-error status, so the layer error
increases by at most $q\binom{n-1}{2}$.
These give the claimed bounds.
\end{proof}

Fix $0<\eta\ll\epsilon\ll\theta\ll1$, and assume
$|C|,|D|\ge n/4$.
When $M\le\eta n^3$, let $X$ consist of the vertices belonging
to at least $\epsilon n^2$ missing normal triples, so
$|X|\le3\eta n/\epsilon$. Vertices outside $X$ are ordinary.
A $CC$ pair is strong if it has at least $(1-\theta)|D|$
completions in $D$; a $CD$ pair is strong if it has at least
$(1-\theta)|C|$ completions in $C$.
For large $n$, each nonstrong pair has at least $\theta n/8$
missing normal completions. Double counting therefore gives
$O(\eta n^2/\theta)$ nonstrong pairs, with
$O(\epsilon n/\theta)$ incident with each ordinary vertex.

\begin{lemma}\label{old:6.3}
For a sufficiently small hierarchy as above and sufficiently
large $n$, let $(C,D)$ be locally optimal,
$\bigl||C|-2n/3\bigr|\le\eta n$, $M\le\eta n^3$, and
$\sigma(H)\ge(2/3-\eta)n$.
If a tight path $Q$ on at most $\eta n/\epsilon$ vertices
has disjoint strong ordinary $CC$ end pairs and
\begin{equation}\label{eq:old:8}
 b(Q)=2|D|-|C|-2,
\end{equation}
then $H$ has a tight Hamilton cycle.
\end{lemma}

\begin{proof}
Any two disjoint ordered strong ordinary pairs of types
$CC,CD,DC$ can be joined by a normal path with four to six
new ordinary vertices, avoiding any prescribed set of at most
$\theta n$ vertices disjoint from the ends.
Give $CC,CD,DC$ phases $0,1,2$ in the periodic word
$CCDCCD\cdots$.
If the end phases are $i,j$, choose $k\in\{4,5,6\}$ with
$k+2\equiv j-i\pmod3$, and choose the $k$ intermediate
vertices independently and uniformly from their full
prescribed classes. Each class has order at least $n/4$.
The two end windows fail with probability at most $\theta$
each; at most two further windows contain one fixed ordinary
vertex and fail with probability $O(\epsilon)$ each;
at most four wholly random windows fail with probability
$O(\eta)$ each.
Collisions, the forbidden set, and $X$ contribute
$O(\theta+\eta/\epsilon)+O(1/n)$.
The total is at most
$2\theta+C(\epsilon+\eta+\theta+\eta/\epsilon)+O(1/n)<1$
for an absolute $C$.
For two $CC$ ends, $k=4$ gives the pattern $DCCD$,
whose balance is two.

Cover each $x\in X\setminus V(Q)$ by a normal five-vertex
path avoiding $Q$, every other exceptional vertex and all
earlier paths.
The role link from Lemma~\ref{old:6.2} has at least
$n^2/50$ edges.
Deleting forbidden vertices and nonstrong pairs costs
$O(\eta n^2/\epsilon+\eta n^2/\theta)$ edges, leaving
more than $n$.
A $P_4$-free graph has at most as many edges as vertices,
so a remaining $P_4$ gives $d_1,c_1,x,d_2,c_2$ when
$x\in C$, and $c_1,c_2,x,c_3,c_4$ when $x\in D$,
with strong ordinary end pairs.
In the first case extend to
$a,b,d_1,c_1,x,d_2,c_2,c$, choosing fresh ordinary
$a,b,c\in C$ outside the forbidden set:
choose $b$ completing $d_1c_1$ with $bd_1$ strong,
$a$ completing $bd_1$ with $ab$ strong, and
$c$ completing $d_2c_2$ with $c_2c$ strong.
Strong neighbourhoods and the ordinary-vertex bounds
permit these choices.
All resulting paths have strong ordinary $CC$ ends and
balance $-2$.
Join them successively to $Q$ by $DCCD$ connectors.
The resulting path $P$ covers $X$, uses
$O(\eta n/\epsilon)$ vertices, and still satisfies
\eqref{eq:old:8}.

Put $A=D\setminus V(P)$, $B=C\setminus V(P)$, and write
$pq,uv$ for the final and initial ordered pairs of $P$.
Then $|B|=2|A|-2$.
Let $W=B\cup\{p,q,u,v\}$ and $a=|A|$, so $|W|=2a+2$.
Every vertex of $A\cup W$ is ordinary and misses at most
$\xi a^2$ normal triples, where $\xi=O(\epsilon)$.
Let $\Gamma$ consist of the pairs on $W$ with at least
$(1-2\theta)a$ completions in $A$. Then
$\Delta(\overline\Gamma)\le\xi a/(2\theta)$, and
$pq,uv\in E(\Gamma)$.
Use a linear word beginning with $p,q$ and ending with $u,v$,
with $a$ unfilled $A$-positions separated by pairs of
$W$-positions. Only the four end vertices are prescribed.
Every triple contains an $A$-position, and the only required
pairs with two prescribed endpoints are $pq$ and $uv$.
Since $\epsilon\ll\theta\ll1$, Lemma~\ref{lem:common-completion}
applies to $H[A\cup W]$, with $U=A$, $B=W$ and $\beta=1/5$.
It gives a spanning normal path with these ordered ends.
Joining it to $P$ along the common end pairs gives a tight
Hamilton cycle.
\end{proof}

\begin{lemma}\label{old:6.4}
Suppose $\bigl||C|-2n/3\bigr|\le\eta n$ and $M\le\eta n^3$. Choose $\kappa,\mu\ll\epsilon\ll\theta$, then $\eta$ sufficiently small. The following constructions have strong ordinary $CC$ end pairs.
\begin{enumerate}
\item If also $|F|\le\eta n^3$ and $\sigma(H)\ge0.65n$, there are $\lfloor\sqrt\eta n\rfloor$ disjoint nine-vertex paths of balance $-3$, chosen greedily while avoiding any prescribed set of at most $\sqrt\eta n$ vertices.
\item If $x\in D$ has $d_F(x)\ge\kappa n^2$, it lies in an eight-vertex path of balance $1$, avoiding any prescribed set of at most $\kappa n/100$ other vertices. Its other vertices are ordinary. If only $x$ is reassigned to $C$, this path has balance $-2$.
\item Let $R\subseteq D$, $|R|\le2\eta n$, and put $C'=C\cup R$, $D'=D\setminus R$. Suppose each $D'$ root misses at most $4\mu n^2$ pairs of $C'$. A matching of at most $4\eta n+3$ rare $C'D'D'$ edges, whose $C'D'$ pairs each have at least $(1-\mu^{1/3})|C'|$ completions in $C'$, extends to disjoint seven-vertex paths of temporary balance $-1$. Their four added vertices lie in $C\setminus X$; the only possible promoted vertex is the prescribed $C'$ vertex of the rare edge.
\end{enumerate}
\end{lemma}
\begin{proof}
For (i), put $t_x=e(L_H(x)[C])$ and
$k_x=e(L_H(x)[D])$ for $x\in C$.
At most $1000\eta n$ roots have $k_x>n^2/1000$,
since $\sum_{x\in C}k_x=|F|\le \eta n^3$. The block norm bound
\[
 0.65n\le \rho(L_H(x)) \le \sqrt{(|C|-1)|D|}
          +\sqrt{2t_x}+\sqrt{2k_x}
\]
gives $t_x>n^2/150$ at every other root: otherwise
the right-hand side is less than $0.64n$ for sufficiently
small $\eta$. There are at least $n/2$ such roots, so
$\sum_{x\in C}t_x=3e(H[C])$ gives
$e(H[C])\ge n^3/900$.
Removing triples meeting $X$, a nonstrong pair, or
$O(\sqrt\eta n)$ forbidden vertices costs
$O((\eta/\epsilon+\eta/\theta+\sqrt\eta)n^3)$ triples.
Thus an ordinary triple $abc$ remains whose pairs are strong.
Join $(b,c)$ to a fresh strong ordinary $CC$ pair by a
$DCCD$ connector, avoiding $a$ and the forbidden vertices,
and prefix $a$. The pattern $CCCDCCDCC$ has balance $-3$.
Throughout the greedy construction the forbidden set has
size $O(\sqrt\eta n)<\theta n$.
The connector argument uses only missing-layer bounds and
ordinary strong ends, so local optimality is unnecessary here.

For (ii), remove exceptional and forbidden vertices and
nonstrong $CD$ pairs from the link of $x$ between $C$ and
$D\setminus\{x\}$. The loss is at most
$O(\eta n^2/\epsilon+\eta n^2/\theta)+\kappa n^2/100$,
leaving more than $n$ edges and hence a path
$d_1,c_1,d_2,c_2$.
Its three link edges give the consecutive triples of the
tight path $d_1,c_1,x,d_2,c_2$, whose end pairs are strong
and ordinary. Extend it to
$a,b,d_1,c_1,x,d_2,c_2,c$ with fresh ordinary $a,b,c\in C$,
as in Lemma~\ref{old:6.3}, still avoiding the forbidden set.
It has five $C$ and three $D$ vertices, hence balance $1$.
Reassigning only $x$ to $C$ decreases its balance by three.

For (iii), write a prescribed edge as $x_1yx_2$ with
$x_i\in D'$ and $y\in C'$.
For each $i$, the missing-completion counts of $x_i z$,
summed over $z\in C'$, total at most $8\mu n^2$.
Thus at most $8\sqrt\mu n$ vertices $z\in C'$ have more
than $\sqrt\mu n$ missing completions of $x_i z$.
Choose ordinary $b\in C\cap N_H(x_1y)$ and
$c\in C\cap N_H(x_2y)$ outside these bad sets, avoiding
other matching vertices, $R$ and previous choices.
The available neighbourhoods contain almost all of $C'$,
whereas the excluded sets have size
$O(\eta n/\epsilon+\sqrt\mu n)$.
Choose ordinary $a\in C\cap N_H(x_1b)$ and
$d\in C\cap N_H(x_2c)$ with $ab,cd$ strong.
Their completion sets have size at least
$|C'|-1-\sqrt\mu n$, and $b,c$ have only
$O(\epsilon n/\theta)$ nonstrong incident pairs.
Then $a,b,x_1,y,x_2,c,d$ is the required path, with five
$C'$ and two $D'$ vertices, hence balance $-1$.
Its four added vertices lie in $C\setminus X$, so only
$y$ can belong to $R$. Greedy choices give disjointness.
\end{proof}

\begin{lemma}\label{old:6.5}
Let $n=3m+r$, $c^*=2m+r-1$, and $|C|=c^*-s$. For sufficiently small $\mu>0$ and large $n$, suppose $d_F(x)\le\mu n^2$ for every $x\in D$. Put $\tau=\mu^{1/3}$. There are $3s+3-r$ disjoint rare edges whose $CD$ pairs
each have at least $(1-\tau)|C|$ completions in $C$
in either of the following cases:
\begin{enumerate}
\item $1\le s\le\mu n$ and $\rho(L_H(x))>c^*-1$ for every $x\in D$;
\item $s=0$, $r\in\{1,2\}$, and the $D$-root radii exceed $V_m$ for $r=1$, or $2m$ for $r=2$.
\end{enumerate}
In all cases each $D$ root misses at most $4\mu n^2$ pairs of $C$.
\end{lemma}
\begin{proof}
Write $c=|C|$, $d=|D|$, and
$\lambda_x=\rho(L_H(x))$ for $x\in D$.
The spectral hypotheses imply $\lambda_x>c-1+h$,
where $h=s$ in case~(1), and $h=1/3$ or $0$ in case~(2)
according as $r=1$ or $2$.
For $r=1$, we use $V_m>2m-2/3$ from
Lemma~\ref{old:2.1}.

Call a $CD$ pair weak if it has fewer than
$(1-\tau)c$ completions in $C$.
For $x\in D$, let $M_x$ count the missing $C$ pairs
in $L_H(x)$, and put
$B_x=\{y\in C:xy\text{ is weak}\}$ and $b_x=|B_x|$.
Each weak $CD$ pair accounts for at least $\tau c/2$
missing completions, and each missing $C$ pair is counted
at most twice. Thus $b_x\le4M_x/(\tau c)$.
Write the link matrix on
$(D\setminus\{x\})\mathbin{\dot\cup}C$ as
$\left(\begin{smallmatrix}
A_x&R_x\\ R_x^T&Q_x
\end{smallmatrix}\right)$.
Then $\rho(A_x)\le d-2$,
$\|R_x\|^2\le\|R_x\|_F^2=d_F(x)$, and
$\rho(Q_x)\le c-1-M_x/c$, by
$\rho(Q_x)(\rho(Q_x)+1)\le2e(Q_x)$.
Moreover, $\lambda_x-\rho(A_x)>c-d+1\ge n/4$.
Take a unit Perron vector $v=(p,q)$, partitioned accordingly.
Its block equations are
\[
 A_xp+R_xq=\lambda_xp,\qquad
 R_x^Tp+Q_xq=\lambda_xq.
\]
Thus $p=(\lambda_xI-A_x)^{-1}R_xq$, and $q\ne0$.
Substituting into the second equation and taking the
inner product with $q$, we obtain
\[
 \begin{aligned}
 \lambda_x\|q\|^2
 &=q^TQ_xq+
   (R_xq)^T(\lambda_xI-A_x)^{-1}(R_xq)\le
 \left(\rho(Q_x)+
 \frac{\|R_x\|^2}{\lambda_x-\rho(A_x)}\right)\|q\|^2.
 \end{aligned}
\]
Dividing by $\|q\|^2$ gives
$\lambda_x\le c-1-M_x/c+4d_F(x)/n$.
Consequently $M_x\le4d_F(x)\le4\mu n^2$ and
$b_x\le\beta n$, where $\beta=32\mu^{2/3}$.

Call a rare edge good if both its $CD$ pairs are nonweak.
Suppose a maximum good matching has size $k<3s+3-r$,
and let $S_C,S_D$ be its vertex sets, of sizes $k,2k$.
Put $D_0=D\setminus S_D$ and $R=\sum_{x\in D_0}b_x$.
Every rare edge either meets $S_D$, has its colour in
$S_C$, or contains a weak pair through $D_0$;
otherwise it could be added to the matching.
These classes have at most $2k\mu n^2$,
$k\binom d2$, and $nR$ edges, respectively, so
$|F|\le nR+kn^2$.
Let $r_{xz}$ be the row sum of $R_x$ indexed by
$z\in D\setminus\{x\}$.
For $x,z\in D_0$, every colour of a rare edge through
$x,z$ belongs to $S_C\cup B_x\cup B_z$.
Hence $r_{xz}\le k+b_x+b_z\le k+2\beta n$.

Since $\lambda_x>c-1\ge n/2$, Cauchy--Schwarz in the
coordinate equations of $v$ gives
$\lambda_xv_i\le\sqrt n$, so every coordinate has
square at most $4/n$.
The formula for $p$ gives
$\|p\|^2\le16\mu\|q\|^2$, and hence $\|q\|^2\ge1/2$.
Using the second block equation again, we obtain
\[
 (h+M_x/c)\|q\|^2
 <q^T(\lambda_xI-Q_x)q
 =(R_xq)^T(\lambda_xI-A_x)^{-1}(R_xq)
 \le\frac4n\|R_xq\|^2
 \le\frac{16}{n^2}\sum_{z\ne x}r_{xz}^2.
\]
The last inequality uses the coordinate bound on $q$.
Consequently
\begin{equation}\label{eq:old:9}
 ch+M_x<\frac{32}{n}\sum_{z\ne x}r_{xz}^2.
\end{equation}
The inequality remains strict when $h=M_x=0$,
since $q\ne0$ and $\lambda_x>\rho(Q_x)$.

Sum \eqref{eq:old:9} over $D_0$.
For rows indexed by $D_0$, use
$r_{xz}\le k+2\beta n$ and
$\sum_{x,z\in D_0,\ x\ne z}r_{xz}\le2|F|$.
For rows indexed by $S_D$, use $r_{xz}^2\le cr_{xz}$ and
$\sum_{x\in D_0,z\in S_D}r_{xz}
 \le\sum_{z\in S_D}d_F(z)\le2k\mu n^2$.
Since $c\ge n/2$, $|D_0|\ge n/5$, and
$\sum_{x\in D_0}M_x\ge\tau nR/8$, we obtain
\begin{equation}\label{eq:old:10}
 \frac{hn^2}{10}+\frac{\tau nR}{8}
 <64(2\beta+k/n)(nR+kn^2)+64\mu kn^2.
\end{equation}

Choose $\mu$ so small that
$64(2\beta+5\mu)\le\tau/16$ and
$320(2\beta+6\mu)<1/10$.
In case~(1), $k\le5s\le5\mu n$.
Absorbing the $nR$ term in \eqref{eq:old:10} gives
$sn^2/10<64(2\beta+6\mu)kn^2
 \le320(2\beta+6\mu)sn^2<sn^2/10$,
a contradiction.
In case~(2) with $r=1$, we have $k\le1$ and
$k/n\le5\mu$ for large $n$.
The same absorption gives
$n^2/30<64(2\beta+6\mu)n^2<n^2/50$.
With $r=2$, we have $k=h=0$, and
\eqref{eq:old:10} gives
$\tau nR/8<128\beta nR$, contradicting
$128\beta<\tau/8$, also when $R=0$.
\end{proof}

\begin{lemma}\label{old:6.6}
For every fixed $K$, a graph partitioned into sets of sizes $2m-1$ and $m$, with at most $K$ cross edges, has radius less than $U_m$ for sufficiently large $m$.
\end{lemma}
\begin{proof}
Denote the parts by $V_1,V_2$, with $|V_1|=2m-1$
and $|V_2|=m$.
Completing both induced graphs can only increase the
spectral radius. Write the resulting adjacency matrix
as $A+E$, where $A$ is the adjacency matrix of
$K_{2m-1}\cup K_m$ and
$E=\left(\begin{smallmatrix}0&B\\B^T&0\end{smallmatrix}\right)$,
with at most $K$ ones in $B$.
The largest eigenvalue of $A$ is $2m-2$, its spectral
gap is $m-1$, and a unit Perron vector $u$ equals
$(2m-1)^{-1/2}$ on $V_1$ and zero on $V_2$.
Thus $u^TEu=0$ and
$\|E\|=\|B\|\le\|B\|_F\le\sqrt K$.
If $t_y$ is the cross degree of $y\in V_2$, then
\[
 \|Eu\|^2
 =\frac{\sum_{y\in V_2}t_y^2}{2m-1}
 \le\frac{K^2}{2m-1},
\]
since $\sum_{y\in V_2}t_y\le K$.
For sufficiently large $m$, we have $m-1>2\sqrt K$,
so Lemma~\ref{old:2.2} gives
\[
 \rho(A+E)
 \le2m-2+\frac{K^2}{(2m-1)(m-1-2\sqrt K)}
 =2m-2+O_K(m^{-2})<U_m,
\]
using $U_m>2m-2+1/(3m)$ from Lemma~\ref{old:2.1}.
\end{proof}

\begin{lemma}\label{old:6.7}
For sufficiently small $\mu>0$ and sufficiently large $m$,
let $n=3m$, $|C|=2m-1$, $|D|=m+1$, and suppose
$d_F(x)\le\mu n^2$ and $\rho(L_H(x))>U_m$ for every $x\in D$.
Put $\tau=\mu^{1/3}$.
Call a $CD$ pair good if it has at least
$(1-\tau)|C|$ completions in $C$, and weak otherwise.
Call a rare edge good if both its $CD$ pairs are good.
At least one of the following occurs:
\begin{enumerate}
\item three disjoint good rare edges;
\item a tight path $c_1,x,y,c_2,z$, with distinct
$c_1,c_2\in C$ and distinct $x,y,z\in D$,
whose end pairs $c_1x,c_2z$ are good;
\item good rare edges $c_1xz,c_2zy,c_3uv$, with three
distinct $C$ vertices and five distinct $D$ vertices,
such that $c_1zc_2\in E(H)$.
\end{enumerate}
If neither (1) nor (2) occurs, at most $1536n$ normal
triples are missing, and there are at most
$12288/\tau$ weak $CD$ pairs.
\end{lemma}

\begin{proof}
Suppose neither (1) nor (2) occurs, and let $F_0$
be the good rare edges.
For $z\in D$, define a bipartite graph $Q_z$ on
$C\mathbin{\dot\cup}(D\setminus\{z\})$ by
$ax\in E(Q_z)$ if $azx\in F$ and $ax$ is good;
there is no condition on $az$.
A path $a,x,b,y$ in $Q_z$ gives the tight path
$a,x,z,b,y$ in (2), so every $Q_z$ is $P_4$-free.

Put $c=|C|$ and retain $M_x,b_x,r_{xz}$
from Lemma~\ref{old:6.5}.
Since $U_m>c-1$, the estimates in its proof give
$M_x\le4\mu n^2$, $b_x\le\beta n$, and
$M_x\le(32/n)\sum_{z\ne x}r_{xz}^2$,
where $\beta=32\mu^{2/3}$.
Put $R=\sum_{x\in D}b_x$ and $M=\sum_{x\in D}M_x$,
so $M\ge\tau nR/8$.
Let $g_{xz}=|\{a\in C:axz\in F_0\}|$.
The graph of pairs with $g_{xz}\ge3$ has matching
number at most two, since three disjoint pairs admit
three distinct colours and give (1).
The endpoints of a maximum matching form a vertex
cover of size at most four.
At a fixed root, good colour lists of size at least two
are pairwise disjoint, since otherwise $Q_x$ contains
a $P_4$. Thus the sum of squares of the large row sizes
at one root is at most $c^2$.
Counting through the cover, and bounding every other
row square by four, gives
$\sum_{x\ne z}g_{xz}^2\le8c^2+4|D|^2\le12n^2$.
Every edge in $F\setminus F_0$ contains a weak pair,
so $a_{xz}:=r_{xz}-g_{xz}$ satisfies
$a_{xz}\le2\beta n$ and
$\sum_{x\ne z}a_{xz}=2|F\setminus F_0|\le2nR$.
Using $(g_{xz}+a_{xz})^2\le2g_{xz}^2+2a_{xz}^2$,
we obtain
\[
 M\le\frac{32}{n}\sum_{x\ne z}r_{xz}^2
 \le768n+256\beta nR.
\]
Choose $256\beta\le\tau/16$.
Since $M\ge\tau nR/8$, this gives $M\le768n+M/2$,
and hence $M\le1536n$ and $R\le R_0:=12288/\tau$.

Suppose also that (3) is absent.
Write $S_{xz}=\{a\in C:axz\in F\}$ and
$S^0_{xz}=\{a\in C:axz\in F_0\}$.
Let $T=\{x\in D:b_x>0\}$, so $|T|\le R$.
For a fixed weak pair $az$, there is at most one
$x\in D\setminus T$ with $a\in S_{xz}$ and
$|S_{xz}|\ge2$: a second such $x'$ and
$b\in S_{xz}\setminus\{a\}$ would give
$b,x,a,x'$ in $Q_z$.
Let $E$ consist of $T$ and these exceptional vertices,
over all weak pairs.
Then $|E|\le2R$, and every full colour list of size
at least two incident with $x\notin E$ consists
entirely of good colours.

Define $J$ on $D$ by $xz\in E(J)$ when $|S^0_{xz}|\ge2$.
Every $x\notin E$ has a neighbour in $J$.
Otherwise each cross row of its link has at most one edge.
Complete both induced sides to cliques, take a Perron
vector of the resulting graph, and move all cross edges
to a vertex of maximum coordinate in $C$, preserving
their $D$ endpoints.
This creates no multiple edges and cannot decrease
the Rayleigh quotient.
Completing the remaining cross edges gives
$K_1\vee(K_{2m-2}\cup K_m)$, contrary to
$\rho(L_H(x))>U_m$.
Edges of $J$ sharing an endpoint have disjoint good
colour lists, since otherwise a $P_4$ occurs.
Three disjoint edges of $J$ must have the same
two-element good colour list: since each list has
size at least two, Hall's theorem otherwise gives
three distinct representatives and hence (1).

If $J$ has three disjoint edges, let $Z$ be their six
endpoints and $\{h_1,h_2\}$ their common list.
Every good rare edge with a $D$ endpoint outside $Z$
has colour in this list, since its pair is disjoint
from at least two fixed pairs.
Choose $x\in D\setminus(E\cup Z)$.
Its $J$ neighbour has exactly this list, and any
additional good rare edge through $x$ would create
a $P_4$ in $Q_x$.
Hence $x$ belongs to exactly two good rare edges.

Otherwise $J$ has a vertex cover $S$ of size at most four.
Since every vertex outside $E$ has a neighbour,
some $z\in S$ has at least $(|D|-|E|-4)/4$ neighbours
in $J$, with pairwise disjoint good colour lists.
If a good edge $c_3uv$ avoids $z$ in its $D$ pair,
discard $u,v$ from this neighbour set and choose
$c_x\in S^0_{zx}\setminus\{c_3\}$ for each remaining
neighbour $x$.
These colours are distinct.
Their quadratically many pairs give distinct normal
triples $c_xzc_y$, at most $1536n$ of which are missing.
A surviving pair gives (3), a contradiction.
Thus every good rare edge contains $z$.
The at least $m-2R$ vertices outside $E\cup\{z\}$
have pairwise disjoint good colour lists on their pairs
with $z$, of total size at most $2m-1$.
Hence some such vertex belongs to at most three
good rare edges.

In both cases, some $x\notin E$ belongs to at most
three good rare edges.
Since $x\notin T$, every edge in $F\setminus F_0$
through $x$ is determined by a weak pair at its other
$D$ endpoint, so there are at most $R$ such edges.
Its link therefore has at most $R_0+3$ cross edges,
contradicting Lemma~\ref{old:6.6}.
\end{proof}

\begin{proof}[Proof of Theorem~\ref{old:6.1}]
Fix completion parameters $\eta_{\rm tr}\ll\epsilon\ll\theta$.
Choose $\kappa\ll\eta_{\rm tr}$ with $\kappa\ll\epsilon^3$,
then $\nu\ll\kappa^3$, and finally $\eta\ll\nu$.
Here $2\kappa$ and $4\nu$ are small enough for
Lemmas~\ref{old:6.5} and~\ref{old:6.7}, respectively.
Choose $\delta$ so that Lemma~\ref{old:6.2} gives a locally
optimal partition with size and layer errors at most
$\eta n$ and $\eta n^3$.
Ordinary and strong initially refer to this partition
and thresholds $\epsilon,\theta$.
All paths and forbidden sets below have $O(\eta n)+O(1)$
vertices; take $\eta$ small and $n$ large enough for the
greedy bounds in Lemma~\ref{old:6.4} and the length allowance
$\eta_{\rm tr}n/\epsilon$ in Lemma~\ref{old:6.3}.

Put $c^*=2m+r-1$ and write $|C|=c^*-s$, so
$\Delta:=2|D|-|C|=3s+3-r$ and $|s|\le\eta n+1$.
It suffices to construct a path of balance $\Delta-2$
with strong ordinary $CC$ ends.
If $\Delta<0$, join $-\Delta$ paths from
Lemma~\ref{old:6.4}(1) by normal connectors.
Their balance is $-3(-\Delta)+2(-\Delta-1)=\Delta-2$.
If $\Delta=0$, the connector construction gives a normal
path with strong ordinary $CC$ ends and balance $-2$.
These cases include every $s\le-1$.
Assume henceforth $s\ge0$ and that $H$ is not Hamiltonian.

Call $x\in D$ rich if $d_F(x)\ge\kappa n^2$.
If there are at least $s+1$ rich vertices, use their paths
from Lemma~\ref{old:6.4}(2), together with $r$ paths from
part~(1). Joining them gives balance
$(s+1)-3r+2(s+r)=\Delta-2$, a contradiction.
Thus the rich set $R$ has size $t\le s$.
Put $C'=C\cup R$, $D'=D\setminus R$, and $s'=s-t$.
Each remaining root has at most
$\kappa n^2+tn\le2\kappa n^2$ rare edges in the new
partition; the new edges were $DDD$ triples meeting $R$.
The temporary imbalance is $\Delta'=3s'+3-r=\Delta-3t$.

A short path of temporary balance $\Delta'-2$ with original
strong ordinary $CC$ ends suffices.
Indeed, cover each unused vertex of $R$ by its rich path,
avoiding the given path and all other vertices of $R$.
Each added path has temporary balance $-2$.
Attach them by original normal connectors avoiding $R$,
each of balance two, so the temporary balance is unchanged.
Every vertex of $R$ is now covered.
Restoring the original roles increases the balance by
$3t$, giving $\Delta-2$, and Lemma~\ref{old:6.3} applies
in the original partition.

If $s'\ge1$, or $s'=0$ and $r\in\{1,2\}$, apply
Lemma~\ref{old:6.5} with parameter $2\kappa$.
It gives $\Delta'$ disjoint rare edges whose $C'D'$ pairs
each have at least $(1-(2\kappa)^{1/3})|C'|$ completions
in $C'$, and each $D'$ root misses at most $8\kappa n^2$
pairs of $C'$.
Lemma~\ref{old:6.4}(3) extends these edges to disjoint
seven-vertex paths of temporary balance $-1$ with
original strong ordinary $CC$ ends.
Join them by original normal connectors avoiding $R$.
The balance is $-\Delta'+2(\Delta'-1)=\Delta'-2$,
so the preceding argument gives a contradiction.
Only $r=0$ and $t=s$ remain.

Every vertex of $D\setminus R$ has rare degree below
$\nu n^2$: otherwise it and the $s$ vertices of $R$
admit paths from Lemma~\ref{old:6.4}(2) at threshold $\nu$,
whose concatenation has balance $3s+1=\Delta-2$.
Thus, in the temporary partition,
$|C'|=2m-1$, $|D'|=m+1$, and every $D'$-root rare degree
is at most $\nu n^2+sn\le\mu n^2$, where $\mu=4\nu$.
The block estimate in the proof of Lemma~\ref{old:6.5}
bounds each $D'$-root normal defect by $4\mu n^2$.
Hence every $D'$ vertex is original ordinary, since
$4\mu<\epsilon$.
Call a $C'D'$ pair good if it has at least
$(1-\tau)|C'|$ completions in $C'$, where
$\tau=\mu^{1/3}\ll\kappa$.

Put $C_0=C\setminus X$.
Every $z\in C'$ has at least $\kappa n^2/3$ link edges
between $C_0\setminus\{z\}$ and $D'$.
For original $C$ vertices, use the role bound $n^2/50$;
for vertices of $R$, use their original rare degree
at least $\kappa n^2$.
Deleting $X$ and the vertices of $R$ from these links
costs $O(\eta n^2/\epsilon)$ edges.

We exclude (1) and (2) of Lemma~\ref{old:6.7}.
In (1), Lemma~\ref{old:6.4}(3) gives three paths of
temporary balance $-1$ with original strong ordinary
$CC$ ends. Two original normal connectors avoiding $R$
give balance one, which equals $\Delta'-2$.
Covering the unused vertices of $R$ as above gives a
contradiction.

For (2), write the core as $c_1,v,u,c_2,w$.
In the role link of $c_1$, restrict the $C_0$ endpoint
to $N_H(c_1v)\cap C$, delete nonstrong original $CD$
pairs, and avoid the core vertices.
The loss is at most $\tau n^2+O(\eta n^2/\theta+n)$,
less than $\kappa n^2/3$.
Choose a surviving edge $ad$, with $a\in C_0$, $d\in D'$.
Thus $c_1ad,ac_1v\in E(H)$ and $ad$ is strong.
Choose fresh ordinary $p\in C\cap N_H(ad)$ with $pd$
strong, then $q\in C\cap N_H(pd)$ with $qp$ strong.
At the other end choose fresh ordinary
$b\in C\cap N_H(c_2w)$ with $wb$ strong, then
$f\in C\cap N_H(wb)$ with $bf$ strong.
These choices are possible since $d,w$ are original
ordinary, the outer pairs are good, and
$\epsilon\ll\theta$.
The path
\begin{equation}\label{eq:old:11}
 q,p,d,a,c_1,v,u,c_2,w,b,f
\end{equation}
has seven $C'$ and four $D'$ vertices, hence temporary
balance one, and original strong ordinary $CC$ ends.
Again, covering the unused vertices of $R$ gives a
contradiction.

Lemma~\ref{old:6.7} now gives (3) and at most $1536n$
missing $C'C'D'$ triples.
From now on, ordinary and strong refer to $(C',D')$.
Every vertex is ordinary and belongs to only $O(1)$
nonstrong pairs.
This partition is locally optimal: by the move formulas
in Lemma~\ref{old:6.2}, the costs for $z\in C'$ and
$x\in D'$ are at least
$(|C'|-1)|D'|-\binom{|D'|}{2}-3072n>0$ and
$|C'|(|D'|-1)-\binom{|D'|-1}{2}-2\mu n^2>0$,
respectively.

Write (3) as $c_1xz,c_2zy,c_3uv$, with
$c_1zc_2\in E(H)$.
Choose fresh $b\in C'\cap N_H(xc_1)$ with $xb$ strong,
then $a\in C'\cap N_H(xb)$ with $ab$ strong.
Similarly choose fresh $c\in C'\cap N_H(yc_2)$ with
$yc$ strong, then $d\in C'\cap N_H(yc)$ with $cd$ strong.
The good outer pairs have linear neighbourhoods, and
each vertex has only $O(1)$ nonstrong partners.
Thus these choices can avoid all prescribed vertices.
The tight path $a,b,x,c_1,z,c_2,y,c,d$ has balance zero;
its middle triple is $c_1zc_2$.
As in Lemma~\ref{old:6.4}(3), extend $c_3uv$ to a
disjoint seven-vertex path of balance $-1$.
A normal connector gives a path of balance one,
equal to $2|D'|-|C'|-2$.
Lemma~\ref{old:6.3} in the new locally optimal partition
gives the final contradiction.
\end{proof}

\begin{remark}\label{rem:twohit-sharpness}
For $r=0$, the bound $U_m$ in Theorem~\ref{old:6.1} is sharp for all sufficiently large $m$. Take $V(H)=\{h\}\mathbin{\dot\cup}C'\mathbin{\dot\cup}D$, with $|C'|=2m-2$ and $|D|=m+1$, and include all triples through $h$ and all triples of types $C'C'C',C'C'D,DDD$. The tight-component argument in Proposition~\ref{prop:sharp} excludes a Hamilton cycle: deleting $h$ would leave a spanning path with its $D$ vertices at distance at least three, allowing at most $m$ of them. The $D$-links are $K_1\vee(K_{2m-2}\cup K_m)$, while the $C'$-links are $K_{2m-2}\vee I_{m+1}$, whose radius exceeds $2m$ for $m\ge4$ because its quadratic at $2m$ equals $-2m^2+6m+2<0$. The $h$-link is complete, so $\sigma(H)=U_m$. With $C=C'\cup\{h\}$, we have $M=0$ and $|F|=\binom{m+1}{2}$, proving local sharpness; this construction attains the global threshold $T_{3m}$ precisely when $3\nmid m$.
\end{remark}

\section{Tight Hamilton cycles near the missing-layer barrier}\label{sec:smallpart}
\subsection{Crossings and vertex roles}

In this section we consider partitions
$V(H)=A\mathbin{\dot\cup}B$ with $|A|$ close to $n/3$.
The missing-layer barrier contains all triples of types
$AAA,AAB,BBB$ and none of type $ABB$.
The $BBB$ triples form a separate tight component,
so every edge of a tight Hamilton cycle would contain
at least two $A$ vertices, forcing $3|A|\ge2n$,
which fails for these part sizes.
We call triples of types $AAA,AAB,BBB$ normal
and those of type $ABB$ rare. Lemma~\ref{old:4.11}(3) describes the corresponding
asymptotic regime, with $o(n^3)$ rare triples and
$o(n^3)$ missing normal triples.
A \emph{crossing} is a vertex-simple tight path whose first
pair lies in the $A$ core and whose last pair lies in the
$B$ core. Theorem~\ref{old:7.2} extends two disjoint crossings to a Hamilton cycle,
incorporating the exceptional vertices. Lemma~\ref{old:7.4} then reduces
the non-Hamiltonian case to a partition with at most one vertex outside
the pointwise role conditions. The subsequent spectral arguments handle
the case with no exceptional vertex; if one remains,
Corollary~\ref{old:8.9} produces a crossing avoiding it, contradicting
the conclusion of Lemma~\ref{old:7.4}. Endpoint allowances are always
fixed before the smaller defect parameters.

\begin{lemma}\label{old:7.1}
Let $1/N\ll\varepsilon\ll\delta\ll\tau\ll1$, and let $V(F)=V_0\mathbin{\dot\cup}Z$, where $|V_0|=N$. Suppose every vertex of $V_0$ belongs to at most $\varepsilon N^2$ missing triples of $F[V_0]$. Let $\mathcal P$ be vertex-disjoint tight paths covering $Z$ and using at most $\delta N$ vertices, with both end pairs in $V_0$ and having at least $(1-\tau)N$ completions there. Given disjoint ordered pairs $(s_1,s_2),(t_1,t_2)$ outside these paths with the same completion bound, $F$ has a spanning tight path with these ordered ends, containing every member of $\mathcal P$ as a segment in either prescribed orientation.
\end{lemma}

\begin{proof}
Put $\theta=\sqrt\varepsilon<\tau$ and call a core pair $uv$
good if it has at most $\theta N$ missing completions in
$V_0\setminus\{u,v\}$.
Each core vertex has at most $2\theta N$ non-good partners.
Regard the oriented members of $\mathcal P$, preceded by
$s_1s_2$ and followed by $t_1t_2$, as $k=|\mathcal P|+2$
fixed blocks. Let $r$ be the number of unused core vertices.
Choose $j\in\{0,1,2\}$ with $r-j\equiv k-1\pmod3$ and
append $j$ fresh core vertices to the first block, each time
making its new last pair good. The completion bound and
the bound on non-good partners permit these choices while
avoiding all other blocks.

Put $g=(r-j-k+1)/3$ and $a=g+k-1$, so $a+2g=r-j$.
For small $\delta$ and large $N$, we have $g\ge k-1$.
Choose $a$ unused core vertices as $U$ and reserve the
remaining $2g$ for $g$ ordinary two-vertex blocks.
Insert the ordinary blocks between the fixed blocks, with
at least one in each gap, and put a $U$-position between
consecutive blocks.
There are exactly $a$ such positions. Let $B$ consist of
the vertices reserved for ordinary blocks and the vertices
in the boundary pairs of the fixed blocks. Both $U$ and $B$
have size at least $|V(F)|/5$ for sufficiently small $\delta$.

Let $\Gamma$ contain all good pairs on $B$ and all fixed
boundary pairs. Then $\Delta(\overline\Gamma)\le2\theta N$,
and each edge misses at most $\tau N$ completions in $U$.
Each $u\in U$ misses at most $\varepsilon N^2$ pairs of $B$.
At most $\delta N+6$ entries are prescribed, and all required
pairs with two prescribed endpoints are fixed boundary pairs.
Lemma~\ref{lem:common-completion} applies to $F$ with
$\beta=1/5$ and defect parameter $2\tau$.
The resulting spanning path has the prescribed ordered ends
and retains every member of $\mathcal P$ in its chosen orientation.
\end{proof}

\begin{theorem}\label{old:7.2}
For sufficiently small $\varepsilon,\delta,\tau>0$ and large $n$, let $V(H)=A_0\mathbin{\dot\cup}B_0\mathbin{\dot\cup}X$, where $0.29n\le|A_0|\le0.37n$ and $|X|\le\delta n$. Suppose $\sigma(H)\ge0.65n$, and every vertex of each core $P\in\{A_0,B_0\}$ belongs to at most $\varepsilon n^2$ missing triples inside $P$. If there are two vertex-disjoint crossings $J_1,J_2$, together using at most $\delta n$ vertices, whose end pairs miss at most $\tau n$ completions in their respective cores, then $H$ has a tight Hamilton cycle containing both crossings, one reversed. Their interiors may contain vertices of $X$.
\end{theorem}

\begin{proof}
Choose $\varepsilon\ll\delta\ll\tau\ll1$ and take $n$
sufficiently large. Reserve $R=V(J_1)\cup V(J_2)$.
For $x\in X\setminus R$, split $L_H(x)$ across $A_0$
and $(B_0\cup X)\setminus\{x\}$.
The off-diagonal part has norm at most
$\sqrt{|A_0|(n-|A_0|)}<0.49n$, so one diagonal block
has radius greater than $0.16n$ and more than
$0.012n^2$ edges. Deleting $X$ removes at most
$\delta n^2$ edges, leaving a core
$P_x\in\{A_0,B_0\}$ with more than $0.011n^2$ link edges.

Call a core pair strong if it misses at most $\tau n$
completions in that core. Double counting missing triples
gives at most $2\varepsilon n^2/\tau$ non-strong pairs
in each core. Process $X\setminus R$ greedily.
Deleting reserved and previously used core vertices,
and all non-strong pairs, from $L_H(x)[P_x]$ removes
at most $5\delta n^2+2\varepsilon n^2/\tau$ edges.
The remaining graph has more than $n$ edges, so it contains
a four-vertex path $p_1p_2p_3p_4$: a graph with no such
path is a disjoint union of stars, triangles and isolated
vertices and has at most as many edges as vertices.
Hence $p_1,p_2,x,p_3,p_4$ is a tight path with strong
ends in $P_x$. These paths are disjoint and avoid $R$;
assign each to its chosen core.

Write the ordered ends of $J_1,J_2$ as
$u_1u_2,v_1v_2$ and $u'_1u'_2,v'_1v'_2$, respectively.
Delete from each core all vertices of $R$ except the
crossing end vertices, obtaining $A',B'$.
Each reduced core $P'$ has order $N\ge0.28n$ and
pointwise defect at most
$\varepsilon n^2\le13\varepsilon N^2$.
Its assigned paths use at most
$5\delta n\le20\delta N$ vertices, and all prescribed
end pairs have at least
$N-(\tau+\delta)n-2\ge(1-5\tau)N$ completions there.
Apply Lemma~\ref{old:7.1} to each reduced core together
with its assigned exceptional vertices, using parameters
$13\varepsilon,20\delta,5\tau$.
Prescribe ordered ends $(u'_2,u'_1),(u_1,u_2)$ in $A'$
and $(v_1,v_2),(v'_2,v'_1)$ in $B'$, obtaining spanning
paths $P_A,P_B$ that preserve the assigned paths.

Concatenate $P_A,J_1,P_B$ and the reverse of $J_2$,
identifying common end pairs. Each vertex occurs once,
and every cyclic triple consists of consecutive vertices
in one of the four paths. This gives the required
Hamilton cycle.
\end{proof}

\begin{lemma}\label{old:7.3}
For every $\varepsilon>0$ there is $\delta>0$ such that the following holds for large $n$. Let $V(H)=A_0\mathbin{\dot\cup}B_0\mathbin{\dot\cup}X$, with $|X|\le\delta n$, $|A_0|=n/3\pm\delta n$, and $|B_0|=2n/3\pm\delta n$. Let $\Gamma_A,\Gamma_B$ be graphs on the two cores and $\Omega\subseteq A_0\times B_0$, where
\[
\Delta(\overline{\Gamma_A})\le |A_0|/4,\qquad e(\overline{\Gamma_B})\le\delta n^2,\qquad |(A_0\times B_0)\setminus\Omega|\le\delta n^2.
\]
Every pair of $\Omega$ has at least $|A_0|/2$ completions in $A_0$. If $x\in X$ satisfies $\rho(L_H(x))\ge(2/3-\delta)n$, then for every $S\subseteq V(H)\setminus\{x\}$ of size at most $\delta n$, at least one of the following holds:
\begin{enumerate}
\item there is a crossing $a_0,a,b_1,x,b_2,b_3$ avoiding $S$, with $a_0a\in\Gamma_A$ and $b_2b_3\in\Gamma_B$;
\item $L_H(x)$ misses at most $\varepsilon n^2$ pairs of $\binom{A_0}{2}\cup(A_0\times B_0)$;
\item $L_H(x)$ misses at most $\varepsilon n^2$ pairs of $\binom{B_0}{2}$.
\end{enumerate}
\end{lemma}

\begin{proof}
Delete $X\setminus\{x\}$ and $S$ from the link, together
with cross edges outside $\Omega$ and $B_0$ edges outside
$\Gamma_B$. Let $G$ be the resulting graph on
$A=A_0\setminus S$, $B=B_0\setminus S$.
At most $4\delta n^2$ edges are removed, so the Frobenius
bound gives $\rho(G)\ge(2/3-\delta-\sqrt{8\delta})n$.
Let $T\subseteq B$ be the vertices with a neighbour in $A$,
and put $J=G[B]$.

If $u,v,w$ is a three-vertex path in $J$ with $u\in T$,
choose $a\in N_G(u)\cap A$. Since $au\in\Omega$, we can
choose $a_0\in A_0\setminus S$ with $a_0a\in\Gamma_A$
and $a_0au\in H$: at least
$|A_0|/2-|A_0|/4-|S|-1>0$ choices remain.
Then $a_0,a,u,x,v,w$ is the crossing in (1).

Otherwise every neighbour in $J$ of a vertex of $T$
has degree one, so every component of $J$ meeting $T$
is a star or an isolated vertex.
Delete all edges of these components, at most $|B|$
in total. With $a=|A|$, $b=|B|$ and $t=|T|$, the
remaining graph $G'$ is contained in
$(K_a\vee I_t)\mathbin{\dot\cup}K_{b-t}$ and satisfies
$\rho(G')\ge(2/3-5\sqrt\delta)n$.
Put $s(a,t)=\rho(K_a\vee I_t)$, so
$s(a,b)=2n/3+O(\delta n)$. For $0\le t\le b$,
\[
 \frac16\le\frac{\partial s}{\partial t}
 =\frac{a}{\sqrt{(a-1)^2+4at}}\le2.
\]
If the clique template reaches the required radius,
then $t=O(\sqrt\delta n)$; otherwise the split template
does, and the lower derivative bound gives
$b-t=O(\sqrt\delta n)$.

In the first case,
$s(a,t)\le a-1+2t=n/3+O(\sqrt\delta n)<n/2$,
so $G'[B\setminus T]$ has radius at least
$(2/3-5\sqrt\delta)n$.
In the second case, the clique block has order
$O(\sqrt\delta n)$, so $G'[A\cup T]$ has this radius.
Restore all other vertices of $V(H)\setminus\{x\}$
as isolated vertices. These graphs are subgraphs of
the clique template on $B_0$ and the split template
with clique $A_0$, respectively.
Apply the uniform form of Lemma~\ref{old:4.3}, choosing
$\delta$ sufficiently small in terms of $\varepsilon$,
to obtain at most $\varepsilon n^2$ missing template
edges. Since the retained graphs are subgraphs of
$L_H(x)$, conclusions (3) and (2), respectively, follow.
\end{proof}

\begin{lemma}\label{old:7.4}
Fix a sufficiently small endpoint allowance $\theta>0$. For every sufficiently small $\xi>0$ there is $\zeta=\zeta(\theta,\xi)>0$ such that the following holds for large $n$. Suppose $\sigma(H)\ge(2/3-\zeta)n$, and $H$ admits a partition $A\mathbin{\dot\cup}B$ with $\bigl||A|-n/3\bigr|\le\zeta n$ and at most $\zeta n^3$ missing normal triples. If $H$ is not Hamiltonian, it has a partition $A'\mathbin{\dot\cup}B'\mathbin{\dot\cup}Z$, where $|Z|\le1$, $|A'|=n/3\pm\xi n$, and $|B'|=2n/3\pm\xi n$. Every vertex of $A'$ belongs to at most $\xi n^2$ missing triples of each type $A'A'A',A'A'B'$, and every vertex of $B'$ to at most $\xi n^2$ missing triples inside $B'$. There are cores $A_0\subseteq A'$, $B_0\subseteq B'$ omitting at most $\xi n$ vertices in total. If $Z=\{h\}$, then $H-h$ has no crossing on at most $20$ vertices whose end pairs miss at most $\theta n$ completions in their respective cores.
\end{lemma}

\begin{proof}
Fix allowances $\varepsilon_c,\delta_c$ for
Theorem~\ref{old:7.2} with endpoint allowance $\theta$.
Choose $\alpha\le\xi/4$ and an allowance
$\delta_r\le1/100$ from Lemma~\ref{old:7.3} for $\alpha$.
Choose
$\kappa\le\min\{\xi/4,\varepsilon_c,\theta/100,\delta_r\}$,
then $\zeta\le\min\{\xi/4,\delta_r/4\}$ so small that
\[
 3\zeta/\kappa\le\min\{\xi,\delta_r,\delta_c\}/4,\qquad
 6\zeta/\theta\le\delta_r,\qquad 20\zeta\le\delta_r.
\]
Take $n$ sufficiently large, with
$20\le\delta_r n$ and $26\le\delta_c n$.

Let $X$ consist of vertices incident with at least
$\kappa n^2$ missing normal triples.
Then $|X|\le3\zeta n/\kappa$.
Put $A_0=A\setminus X$ and $B_0=B\setminus X$.
Let $\Gamma_A,\Gamma_B$ consist of core pairs missing
at most $\theta n$ completions in their cores, and let
$\Omega$ consist of cross-pairs with at least
$|A_0|/2$ completions in $A_0$.
Counting missing triples gives
\[
 \Delta(\overline{\Gamma_A})\le4\kappa n/\theta<|A_0|/4,\quad
 e(\overline{\Gamma_B})\le6\zeta n^2/\theta,\quad
 |(A_0\times B_0)\setminus\Omega|\le20\zeta n^2.
\]
The last two bounds are at most $\delta_r n^2$.
Together with $|X|\le\delta_r n$, part-size errors at most
$(\zeta+3\zeta/\kappa)n\le\delta_r n$, and $\zeta\le\delta_r$,
these allow Lemma~\ref{old:7.3} at every $x\in X$.

Let $Z\subseteq X$ consist of vertices satisfying neither
conclusion (2) nor (3) of that lemma.
Each gives a six-vertex crossing with ends in
$\Gamma_A,\Gamma_B$, avoiding any prescribed set of at most
$\delta_r n$ other vertices.
Two give disjoint crossings by first avoiding the other
root, then the first crossing.
Since $\kappa\le\varepsilon_c$, $|X|\le\delta_c n$,
$0.29n\le|A_0|\le0.37n$ and $\sigma(H)\ge0.65n$,
Theorem~\ref{old:7.2} would make $H$ Hamiltonian.
Hence $|Z|\le1$.

Assign $X\setminus Z$ to $A'$ or $B'$ according to
conclusion (2) or (3), keeping the cores fixed.
Among pairs wholly in the cores, a core root misses at most
$\kappa n^2$ required pairs, and a reassigned root at most
$\alpha n^2$.
The remaining relevant pairs meet $X$ and number at most
$|X|n$. Thus every required pointwise defect is at most
$(\alpha+\kappa+3\zeta/\kappa)n^2\le\xi n^2$.
The part-size errors are at most
$(\zeta+3\zeta/\kappa)n\le\xi n$,
and the cores omit at most $|X|\le\xi n$ vertices.

If $Z=\{h\}$ and a crossing $Q\subseteq H-h$ on at most
$20$ vertices had the stated end pairs,
Lemma~\ref{old:7.3} would give a crossing through $h$
avoiding $V(Q)$.
Together they use at most $26$ vertices and have endpoint
allowance $\theta n$, so Theorem~\ref{old:7.2} again
gives a contradiction.
\end{proof}

\subsection{Pair families and spectral estimates}

\begin{lemma}\label{old:8.1}
Let $(\mathcal F_i)$ be pairwise cross-intersecting families of pairs on a set $Q$. Then one of the following holds.
\begin{enumerate}
\item Some $\mathcal F_{i_0}$ contains two disjoint pairs on a four-set $W$, and every other family is supported on $W$. At each vertex $v\notin W$, all incident pairs consequently have colour $i_0$. If another colour is nonempty, at most two pairs of $\mathcal F_{i_0}$ contain $v$; if another colour has at least two pairs, at most one pair of $\mathcal F_{i_0}$ contains $v$.
\item The union of the families is intersecting, and is contained in a star or a triangle.
\end{enumerate}
\end{lemma}

\begin{proof}
If $\mathcal F_{i_0}$ contains disjoint pairs with union $W$,
every pair of another colour meets both, hence lies in $W$.
For $v\notin W$, any $vu\in\mathcal F_{i_0}$ has
$u\in f$ for every pair $f$ of another colour.
One such pair leaves at most two choices for $u$,
and two distinct pairs leave at most one.

Otherwise cross-intersection makes the union intersecting.
If it is not contained in a star, choose $ab,ac$ and a pair
avoiding $a$. This pair is $bc$, and every pair meeting
$ab,ac,bc$ lies in $\{a,b,c\}$.
\end{proof}

The next estimates retain the constant and reciprocal-linear error terms needed below.

\begin{lemma}\label{old:8.2}
Let $G$ have a clique $P$ of order $a\ge1$ and a complementary part $Q$ of order $b$. Suppose $J=G[Q]$ has $e$ edges and maximum degree $\Delta$, and $M$ cross-pairs are missing. If $\lambda=\rho(G)>\Delta$, then
\begin{equation}\label{eq:old:12}
M\le ab-\lambda(\lambda-a+1)+\frac{2ae}{\lambda-\Delta}.
\end{equation}
In particular, if $F$ has $b$ vertices, $e$ edges, and maximum degree $\Delta<\Lambda$, then
\begin{equation}\label{eq:old:13}
\rho(K_a\vee F)>\Lambda\quad\Longrightarrow\quad e>\frac{\Lambda-\Delta}{2a}\bigl(\Lambda^2-(a-1)\Lambda-ab\bigr).
\end{equation}
\end{lemma}

\begin{proof}
Since $\lambda>\Delta\ge\rho(J)$, a nonnegative Perron
vector $(x,y)$ has nonzero $P$ part; normalize
$\max_{u\in P}x_u=1$.
Put $Y=\sum_{v\in Q}y_v$ and $z=\|y\|_\infty$.
The eigenvalue equations give
$\lambda\le a-1+Y$ and $(\lambda-\Delta)z\le a$.
Summing over $Q$, each cross edge contributes at most $1$
and each edge of $J$ at most $2z$, so
\[
 \lambda(\lambda-a+1)\le\lambda Y
 \le ab-M+2ez
 \le ab-M+\frac{2ae}{\lambda-\Delta}.
\]
This proves~\eqref{eq:old:12}.

For the second assertion, put $f(t)=t^2-(a-1)t-ab$ and
$\lambda=\rho(K_a\vee F)>\Lambda$.
If $f(\Lambda)<0$, the conclusion follows from $e\ge0$.
Otherwise $\Lambda\ge a-1$, and $f$ is strictly increasing
on $[\Lambda,\infty)$.
Applying~\eqref{eq:old:12} with $M=0$ gives
\[
 2ae\ge(\lambda-\Delta)f(\lambda)
 >(\Lambda-\Delta)f(\Lambda),
\]
proving~\eqref{eq:old:13}.
\end{proof}

\begin{lemma}\label{old:8.3}
Let $J=K_q\vee F$, where $q\ge1$, $\Delta(F)\le D$, and $\rho(J)\ge q+D$. Let $W\subseteq V(F)$. Suppose $G$ is obtained from a subgraph of $J$ by adding edges only inside $W$, and every vertex of $W$ misses at least $|W|-1$ neighbours in the clique. Then $\rho(G)\le\rho(J)$.
\end{lemma}

\begin{proof}
By symmetry, normalize a positive Perron vector $x$ of $J$
to $1$ on the clique. Put $\lambda=\rho(J)$ and
$z=\max_{v\in V(F)}x_v$.
The $F$ equations give $(\lambda-D)z\le q$, hence $z\le1$.
At each retained vertex of $W$, at least $|W|-1$ unit clique
contributions are lost, while at most $|W|-1$ contributions
bounded by $1$ are added.
Other vertices only lose contributions.
Thus $\sum_{u\in N_G(v)}x_u\le\lambda x_v$ for every
$v\in V(G)$, and the weighted row-sum bound gives
$\rho(G)\le\lambda$.
\end{proof}

\begin{lemma}\label{old:8.4}
Let $G$ have a partition $P\mathbin{\dot\cup}Q$, with $|P|=p$, $|Q|=q$, and cross-adjacency matrix $R$. Put $S_2=\sum_{j\in Q}d_j^2$, where $d_j=\sum_{i\in P}R_{ij}$. If $1.9m\le p\le2.1m$, $0.9m\le q\le1.1m$, and $d_j\le m/100$ for every $j$, then
\begin{equation}\label{eq:old:14}
 \rho(G)\le p-1+\frac{2S_2}{m^2}.
\end{equation}
\end{lemma}

\begin{proof}
Complete both parts to cliques and let $\lambda$ be the
resulting radius. Since $\lambda\ge p-1>q-1$, take a
nonnegative Perron vector $(x,y)$ with $\max_{i\in P}x_i=1$.
Put $X=\sum_i x_i$, $Y=\sum_j y_j$, $S_1=\sum_jd_j$ and
$h=\lambda-q+1$.
The equations on $Q$ give
\[
 (\lambda+1)y_j=Y+\sum_iR_{ij}x_i,\qquad hY\le S_1.
\]
Using $S_1^2\le qS_2$ and $h+q=\lambda+1$, we get
$\sum_jd_jy_j\le(S_2+S_1^2/h)/(\lambda+1)\le S_2/h$.
At a vertex with $x_i=1$, we have $X\ge\lambda+1-Y$.
Hence
\[
 hX\ge(\lambda+1)h-S_1
 \ge p(p-q)-S_1
 \ge m^2/2,
\]
by the size bounds and $S_1\le qm/100$.
Finally, summing the equations on $P$ gives
$(\lambda-p+1)X=\sum_jd_jy_j\le S_2/h$.
Thus $\rho(G)\le\lambda\le p-1+2S_2/m^2$.
\end{proof}

\subsection{The nonzero residues and a prescribed exceptional vertex}

\begin{theorem}\label{old:8.5}
There are constants $\eta,\varepsilon>0$ such that the following holds for all sufficiently large $m$. Let $r\in\{1,2\}$ and let $V(H)=A\mathbin{\dot\cup}B$, where $|A|=m+s$, $|B|=2m+r-s$, and $|s|\leq\eta m$. Suppose every vertex of $A$ belongs to at most $\varepsilon m^2$ missing triples of each type $AAA,AAB$, and every vertex of $B$ belongs to at most $\varepsilon m^2$ missing triples of type $BBB$. Put $\Lambda_1=V_m$ and $\Lambda_2=2m$. If $\sigma(H)>\Lambda_r$, then $H$ has a tight Hamilton cycle.
\end{theorem}

\begin{proof}
Choose $0<\tau,\gamma\ll1$, then $\eta,\varepsilon>0$
sufficiently small, and suppose that $H$ is not Hamiltonian.
Write $a=|A|$, $b=|B|$, $\Lambda=\Lambda_r$ and
$\theta=2m-\Lambda$.
Thus $\theta=0$ for $r=2$, and $0<\theta<2/3$ for $r=1$.
All constants denoted by $C$ are absolute.
For $P=A,B$, let $\Gamma_P$ consist of pairs with at least
$(1-\tau)|P|$ completions in $H[P]$.
For $x\in A$, put
\[
 W_x=\{y\in B:|N_H(xy)\cap A|<(1-\tau)a\},
 \qquad w_x=|W_x|.
\]
The pointwise hypotheses give
$\Delta(\overline{\Gamma_P}),w_x\le C\varepsilon m/\tau=\omega m$.
Let $D$ consist of the $B$-pairs with fewer than $\gamma m$
completions in $B$, and put $E=e(D)$ and $d=\Delta(D)$.
Counting missing triples through a $B$-vertex gives
$d(b-2-\gamma m)\le2\varepsilon m^2$, hence $d\le2\varepsilon m$.
Take $\varepsilon$ small enough that $\omega\ll\tau,\gamma$.

For each $x\in A$, let $G_x$ consist of the pairs
$uv\in\binom B2$ with $xuv\in H$, $uv\notin D$ and
$\{u,v\}\not\subseteq W_x$.
Call these triples \emph{usable}.
Such a triple extends to a crossing $p,x,u,v,z$:
orient $uv$ so that $u\notin W_x$, then choose
$p\in A$ with $pxu\in H$ and $px\in\Gamma_A$, and
$z\in N_H(uv)\cap B$ with $vz\in\Gamma_B$.
There are at least $(1-\tau)a-\omega m-O(1)$ and
$\gamma m-\omega m-O(1)$ choices.
Thus disjoint pairs in distinct families give disjoint
crossings by avoiding the reserved vertices.
Theorem~\ref{old:7.2}, with cores $A,B$ and no exceptional
vertices, would give a Hamilton cycle.
Hence the families $G_x$ are pairwise cross-intersecting.

By Lemma~\ref{old:8.1}, either there are $h\in A$ and
$S\subseteq B$, $|S|\le4$, supporting every $G_x$ with
$x\ne h$, or all $G_x$ lie in a star centred at $z\in B$.
The triangle case is included in the former by choosing
$h$ arbitrarily.
Put $j=0$, $A^*=A\setminus\{h\}$ in the first case, and
$j=1$, $A^*=A$ in the second.
For every $x\in A^*$,
\begin{equation}\label{eq:old:17}
 L_H(x)[B]\subseteq G_x\cup D\cup\binom{W_x}{2}.
\end{equation}
In the first case set $J_x=K_{a-1}\vee(D\cup K_S)$.
In the second set $J_x=K_a\vee D[B\setminus\{z\}]$,
with clique $(A\setminus\{x\})\cup\{z\}$.
Its clique order is $q=a-1+j$ and the complementary
maximum degree is at most $d+3$.
Every weak column outside the clique misses at least
$\tau a-1\ge\tau a/2$ clique neighbours, while $w_x\le\omega m$.
Also $\rho(J_x)\ge1.8m>q+d+3$.
Lemma~\ref{old:8.3}, applied to $W_x$ or $W_x\setminus\{z\}$,
therefore gives
\begin{equation}\label{eq:old:18}
 \rho(L_H(x))\le\rho(J_x).
\end{equation}

We next show
\begin{equation}\label{eq:old:19}
 W:=\sum_{x\in A^*}w_x
 \le\frac C\tau\bigl(m(s_++1)+E\bigr),
 \qquad s_+=\max\{s,0\}.
\end{equation}
Complete only the above clique in $L_H(x)$, retaining its
cross-defects, and call the resulting graph $\widehat L_x$.
Its complementary part has at most $E+6+w_x^2/2$ edges
and maximum degree at most $d+3+w_x$.
Its radius exceeds $\Lambda$.
The expression $q(b-j)-t(t-q+1)$ decreases for $t\ge\Lambda$,
and its value at $\Lambda$ is at most $Cm(s_++1)$.
Moreover, $2q/(t-d-3-w_x)$ is bounded in this range.
Lemma~\ref{old:8.2} thus bounds the number $M_x$ of missing
cross-pairs by
\[
 M_x\le C[m(s_++1)+E+w_x^2].
\]
On the other hand, $M_x\ge(\tau a/2)(w_x-j)$, since only
the promoted vertex can be omitted.
Using $w_x^2\le\omega m w_x$ and $\omega\ll\tau$, absorb
this term to obtain
$w_x\le(C/\tau)(s_++1+E/m)$.
Summing proves~\eqref{eq:old:19}.

For $y\in B$, put $T=N_D(y)$ and $t=|T|$.
In the first case take $y\notin S$ and
\[
 P=(B\setminus(T\cup\{y\}))\cup\{h\},
 \qquad Q=(A\setminus\{h\})\cup T.
\]
In the star case take $y\ne z$ and
$P=B\setminus(T\cup\{y\})$, $Q=A\cup T$.
Then $|P|=b-j-t$ and $|Q|=a-1+j+t$, within the ranges
of Lemma~\ref{old:8.4}.
An $A^*$-column has at most
$1+\mathbf1_{\{y\in W_x\}}w_x$ cross-neighbours:
apart from $h$ or $z$, a cross-edge represents a rare pair
outside $D$ and not in $G_x$, so both its ends lie in $W_x$.
A column in $T$ has at most $\gamma m+1$ cross-neighbours.
All column degrees are at most $m/100$, and their squared
sum is at most
\[
 |A^*|+3\sum_{\substack{x\in A^*\\y\in W_x}}w_x^2
       +4\gamma^2m^2t.
\]
Lemma~\ref{old:8.4} gives
\begin{equation}\label{eq:old:20}
 \rho(L_H(y))\le b-j-t-1+\frac C m+
 \frac C{m^2}\sum_{\substack{x\in A^*\\y\in W_x}}w_x^2
 +C\gamma^2t.
\end{equation}
At most four roots are excluded, losing at most $4d$
from $\sum_y d_D(y)=2E$.
Also
\[
 \sum_{y\in B}\sum_{\substack{x\in A^*\\y\in W_x}}w_x^2
 =\sum_{x\in A^*}w_x^3
 \le\omega^2m^2W.
\]
Put $\alpha=s+1+j-r-\theta$.
Since every link radius exceeds $\Lambda$,
summing~\eqref{eq:old:20} and using~\eqref{eq:old:19} yields
\[
 \alpha b+
 \bigl(2-O(\gamma^2)-O(\omega^2/\tau)\bigr)E
 <O(\omega^2/\tau)m(s_++1)+O(d+|s|+1).
\]
The last error includes the excluded roots and their
contribution to $\alpha b$.
Taking $\gamma$ small, then $\eta,\varepsilon$ small and
$m$ large, makes the coefficient of $E$ at least $1.9$
and the right side less than $m(s_++1)/20$.
Thus
\begin{equation}\label{eq:old:21}
 1.9E+(s+1+j-r-\theta)b<\frac1{20}m(s_++1).
\end{equation}

Put $k=1-j-s\in\mathbb Z$ and $\chi_r=2-r-\theta$.
The comparison join in~\eqref{eq:old:18} is $K_q\vee F$,
where
\[
 q=m-k,\qquad |F|=2m+r-1+k,\qquad
 e(F)\le E+6,\qquad \Delta(F)\le d+3.
\]
Equation~\eqref{eq:old:21} becomes
\[
 1.9E+(\chi_r-k)b<\frac1{20}m(s_++1).
\]
Here $\chi_2=0$ and $\chi_1=1-\theta>1/3$.
If $k<0$, then $s_++1\le3(-k)$, contradicting this inequality
and $b\ge1.9m$.
For $r=1$, $\chi_1>1/3$ also excludes $k=0$.
Thus $k\ge0$ for $r=2$ and $k\ge1$ for $r=1$,
and the right side is at most $m/10$.

Set $P_r=\Lambda^2-(q-1)\Lambda-q|F|$.
Direct calculation gives
\[
 P_2=(1+3k)m+k(k+1),\qquad
 P_1=k(\Lambda+m+k)=3km+k(k-\theta),
\]
using $\Lambda^2-(m-1)\Lambda-2m^2=0$ when $r=1$.
Hence $P_r\ge(3k+r-1)m\ge m$.
Since $q\le m$ and $\Delta(F)\le2\varepsilon m+3$,
the prefactor in~\eqref{eq:old:13} is at least $0.9$.
By~\eqref{eq:old:18}, $\rho(K_q\vee F)>\Lambda$, so
\[
 E>0.9P_r-6\ge0.8(3k+r-1)m
\]
for large $m$.
Using $b\le2.1m$ and $\chi_r\ge0$, we obtain
\[
 1.9E+(\chi_r-k)b
 >1.52(r-1)m+2.46km>m,
\]
contradicting the upper bound $m/10$.
\end{proof}

We next allow one unrestricted vertex $h$. We seek a crossing avoiding $h$, since a second disjoint crossing will later pass through this vertex. Throughout the following results put $\Lambda_0=U_m$, $\Lambda_1=V_m$, and $\Lambda_2=2m$.

\begin{lemma}\label{old:8.6}
For sufficiently small positive constants $\eta,\gamma,\tau$ and large $m$, let $V(H)=A\mathbin{\dot\cup}B\mathbin{\dot\cup}\{h\}$, where $a=|A|=m+s$, $b=|B|=2m+r-1-s$, $r\in\{0,1,2\}$, and $|s|\leq\eta m$. Suppose every $A$-vertex belongs to at most $\eta m^2$ missing $AAB$ triples and every $B$-vertex to at most $\eta m^2$ missing $BBB$ triples. If $\sigma(H)>\Lambda_r$, then either $r=0$ and $s\in\{-1,0\}$, or there is a rare triple $xbc$ such that, after orienting $b,c$, $d_B(bc)\geq\gamma m$ and $d_A(xb)\geq(1-\tau)(a-1)$. No assumption is made on $L_H(h)$ or the $AAA$ layer.
\end{lemma}

\begin{proof}
Choose $\eta\ll\omega\ll\gamma\ll\tau\ll1$, with
$\omega^2<\tau/1000$ and $\gamma<10^{-3}$.
For $x\in A$, put
$W_x=\{y\in B:d_A(xy)<(1-\tau)(a-1)\}$,
$w_x=|W_x|$, and $W=\sum_{x\in A}w_x$.
Let
$D=\{uv\in\binom B2:d_B(uv)<\gamma m\}$
and $d=e(D)$.
Counting missing $AAB$ and $BBB$ triples gives
$w_x,\Delta(D)\le\omega m$.
Suppose that the asserted junction does not exist.
Then
$J_x=L_H(x)[B]\subseteq D\cup\binom{W_x}{2}$,
so
$\Delta(J_x)\le2\omega m$ and
$e(J_x)\le d+w_x^2/2$.

In $L_H(x)$, complete the clique
$(A\setminus\{x\})\cup\{h\}$ and all edges incident with $h$.
The resulting graph has a clique of order $a$,
cross-defect $M_x\ge\tau(a-1)w_x$, and radius
\[
 \Lambda_r<\lambda
 \le\rho(K_a\vee I_b)+\Delta(J_x)<2.02m.
\]
Put $F_r(t)=t(t-a+1)-ab$.
Thus $\kappa_x=(\lambda-\Delta(J_x))/(2a)\in[0.9,1.1]$,
and Lemma~\ref{old:8.2} gives
\[
 e(J_x)\ge\kappa_x\bigl(F_r(\lambda)+M_x\bigr).
\]
Set $t_0=2m-2$, $t_1=V_m$, and $t_2=2m$.
We have $\lambda>t_r$, and $F_r$ is increasing on
$[t_r,\infty)$. Direct calculation gives
\[
\begin{aligned}
 F_2(t_2)&=m-(3m+1)s+s^2,\\
 F_1(t_1)&=s^2-(V_m+m)s,\\
 F_0(t_0)&=-3m+2-(3m-3)s+s^2.
\end{aligned}
\]
Let $f=\min\{0.9F_r(t_r),1.1F_r(t_r)\}$.
This choice gives $\kappa_xF_r(\lambda)\ge f$
regardless of the sign of $F_r(t_r)$.
Since $0.9M_x-w_x^2/2\ge(\tau/4)mw_x$, we obtain
$d\ge f+(\tau/4)mw_x$.
Summing over $A$ and dividing by $a$ gives
\begin{equation}\label{eq:old:24}
 d\ge f+\frac{\tau m}{4a}W.
\end{equation}

For $y\in B$, put $Z=N_D(y)$ and $k_y=|Z|$.
In $L_H(y)$, take
$P=(B\setminus(\{y\}\cup Z))\cup\{h\}$ and $Q=A\cup Z$,
of orders $p=b-k_y$ and $q=a+k_y$.
Complete both parts and all edges incident with $h$.
An $A$-column $x$ has at most
$1+w_x\mathbf1_{\{y\in W_x\}}$ cross-neighbours,
and a $Z$-column has at most $\gamma m+1$.
Put $Q_y=\sum_{x:y\in W_x}w_x$.
The squared column sum is at most
$a+3\omega mQ_y+4\gamma^2m^2k_y$.
Here $p\in[1.9m,2.1m]$, $q\in[0.9m,1.1m]$,
and every column degree is at most $m/100$.
Lemma~\ref{old:8.4}, using $2a/m^2\le3/m$, gives
\begin{equation}\label{eq:old:25}
 \rho(L_H(y))\le2m+r-2-s-k_y+\frac3m
 +\frac{6\omega Q_y}{m}+8\gamma^2k_y.
\end{equation}
Set $c=1-8\gamma^2>0.99$ and
$\chi_r=\Lambda_r-(2m+r-2)$.
Sum $\Lambda_r<\rho(L_H(y))$ over $B$, using
$\sum_yk_y=2d$,
$\sum_yQ_y=\sum_xw_x^2\le\omega mW$, and $3b/m<7$,
to obtain
\begin{equation}\label{eq:old:26}
 b(s+\chi_r)+2cd<7+6\omega^2W,
 \qquad \chi_2=0,\quad\chi_1>1/3,\quad\chi_0>0.
\end{equation}

Put $\delta=24a\omega^2/(\tau m)<1/40$.
Eliminating $W$ by~\eqref{eq:old:24} gives
$b(s+\chi_r)+(2c-\delta)d<7-\delta f$.
Since $d\ge\max\{f,0\}$ and $2c>\delta$,
\[
 b(s+\chi_r)+2c\max\{f,0\}
 <7+\delta\max\{-f,0\}.
\]
If $s\ge1$, the displayed polynomials give
$-7ms\le f\le0$, whereas $b(s+\chi_r)\ge1.9ms$.
Thus $1.9ms<7+7\delta ms$, a contradiction for large $m$.

For $s=-k\le0$, the same polynomials give
$f\ge0.9(1+3k)m$ when $r=2$,
$f\ge2.7km$ when $r=1$, and
$f\ge2.7(k-1)m$ when $r=0$ and $k\ge2$.
Using $b\le2.1m$ and $2c>1.98$, the left side exceeds
$m$ in these ranges except $r=1,k=0$.
In particular, for $r=0,k\ge2$, it exceeds
$(3.246k-5.346)m\ge1.146m$.
When $r=1,k=0$, we have $f=0$ and the left side is
$b\chi_1>b/3$, again a contradiction.
Only $r=0,s\in\{-1,0\}$ remain.
\end{proof}

We will repeatedly use the following attachment observation. Suppose the pointwise $AAA$ and $BBB$ defects are at most $\eta m^2$, and the cores $A_0,B_0$ omit at most $\zeta m$ vertices. Call a pair in either core strong if it misses at most $\tau m$ completions in that core. If $xbc$ is rare with $d_A(xb),d_B(bc)\geq\gamma m$, it extends to the crossing
\begin{equation}\label{eq:old:27}
 q,p,x,b,c,z,u.
\end{equation}
Choose $p\in A_0\cap N_H(xb)$ so that $xp$ misses at most $\tau m/2$ completions in $A$; then choose $q\in A_0\cap N_H(xp)$ with $pq$ strong. Choose $z\in B_0\cap N_H(bc)$ so that $cz$ misses at most $\tau m/2$ completions in $B$, and then $u\in B_0\cap N_H(cz)$ with $zu$ strong. At each step, bad pairs exclude $O(\eta m/\tau)$ candidates; core omissions and a bounded prescribed forbidden set disjoint from $\{x,b,c\}$ exclude at most $\zeta m+O(1)$ more. Thus, when $\eta/\tau,\zeta\ll\gamma$, we can make distinct choices and obtain a path with both end pairs in the specified cores.

We next treat the case $r=0$, $s=0$ left open by
Lemma~\ref{old:8.6}.

\begin{lemma}\label{old:8.7}
For sufficiently small $\eta,\gamma>0$ and large $m$, let $V(H)=A\mathbin{\dot\cup}B\mathbin{\dot\cup}\{h\}$ with $|A|=m$ and $|B|=2m-1$. Assume the pointwise $AAB$ and $BBB$ hypotheses of Lemma~\ref{old:8.6}. If $\sigma(H)>U_m$, then there is a rare triple $xbc$ with $d_B(bc),d_A(xb)\geq\gamma m$, after orienting $b,c$.
\end{lemma}

\begin{proof}
Choose $\eta\ll\omega\ll\gamma<10^{-3}$ and put
$W_x=\{y\in B:d_A(xy)<\gamma m\}$,
$w_x=|W_x|$, and $W=\sum_xw_x$.
Let $D=\{uv\in\binom B2:d_B(uv)<\gamma m\}$ and $d=e(D)$.
The pointwise hypotheses give $w_x,\Delta(D)\le\omega m$.
If no asserted junction exists, then
$J_x=L_H(x)[B]\subseteq D\cup\binom{W_x}{2}$, so
$\Delta(J_x)\le2\omega m$ and $e(J_x)\le d+w_x^2/2$.

In $L_H(x)$, complete $(A\setminus\{x\})\cup\{h\}$
to a clique and all edges incident with $h$.
Let $M_x$ be the remaining cross-defect and $\lambda$
the resulting radius. Since $\lambda>U_m>2m-2$ and
$(2m-2)(m-1)-m(2m-1)=-3m+2$,
Lemma~\ref{old:8.2} gives
\[
 (m-1-\gamma m)w_x\le M_x
 \le3m+1.1d+0.55w_x^2.
\]
Using $w_x^2\le\omega m w_x$ and summing gives
$W\le4m+2d$.
The $B$-root argument for~\eqref{eq:old:25}
applies to these weak sets and gives
$2(1-8\gamma^2)d<7+6\omega^2W$.
Hence $d\le4+4\omega^2W$, and therefore
$W\le5m$ and $d\le4+20\omega^2m$.

Let $S$ be the vertices of $B$ outside the support of $D$.
Then $|S|\ge2m-1-2d\ge1.99m$.
Each $y\in S$ belongs to a rare triple $x(y)yz$;
otherwise
$L_H(y)\subseteq K_1\vee(K_m\cup K_{2m-2})$,
contrary to $\rho(L_H(y))>U_m$.
Since $yz\notin D$, we have $y\in W_{x(y)}$.
Thus at least $m-1-\gamma m$ triples of type $AAB$
with $B$-vertex $y$ are missing.
These sets are disjoint for distinct $y$, and each missing
triple contributes to exactly two of the $M_x$. Hence
\[
 \sum_xM_x\ge2|S|(m-1-\gamma m)>3.9m^2.
\]
On the other hand, summing the $A$-root upper bound gives
\[
 \sum_xM_x\le3m^2+1.1md+0.55\omega mW<3.1m^2,
\]
a contradiction.
\end{proof}

We now address the remaining case $r=0$, $s=-1$ from
Lemma~\ref{old:8.6}, obtaining a crossing on at most eight vertices.

\begin{theorem}\label{old:8.8}
For sufficiently small constants $\eta,\zeta,\gamma,\tau>0$ and all sufficiently large $m$, let $V(H)=A\mathbin{\dot\cup}B\mathbin{\dot\cup}\{h\}$, with $|A|=m-1$ and $|B|=2m$. Suppose each $A$-vertex belongs to at most $\eta m^2$ missing $AAA$ triples and at most $\eta m^2$ missing $AAB$ triples, and each $B$-vertex belongs to at most $\eta m^2$ missing $BBB$ triples. Let $A_0\subseteq A,B_0\subseteq B$ omit at most $\zeta m$ vertices in total, and let $\Gamma_A,\Gamma_B$ consist of the pairs in these cores with at most $\tau m$ missing completions in their respective cores. If $\sigma(H)>U_m$, then $H-h$ has a vertex-simple tight path on at most eight vertices, with one end pair in $\Gamma_A$ and the other in $\Gamma_B$. No condition is imposed on the link of $h$.
\end{theorem}

\begin{proof}
Choose $0<\eta,\zeta\ll\omega\ll\gamma\ll\tau\ll1$,
with $\gamma<10^{-3}$, and suppose that no such path exists.
Put $W_x=\{y\in B:d_A(xy)<\gamma m\}$,
$w_x=|W_x|$, and $W=\sum_xw_x$.
Let $D=\{uv\in\binom B2:d_B(uv)<\gamma m\}$,
$d=e(D)$, and $k_y=d_D(y)$.
The pointwise hypotheses give $w_x,\Delta(D)\le\omega m$.

Let $J_x=L_H(x)[B]$.
The attachment observation~\eqref{eq:old:27} gives
$J_x\subseteq D\cup\binom{W_x}{2}$, since otherwise
there is a seven-vertex core crossing.
Thus $\Delta(J_x)\le2\omega m$ and
$e(J_x)\le d+w_x^2/2$.
In $L_H(x)$, complete the $(m-1)$-vertex clique
$(A\setminus\{x\})\cup\{h\}$ and all edges incident with $h$.
Its cross-defect satisfies $M_x\ge(m-2-\gamma m)w_x$.
Since
$t(t-m+2)-2m(m-1)=(t-2m+2)(t+m)$,
Lemma~\ref{old:8.2}, applied at the resulting radius
$\lambda>U_m>2m-2$, gives
\[
 M_x<1.02e(J_x)\le1.02(d+w_x^2/2).
\]
In particular, $e(J_x)>0$.
Using $w_x^2\le\omega m w_x$ and absorbing gives
\[
 w_x<1.1d/m\quad\text{for every }x\in A,
 \qquad W<1.1d.
\]

For $y\in B$, put $Z=N_D(y)$ and partition its link into
$P=(B\setminus(\{y\}\cup Z))\cup\{h\}$ and $Q=A\cup Z$,
of orders $2m-k_y$ and $m-1+k_y$.
After completing both parts and all edges incident with $h$,
an $A$-column $x$ has sum at most
$1+w_x\mathbf1_{\{y\in W_x\}}$,
and a $Z$-column has sum at most $\gamma m+1$.
Put $Q_y=\sum_{x:y\in W_x}w_x$.
The squared column sum is at most
$m-1+3\omega mQ_y+4\gamma^2m^2k_y$.
The part sizes and column degrees satisfy the hypotheses
of Lemma~\ref{old:8.4}, which gives, with $c=1-8\gamma^2$,
\[
 \rho(L_H(y))\le2m-1-ck_y+3/m+6\omega Q_y/m.
\]
Summing over $B$, using $\rho(L_H(y))>2m-2$,
$\sum_yk_y=2d$, and
$\sum_yQ_y=\sum_xw_x^2\le\omega mW$, gives
$(2c-6.6\omega^2)d<2m+6$.
Hence $d<1.1m$ and $w_x<1.21<2$.
By integrality, $w_x\le1$, so $J_x\subseteq D$.
Consequently each $A$-column in the same partition
has only $h$ as a cross-neighbour.
Repeating the estimate gives
\[
 \rho(L_H(y))\le2m-1-ck_y+3/m.
\]
If $k_y\ge2$, the right side is less than $2m-2$,
a contradiction. Thus $D$ is a matching.

Choose a rare triple $xuv$, which exists by the $A$-root
bound above. We have $uv\in D$ and may orient it so that
$u\notin W_x$, since $w_x\le1$.
If $uvz\in H[B]$, then $vz\notin D$, so
$d_A(xu),d_B(vz)\ge\gamma m$.
The attachment choices from~\eqref{eq:old:27},
applied at $xu$ and $vz$, extend $x,u,v,z$ to
$q,p,x,u,v,z,t,w$, with $pq\in\Gamma_A$ and $tw\in\Gamma_B$.
At each step at most $O(\eta m/\tau)+\zeta m+O(1)$
candidates are excluded, so all choices can be made
distinct. This contradicts our assumption.
Thus $uv$ has no $BBB$ completion.

In $L_H(u)$, no edge joins $A$ to $B\setminus\{u,v\}$,
since $J_x\subseteq D$ for every $x\in A$ and $D$ is a matching.
Nor does $v$ have a neighbour in $B\setminus\{u,v\}$,
since $d_B(uv)=0$.
Complete the two blocks $A\cup\{v\}$ and $B\setminus\{u,v\}$,
of orders $m$ and $2m-2$, and all edges incident with $h$.
Then $L_H(u)\subseteq K_1\vee(K_m\cup K_{2m-2})$,
contrary to $\rho(L_H(u))>U_m$.
\end{proof}

\begin{corollary}\label{old:8.9}
For sufficiently small $\eta,\varepsilon,\zeta,\tau>0$ and large $m$, let $V(H)=A\mathbin{\dot\cup}B\mathbin{\dot\cup}\{h\}$, where $|A|=m+s$, $|B|=2m+r-1-s$, $r\in\{0,1,2\}$, and $|s|\leq\eta m$. Suppose every $A$-vertex belongs to at most $\varepsilon m^2$ missing triples of each type $AAA,AAB$, and every $B$-vertex to at most $\varepsilon m^2$ missing $BBB$ triples. Let $A_0\subseteq A,B_0\subseteq B$ omit at most $\zeta m$ vertices in total, and let $\Gamma_A,\Gamma_B$ consist of their pairs missing at most $\tau m$ completions in the respective cores. If $\sigma(H)>\Lambda_r$, then $H-h$ has a crossing on at most eight vertices with ends in $\Gamma_A,\Gamma_B$.
\end{corollary}

\begin{proof}
Choose $\gamma\ll\tau$, then take
$\alpha=\max\{\eta,\varepsilon\}$ and $\zeta$ sufficiently
small for the preceding results and~\eqref{eq:old:27}.
Lemma~\ref{old:8.6}, with allowance $\alpha$, and~\eqref{eq:old:27}
give a seven-vertex crossing unless $r=0,s\in\{-1,0\}$.
For $r=0,s=0$, use Lemma~\ref{old:8.7} and~\eqref{eq:old:27};
for $r=0,s=-1$, use Theorem~\ref{old:8.8}.
\end{proof}

\subsection{The divisible orders and completion of the proof}

\begin{lemma}\label{old:8.11}
Fix $K\ge1$. For sufficiently small
$\eta,\gamma,\tau,\varepsilon>0$ and large $m$, let
$|A|=m+s$, $|B|=2m-s$, where $|s|\le\eta m$,
and assume the pointwise role hypotheses of
Theorem~\ref{old:8.5}.
Define $D,W_x,G_x$ as in that theorem.
Suppose there are sets $F\subseteq A$ and $T\subseteq B$,
both of size at most $K$, such that $e(G_x)\le K$
for $x\notin F$, $e(G_x)\le K|B|$ for $x\in F$,
and $\sum_xd_{G_x}(y)\le K$ for $y\notin T$.
Then $\sigma(H)\le U_m$.
\end{lemma}

\begin{proof}
Put $a=|A|$, $b=|B|$, $\lambda_0=2m-2$,
$E=e(D)$, $W=\sum_xw_x$, and $s_+=\max\{s,0\}$.
As before, $w_x\le\omega m$ and
$\Delta(D)\le2\varepsilon m$, where $\omega=C\varepsilon/\tau$.
Choose $\gamma,\tau$ first, and then $\eta,\varepsilon$
sufficiently small in terms of them and $K$, with $\omega\ll\tau$.
Suppose $\sigma(H)>U_m$.
Let $J_x=L_H(x)[B]$, and let $M_x$ count missing cross-pairs
between $A\setminus\{x\}$ and $B$. Then
\begin{equation}\label{eq:old:37}
 J_x\subseteq G_x\cup D\cup\binom{W_x}{2},
 \qquad M_x\ge\frac{\tau a}{2}w_x.
\end{equation}

The polynomial of $K_{a-1}\vee I_b$ is
$P_s(z)=z^2-(m+s-2)z-(m+s-1)(2m-s)$.
Since $s$ is an integer, the identities
$P_s(\lambda_0)=s(s+1-3m)$ and
$P_s(\lambda_0+s)=s(s-1)$ show that its radius is at most
$\lambda_0+s_+$.
The graph argument in Lemma~\ref{old:5.3}, with clique
order $a-1$ and spectral allowance $s_+$, gives
\begin{equation}\label{eq:old:38}
 M_x\le C\bigl(ms_++e(J_x)\bigr).
\end{equation}
This includes $x\in F$, since that argument imposes
no maximum-degree condition on the added edges.
Using $\sum_xe(G_x)\le Ka+K^2b=O_K(m)$ and
$\sum_xw_x^2\le\omega mW$, sum
\eqref{eq:old:37}--\eqref{eq:old:38} and absorb the quadratic term:
\begin{equation}\label{eq:old:39}
 W\le\frac{C_K}{\tau}\bigl(ms_++E+1\bigr).
\end{equation}

For $y\in B\setminus T$, put $Z=N_D(y)$, $t=|Z|$,
$P=B\setminus(\{y\}\cup Z)$ and $Q=A\cup Z$.
An $A$-column of the cross matrix of $L_H(y)$ has at most
$d_{G_x}(y)+\mathbf1_{\{y\in W_x\}}w_x$ neighbours in $P$;
a $Z$-column has fewer than $\gamma m$.
Since $\sum_xd_{G_x}(y)\le K$, their squared degrees sum to at most
$2K^2+2\sum_{x:y\in W_x}w_x^2+\gamma^2m^2t$.
The part sizes and column degrees satisfy Lemma~\ref{old:8.4},
which gives
\begin{equation}\label{eq:old:40}
 \rho(L_H(y))\le\lambda_0-s-t+\frac{C_K}{m^2}
 +\frac C{m^2}\sum_{x:y\in W_x}w_x^2+C\gamma^2t.
\end{equation}
Choose $c=1-C\gamma^2\ge0.99$.
Every edge of $D$ not contained in $T$ contributes at least
once to $\sum_{y\notin T}d_D(y)$, so this sum is at least
$E-\binom K2$.
Using $\sum_xw_x^3\le\omega^2m^2W$ and summing over $B\setminus T$
gives
\begin{equation}\label{eq:old:41}
 s(b-|T|)+c\left(E-\binom K2\right)
 <\frac{C_K}{m}+C\omega^2W.
\end{equation}
Taking $\varepsilon$ sufficiently small, substitution
of~\eqref{eq:old:39} yields
\[
 0.98E+s(b-|T|)<\frac{ms_+}{100}+O_K(1).
\]
This excludes $s\ge1$.
For $s=0$, we have $E,W=O_{K,\tau}(1)$.
Since $|T|+2E+W=O_{K,\tau}(1)$, choose $y$ outside $T$,
all endpoints of $D$, and all sets $W_x$.
Then~\eqref{eq:old:40} gives
\[
 \rho(L_H(y))\le\lambda_0+\frac{C_K}{m^2}<U_m,
\]
where the last inequality uses
$U_m>\lambda_0+1/(3m)$ from Lemma~\ref{old:2.1}.
This is a contradiction.

Finally, let $s=-k\le-1$.
The same inequality gives $E\le2.1km$ for small $\eta$
and large $m$.
Choose $x\in A\setminus F$ and put $R=D\cup G_x$.
Then $e(R)\le E+K$ and $\Delta(R)\le2\varepsilon m+K$.
Every weak column misses at least
$\tau a-1\ge w_x-1$ clique neighbours, and all added
edges outside $R$ lie in $\binom{W_x}{2}$.
Also
\[
 \rho(K_{m-k-1}\vee R)
 \ge\sqrt{(m-k-1)(2m+k)}>1.3m,
 \qquad m-k-1+\Delta(R)<1.1m.
\]
Lemma~\ref{old:8.3} therefore gives
$\rho(L_H(x))\le\rho(K_{m-k-1}\vee R)$.
At $\lambda_0$,
$P_{-k}(\lambda_0)=k(3m+k-1)\ge3km$.
Since $(\lambda_0-\Delta(R))/(2(m-k-1))\ge0.9$,
equation~\eqref{eq:old:13} gives
$E+K\ge e(R)>2.7km$, contradicting $E\le2.1km$
for sufficiently large $m$.
\end{proof}

\begin{lemma}\label{old:8.12}
Fix a sufficiently small endpoint allowance $\theta>0$.
For sufficiently small $\eta,\varepsilon>0$ and large $m$,
let $|A|=m+s$, $|B|=2m-s$, $|s|\le\eta m$,
and assume the pointwise role hypotheses of
Theorem~\ref{old:8.5}.
Suppose $\sigma(H)>U_m$ and $H$ is not Hamiltonian.
There is a vertex $h$ such that, on putting
$A'=A\setminus\{h\}$ and $B'=B\setminus\{h\}$,
every set $X\subseteq V(H)\setminus\{h\}$ of size at most
$20$ can be avoided by a crossing through $h$ on at most
seven vertices, whose end pairs lie in $A',B'$ and miss
at most $\theta m$ completions in their respective parts.
Consequently, $H-h$ has no crossing on at most $20$
vertices with such end pairs.
Moreover, $|A'|=m+t$, $|B'|=2m-1-t$ for
$t\in\{s-1,s\}$.
\end{lemma}

\begin{proof}
Choose the usability and defect parameters sufficiently small
after $\theta$, in particular for Lemma~\ref{old:8.11} with $K=100$.
Disjoint usable pairs of distinct colours extend successively
by~\eqref{eq:old:27} to disjoint crossings, which would give
a Hamilton cycle by Theorem~\ref{old:7.2}.
Thus the families $G_x$ are pairwise cross-intersecting.
Lemma~\ref{old:8.1} gives the finite-support or star
alternative, with the triangle case included in the former.

In the finite-support case, choose $h\in A$ and
$S\subseteq B$, $|S|\le4$, supporting every $G_x$ with $x\ne h$.
If $\nu(G_h)\le40$, the endpoints of a maximum matching
form a vertex cover $U$ of size at most $80$.
With $F=\{h\}$ and $T=U\cup S$, we have
$e(G_x)\le6$ for $x\ne h$, $e(G_h)\le80|B|$,
and $\sum_xd_{G_x}(y)\le80$ for $y\notin T$.
This contradicts Lemma~\ref{old:8.11} with $K=100$.
Thus $\nu(G_h)\ge41$.

In the star case, let $h\in B$ be the centre.
Form the bipartite graph on $A,B\setminus\{h\}$ with
edge $xy$ whenever $xhy$ is usable.
A maximum matching of size at most $40$ would have
an endpoint cover $F\cup T_0$ of size at most $80$.
Taking $T=T_0\cup\{h\}$, each colour outside $F$ has
at most $80$ usable edges, each colour in $F$ has
at most $|B|$, and every $y\notin T$ has total usable
incidence at most $80$.
Lemma~\ref{old:8.11} again gives a contradiction.

In either case, we have $41$ usable triples whose pairwise
intersections equal $\{h\}$.
Given $X\subseteq V(H)\setminus\{h\}$ with $|X|\le20$,
one of these triples avoids $X$.
Write it as $xuv$, oriented with $u\notin W_x$.
The cores $A',B'$ omit only $h$.
Apply~\eqref{eq:old:27} with these cores, endpoint allowance
$\theta/2$, and forbidden set $X$.
It gives a crossing $q,p,x,u,v,z,w$ through $h$,
with end pairs in $A',B'$ missing at most $\theta m/2$
completions there.

If a crossing $Q\subseteq H-h$ on at most $20$ vertices
had the stated end-pair strength, apply this construction
with $X=V(Q)$.
The two disjoint crossings use at most $27$ vertices.
The inherited pointwise bounds and $\sigma(H)>U_m>0.65(3m)$
allow Theorem~\ref{old:7.2} with cores $A',B'$ and
exceptional set $\{h\}$, giving a Hamilton cycle, a contradiction.
The size formula holds with $t=s-1$ for $h\in A$
and $t=s$ for $h\in B$.
\end{proof}

\begin{theorem}\label{old:8.13}
For suitable $\eta,\varepsilon>0$ and large $m$, suppose $|A|=m+s$, $|B|=2m-s$, $|s|\leq\eta m$, with at most $\varepsilon m^2$ missing triples of each type $AAA,AAB$ through every $A$ vertex and at most $\varepsilon m^2$ missing $BBB$ triples through every $B$ vertex. If $\sigma(H)>U_m$, then $H$ is Hamiltonian.
\end{theorem}
\begin{proof}
Fix an endpoint allowance $\theta>0$ sufficiently small for
Lemma~\ref{old:8.12} and Corollary~\ref{old:8.9}.
Choose $\eta,\varepsilon$ sufficiently small for both results,
with the corollary allowing size error $2\eta m$.

Suppose $H$ is not Hamiltonian.
Lemma~\ref{old:8.12} gives $h,A',B'$ with
$|A'|=m+t$ and $|B'|=2m-1-t$, where
$|t|\le\eta m+1\le2\eta m$ for large $m$.
It prohibits every crossing in $H-h$ on at most $20$ vertices
whose end pairs miss at most $\theta m$ completions
in their respective parts.
The pointwise defect bounds pass to $A',B'$.
Apply Corollary~\ref{old:8.9} to the original $H$ with
partition $A'\mathbin{\dot\cup}B'\mathbin{\dot\cup}\{h\}$,
$r=0$, and cores $A',B'$ themselves.
It gives a crossing in $H-h$ on at most eight vertices
with precisely this end-pair strength, a contradiction.
\end{proof}

\begin{theorem}\label{old:8.14}
There are $\delta>0$ and $n_0$ such that the following holds for $n=3m+r\geq n_0$, where $r\in\{0,1,2\}$. Suppose $H$ has a partition $A\mathbin{\dot\cup}B$ with $\bigl||A|-n/3\bigr|\leq\delta n$ and at most $\delta n^3$ missing triples of types $AAA,AAB,BBB$. If $\sigma(H)>\Lambda_r$, where $\Lambda_0=U_m$, $\Lambda_1=V_m$ and $\Lambda_2=2m$, then $H$ has a tight Hamilton cycle. There is no restriction on exceptional vertices, individual pair codegrees, or the number of $ABB$ triples.
\end{theorem}
\begin{proof}
Fix an endpoint allowance $\theta>0$ sufficiently small for
Lemma~\ref{old:7.4} and Corollary~\ref{old:8.9}.
Choose admissible size, pointwise-defect and core-omission
allowances $\eta_*,\varepsilon_*,\zeta_*$ for the corollary,
with $\eta_*,\varepsilon_*$ also admissible in
Theorems~\ref{old:8.5} and~\ref{old:8.13}.
Choose $\xi>0$ sufficiently small for Lemma~\ref{old:7.4}, with
\[
 \xi\le\min\{1/100,\eta_*/5,\varepsilon_*/16,\zeta_*/4\},
\]
and let $\zeta_0=\zeta(\theta,\xi)$ be supplied by that lemma.
Take $\delta\le\zeta_0$ and $n_0$ sufficiently large for all
cited results, with $m\ge1/\xi$ and $n\ge2/\zeta_0$.
Since
$\Lambda_r>2n/3-2\ge(2/3-\zeta_0)n$,
the spectral hypothesis of Lemma~\ref{old:7.4} holds.

Suppose $H$ is not Hamiltonian.
That lemma gives a partition
$A'\mathbin{\dot\cup}B'\mathbin{\dot\cup}Z$, where $|Z|\le1$,
and cores $A_0,B_0$.
Put $s'=|A'|-m$.
As $n\le4m$, the pointwise defects, size error and core omissions
satisfy, respectively,
\[
\begin{aligned}
 \xi n^2&\le16\xi m^2\le\varepsilon_*m^2,\\
 |s'|&\le\xi n+r/3\le5\xi m\le\eta_*m,\\
 \xi n&\le4\xi m\le\zeta_*m.
\end{aligned}
\]
If $Z=\varnothing$, then $|B'|=2m+r-s'$.
Theorem~\ref{old:8.5} applies for $r=1,2$, and
Theorem~\ref{old:8.13} for $r=0$, a contradiction.
If $Z=\{h\}$, then $|B'|=2m+r-1-s'$.
Apply Corollary~\ref{old:8.9} to the original $H$ with
parts $A',B'$, unrestricted vertex $h$, and cores $A_0,B_0$.
It gives a crossing in $H-h$ on at most eight vertices,
whose end pairs miss at most $\theta m\le\theta n$
completions in these cores.
This contradicts the conclusion of Lemma~\ref{old:7.4}.
\end{proof}

\begin{proof}[Proof of Theorem~\ref{thm:tight}]
Proposition~\ref{prop:sharp} establishes sharpness.
Suppose there is a sequence of non-Hamiltonian $H_n$ with
$\sigma(H_n)>T_n$ and $n\to\infty$.
Since $T_n=2n/3+O(1)$, Corollary~\ref{old:4.10} applies.
Passing to a subsequence, fix $i\in\{1,2,3\}$ such that
$b_i(H_n)=o(1)$, and choose partitions attaining this density.

For $i=1$, these partitions satisfy the space hypothesis
of Theorem~\ref{old:5.7} for all sufficiently large $n$.
For $i=2$, Lemma~\ref{old:4.11}(2) supplies the missing-layer
bound required by Theorem~\ref{old:6.1}.
For $i=3$, Lemma~\ref{old:4.11}(3) supplies the normal-layer
bounds required by Theorem~\ref{old:8.14}.
The uniform forms of these estimates verify the fixed defect
allowances in the respective theorems.

Their spectral thresholds are at most $T_n$: the space case
uses $\widehat T_n\le T_n$, and the other two use $U_m,V_m,2m$
for $n=3m,3m+1,3m+2$, respectively. Here $T_{3m}\ge U_m$
by Lemma~\ref{old:2.1}. Each case therefore gives a tight
Hamilton cycle for all sufficiently large members of the
sequence, a contradiction.
\end{proof}

\begin{remark}\label{old:9.1}
The large-order qualification in Theorem~\ref{thm:tight} is necessary. For $n=9$, take $|A|=2$, $|B|=7$, include all triples meeting $A$, and put a Fano plane on $B$. The two $A$ vertices of a tight Hamilton cycle leave at most two $B$-runs, one of which contains four consecutive vertices. Their two consecutive triples share a pair, contrary to the linearity of the Fano plane. Every $B$-link is $K_2\vee3K_2$, and the $A$-links are complete. Hence $\sigma(H)=1+2\sqrt3>W_3=(3+\sqrt{33})/2$. Determining all exceptional small orders and classifying the hypergraphs attaining the large-order thresholds remain separate questions.
\end{remark}

\section*{Acknowledgments}
The author is supported by the China Scholarship Council. The author sincerely thanks Professor Suil O for his helpful comments and suggestions.

\section*{Declaration of generative AI and AI-assisted technologies in the manuscript preparation process}

During the preparation of this work, the author used ChatGPT (OpenAI) for language refinement, technical editing, and computational verification of examples. The author reviewed and edited the output as needed and takes full responsibility for the content of the published article.

\section*{Data availability}
Data sharing is not applicable to this article as no datasets were
generated or analysed during the current study.

\end{document}